\documentclass[12pt]{extreport}

\usepackage{tikz}
\usepackage{subfigure}
\usepackage{hyperref}
\usepackage[text={450pt,600pt}]{geometry}
\usepackage{titlesec}
\usepackage{amssymb}
\usepackage{xcolor}
\usepackage{amsmath,amssymb,amscd,mathrsfs}
\usepackage{color}
\usepackage{amssymb}
\usepackage{amsmath}
\usepackage{amscd}
\usepackage{ntheorem}
\usepackage{xcolor}

\DeclareMathOperator{\Prog}{Prog}

 \newcommand{\sgn}{\operatorname{sgn}}
 \newcommand{\esssup}{\operatorname{ess\ sup}}
\usepackage{cleveref}

\hypersetup{
  colorlinks   = true,    
  urlcolor     = blue,    
  linkcolor    =  blue,  
  citecolor    = blue      
}
\newcommand*{\qed}{\par\vspace{5mm} \hfill\ensuremath{\Box}}
\usepackage[latin1]{inputenc}   
\usepackage[T1]{fontenc}   
\usepackage{mathtools}
\DeclarePairedDelimiterX{\norm}[1]{\lVert}{\rVert}{#1}

\begin{document}


\newtheorem{thm}{Theorem} [chapter]
  \newtheorem{Conjecture}[thm]{Conjecture}
   \newtheorem{Lemma}[thm]{Lemma} 
\newtheorem{prop}[thm]{Proposition}
 \newtheorem{Corollary}[thm]{Corollary}
     \newtheorem{ex}[thm]{Example}
    \theorembodyfont{\normalfont}
  \newtheorem{rem}{Remark} 
     \newtheorem{as}{Assumption}[chapter]
\renewcommand{\theas}{{{(A.\arabic{as}})}}
\newtheorem{defi}{Definition}[chapter]
  \renewcommand{\thedefi}{\arabic{chapter}.\arabic{defi}.}
  \renewcommand{\therem}{\arabic{rem}.}
   \theorembodyfont{\normalfont}\theoremstyle{nonumberplain}  
     \newtheorem{Proof}{Proof.}
  \newtheorem{remm}{Remark.} 

  \renewcommand{\thethm}{\arabic{chapter}.\arabic{thm}}%

\newenvironment{thmbis}[1]
  {\renewcommand{\thethm}{\ref{#1}$'$}%
   \addtocounter{thm}{-1}%
   \begin{thm}}
  {\end{thm}}

\newenvironment{asp}[1]
 {\addtocounter{as}{-1}%
  \renewcommand{\theas}{(A.\arabic{as}$'$)}%
   \begin{as}}
  {\end{as}}

\newenvironment{aspp}[1]
 {\addtocounter{as}{-1}%
  \renewcommand{\theas}{(A.\arabic{as}$''$)}%
   \begin{as}}
  {\end{as}}

\makeatletter
\newcommand{\neutralize}[1]{\expandafter\let\csname c@#1\endcsname\count@}
\makeatother

\begin{titlepage}
\title{
\vspace{-2cm}
\Huge{\textbf{On Quadratic BSDEs with Final Condition in $\mathbb{L}^2$}} 
}
\author{\Large{\bf Hanlin Yang} 
\\
\\
\Large{Master Thesis}\\
\\
\Large{University of Z\"urich and ETH Z\"urich} 
\\
\\
\includegraphics[width=8cm]{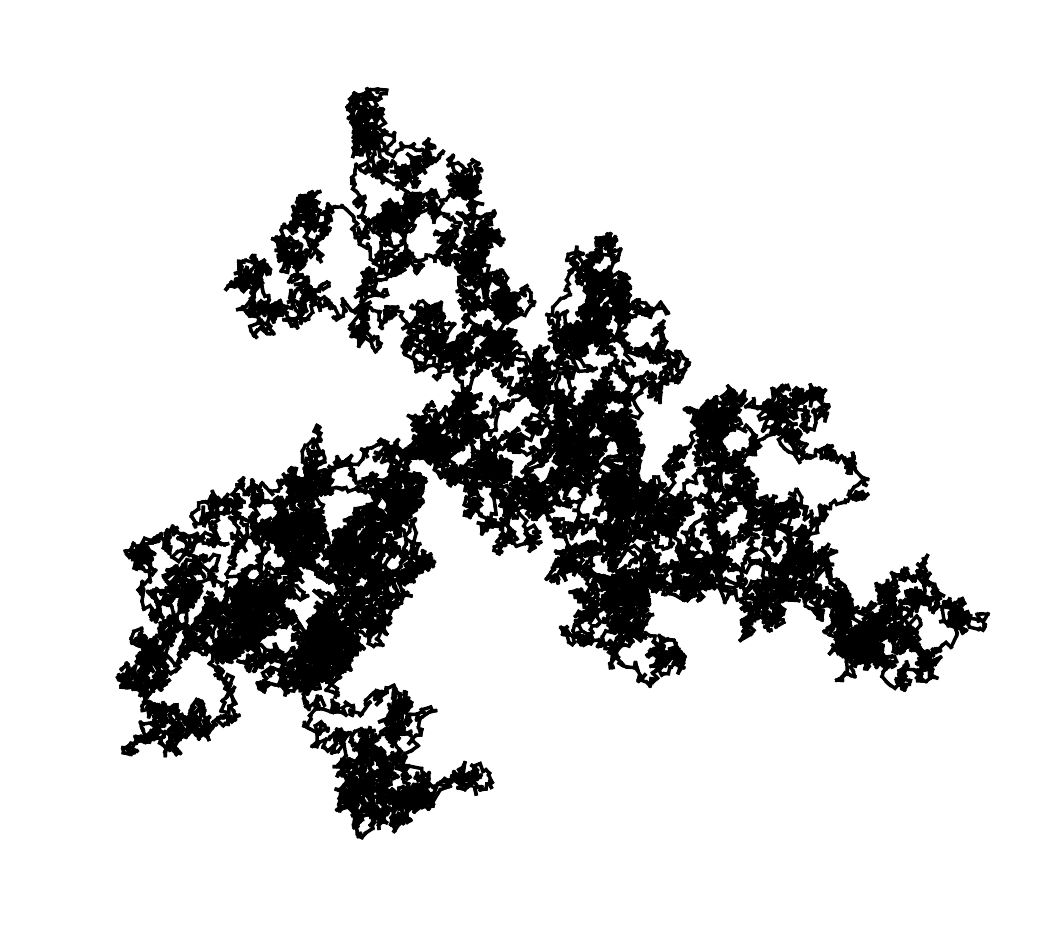}
\\
\Large{Supervised by}
\\
\\
\Large{\bf Professor Martin Schweizer}
\\
\\
\Large{ETH Z\"urich}
}
\date{\Large{Submitted in April, 2015\\
 Revised in May, 2015
}\textcopyright} 
\maketitle
\chapter*{Acknowledgements}
First of all, I would like to express my gratitude to my supervisor Professor Martin Schweizer, for being generous with his time in advising my master thesis. 
Professor Martin Schweizer gave me a  fascinating thesis topic from which I understand  probability theory and mathematical finance better. I am also grateful  for   numerous insightful comments which  
 enrich the contents of my master  thesis. 
Finally, I am grateful to him for supporting my Ph.D study.

Special thanks belong to Meng Chen and Xinyi Li, for carefully reading this manuscript and providing valuable comments. 

I would also like to thank Professor Markus Leippold for funding my Ph.D research.

Finally, to my parents, I am thankful for their love and support.

\tableofcontents
\end{titlepage}
\chapter{Introduction}
We are concerned with $\mathbb{R}$-valued backward stochastic differential equations (BSDEs) in the continuous semimartingale framework
\begin{align}
Y_t  = \xi +\int_t^T \big(F(s, Y_s, Z_s)dA_s +  g_sd\langle N\rangle_s\big)-\int_t^T \big(Z_s dM_s + dN_s\big), \label{c1:bsde}
\end{align}
where the generator $F(t, y, z)$ has at most  quadratic growth in $z$ and $g$ is a progressively measurable integrable process. For this reason, 
 (\ref{c1:bsde}) is called \emph{quadratic}.
 We call an adapted process $(Y, Z, N)$ a \emph{solution} to (\ref{c1:bsde}) if $Y$ is continuous, $Z$ is progressively measurable and integrable with respect to the fixed continuous local martingale $M$, and $N$ is a continuous local martingale strongly orthogonal to $M$. In particular, if the filtration is generated by a Brownian motion $W$, (\ref{c1:bsde})  becomes the classic BSDE with $A_t=t$, $M=W$ and $N_\cdot =0$.

Let us recall that, 
quadratic BSDEs
 are first studied by Kobylanski \cite{K2000}. Existence and uniqueness, comparison theorem and stability results are proved, when the terminal value is bounded. 
  Later, Briand and Hu \cite{BH2006}, \cite{BH2008}
extend the  existence result  by assuming that the terminal value has  exponential moments integrability. Moreover, a uniqueness result is obtained given a convexity condition as an additional requirement. Afterwards, Morlais \cite{M2009} and Mocha and Westray \cite{MW2012} extend all these results to continuous semimartingale setting under rather strong assumptions on the generator. 
 Recently, for Brownian framework, Bahlali et al \cite{BEO2014} constructs a solution to quadratic BSDEs with  the terminal value in $\mathbb{L}^2$ and the  generator 
$F(t, y, z)$  satisfying
 $\mathbb{P}$-a.s. for all $(t, y, z)\in [0, T]\times\mathbb{R}\times\mathbb{R}^d$, 
\begin{align}
|F(t, y, z)|\leq \alpha + \beta|y| + \gamma |z| + f(|y|)|z|^2,  \label{c1:bsde2}
\end{align}
for some $\alpha, \beta, \gamma \geq 0$ and $f(|\cdot|):\mathbb{R}\rightarrow\mathbb{R}$ which is integrable and bounded on compact subsets of $\mathbb{R}$.
However, as to the uniqueness of a solution, only purely quadratic BSDEs are studied.

As a natural extension of these works, this paper is devoted to answering the following questions:  {\bf 1.}
Does existence and uniqueness hold for  BSDEs satisfying (\ref{c1:bsde2}) with  terminal value in $\mathbb{L}^p$ for a cetain $p>0$ ? 
 {\bf 2.} Can one establish the solvability of quadratic semimartingale BSDEs in a more general way under weaker assumptions ?

In {\bf Chapter \ref{lp}} we address the first question. 
We prove an existence result, by merely assuming that the generator is monotonic at $y=0$ and has a linear-quadratic growth in $z$ of type (\ref{c1:bsde2}), and that the terminal value belongs to $\mathbb{L}^p$ for a certain $p>1$. To establish the a priori estimates, we use a combination of the estimates developed by Bahlali et al \cite{BEO2014}  and $\mathbb{L}^p$-type estimates developed by Briand et al \cite{B2003}.  Thanks to the estimates,  we prove an existence result  based on  the localization procedure developed by Briand and Hu \cite{BH2006}, \cite{BH2008}.
The second contribution  of this chapter is the uniqueness result. In the spirit of Da Lio and Ley \cite{LL2006} or  Briand and Hu \cite{B2008}, we  prove comparison theorem, uniqueness and a stability result via $\theta$-technique under a convexity assumption.  It turns out that our results of existence and uniqueness not simply provide wider perspectives on quadratic BSDEs but also, by setting $f(|\cdot|)=0$,  concern non-quadratic BSDEs studied in \cite{P1999}, \cite{BC2000}, \cite{B2003}, \cite{BEK2013}, etc. 

{\bf Chapter \ref{lp}} is organized as follows.  In Section \ref{lpp}, \ref{I}, 
we introduce  basic notions   and present auxiliary results. Section \ref{itokrylov}
 proves  the It\^{o}-Krylov formula  and a generalized It\^{o} formula for $y\mapsto |y|^p (p\geq 1)$. The former one is used to treat discontinuous quadratic generators or discontiuous quadratic growth and the later one is used for a $\mathbb{L}^p$-type estimate. 
  Section \ref{degenerate} reviews  purely quadratic BSDEs and their natural extensions, based on
 Bahlali et al \cite{BEK2013}.
 Section \ref{lpsolution}  studies existence, uniqueness and a stability result. 
  Finally, in Section \ref{qpde},  we derive the probabilistic representation for the viscosity solution to the associated quadratic PDEs.

{\bf Chapter \ref{qsbsde}} 
addresses the second question by  
using a  regularization procedure which is different from   Morlais \cite{M2009} and Mocha and Westray \cite{MW2012}. The first contribution is to obtain an existence and uniqueness result given a Lipschitz-continuous generator and a bounded integrand $g$. BSDEs of this type are called \emph{Lipschitz-quadratic}, and serve as an basic ingredient for the study of  quadratic BSDEs.  In the second step, we prove  a more general version of monotone stability result which allows one to construct solutions to  quadratic BSDEs via
Lipschitz-quadratic regularizations. Finally, we rely on a convexity assumption to obtain the uniqueness result via $\theta$-technique.
 
{\bf Chapter \ref{qsbsde}} is organized as follows. Section \ref{section20} presents the basic notions of semimartingale BSDEs.  Section \ref{section21} concerns the existence and uniqueness result for Lipschitz-quadratic BSDEs. In Section \ref{section22}, we prove a general version of monotone stability. As an application,   the existence of a bounded solution is immediate. In Section \ref{section23}, existence and uniqueness of unbounded solution are proved. Finally, we show in Section that the martingale part of a solution defines an equivalent change of measure.

{\bf Chapter \ref{qs}} is a survey of the stability result of  quadratic semimartingales studied in  Barrieu and El Karoui \cite{BEK2013}. Section \ref{qs2}
introduces the notion of quadratic semimartingales and their characterizations. 
In Section \ref{qs3}, 
we use a forward point of view to address the issue of convergence:  the stability of quadratic semimartingales is  proved in the first step; it is then used to deduce the convergence of the martingale parts. 
Finally, in Section \ref{section24},  the solutions to quadratic  BSDEs are characterized as quadratic semimartingales. As a counterpart, a corresponding monotone stability result for BSDEs are formulated.   The prime advantage of this stability result, in contrast to others, is that  the boundedness is no longer needed.  

\chapter{ $\mathbb{L}^p (p \geq 1)$ Solutions to Quadratic BSDEs}
\label{lp}
\section{Preliminaries}
\label{lpp}
In this chapter, we   study  a class of quadratic BSDEs driven by Brownian motion.  We fix  the time horizon $T>0$ and   a $d$-dimensional Brownian motion  $(W_t)_{0\leq t\leq T}$ defined on a complete probability space $(\Omega, \mathcal{F},\mathbb{P})$.  $(\mathcal{F}_t)_{0\leq t \leq T}$ is the  filtration generated by $W$  and augmented by $\mathbb{P}$-null sets of $\mathcal{F}$. Any measurability will refer to this filtration. In particular, $\Prog$  denotes the progressive $\sigma$-algebra  on $\Omega \times [0, T]$.  Let us  introduce the notion of BSDEs and their solutions in the following paragraph. As mentioned in the introduction, we  exclusively study $\mathbb{R}$-valued BSDEs.

{\bf BSDEs: Definition and Solutions.}
Let $\xi$ be an $\mathbb{R}$-valued $\mathcal{F}_T$-measurable random variable, $F: \Omega\times[0, T]\times \mathbb{R}\times\mathbb{R}^d \rightarrow \mathbb{R}$  a $\Prog\otimes \mathcal{B}(\mathbb{R})\otimes \mathcal{B}(\mathbb{R}^d)$-measurable random function. The BSDEs of our study can be written as 
\begin{align}
Y_t = \xi +\int_t^T F(s, Y_s, Z_s)ds -\int_t^T Z_s dW_s,\label{bsde}
\end{align}
where $\int_0^\cdot Z_sdW_s$, sometimes denoted by $Z\cdot W$,  refers to the vector stochastic integral; see, e.g.,
Shiryaev and  Cherny \cite{SC2002}. 
We call a process  $(Y, Z)$ valued in $\mathbb{R}\times\mathbb{R}^d$
a \emph{solution} to   (\ref{bsde}), if  $Y$ is a continuous adapted process and $Z$ is a $\Prog$-measurable process such that  $\mathbb{P}$-a.s.
$
\int_0^T |Z_s|^2 ds <+\infty$ and  $\int_0^T |F(s, Y_s, Z_s)| ds <+\infty,
$
and (\ref{bsde}) holds  $\mathbb{P}$-a.s. for any   $t\in[0, T]$.
The first inequality above ensures  that  $Z$ is integrable with respect to $W$ in the sense of vector stochastic integration.
 As a result, $Z\cdot W$ is a continuous local martingale. 
We call $F$ the \emph{generator}, $\xi$ the \emph{terminal value} and $(\xi, \int_0^T |F(s, 0, 0)|ds)$ the \emph{data}. In our study, 
the integrability property of the data determines estimates for a solution. 
The conditions imposed on the generator are called the \emph{structure conditions}.
For notational convenience,  we sometimes write $(F, \xi)$ instead of (\ref{bsde}) to denote the BSDE with generator $F$ and terminal value $\xi$.

We are interested in BSDEs satisfying,  $\mathbb{P}$-a.s. for any $(t, y ,z)\in [0, T]\times\mathbb{R}\times\mathbb{R}^d$,
\begin{align}
\sgn(y)F(t, y, z) &\leq \alpha_t + \beta|y| + \gamma |z| +f(|y|)|z|^2,\nonumber \\
|F(t, y, z)| &\leq \alpha_t + \varphi(|y|) + \gamma |z| + f(|y|)|z|^2,  \label{sc}
\end{align}
where  $\mathbb{P}$-a.s. for any $t\in[0, T]$, $(y, z)\longmapsto F(t, y, z)$ is continuous, 
 $\alpha$ is an $\mathbb{R}^+$-valued $\Prog$-measurable  process, 
 $\varphi: \mathbb{R}^+ \rightarrow \mathbb{R}^+$ is a continuous nondecreasing function with $\varphi(0)=0$, $f(|\cdot|): \mathbb{R}\rightarrow \mathbb{R}$ is a measurable function and $\gamma \geq 0$.  
 As will be seen later,  
 the BSDEs satisfying (\ref{sc}) are solvable if $f(|\cdot|)$ belongs to  $\mathcal{I}$, the set of 
 integrable functions  from $\mathbb{R}$ to $\mathbb{R}$  which are  bounded on any compact subset of $\mathbb{R}$. 
Note that (\ref{sc})  has an even more general growth in $y$, compared to the assumption (\ref{c1:bsde2}) which is studied  by Bahlali et al \cite{BEO2014}.

Let us close this section by introducing all required notations for this chapter. 
For any random variable or process $Y$, we say $Y$ has some property if this is true except on a $\mathbb{P}$-null subset of $\Omega$. Hence we 
omit ``$\mathbb{P}$-a.s.'' in situations without ambiguity.  Define $\sgn(x): = \mathbb{I}_{\{x\neq 0 \}}\frac{x}{|x|}$.
For any c\`{a}dl\`{a}g adapted process $Y$, set $Y_{s, t}: = Y_t -Y_s$ and $Y^* : = \sup_{t\in [0, T]} |Y_t|$.
For any $\Prog$-measurable process $H$, set $|H|_{s,t}:= \int_s^t H_u du$ and $|H|_t : = |H|_{0, t}$. $\mathcal{T}$ stands for the set of  stopping times valued in $[0, T]$ and $\mathcal{S}$ denotes the space of continuous adapted processes. 
For any  local martingale $M$, we call $\{\sigma_n\}_{{n\in\mathbb{N}^+ }}\subset \mathcal{T}$ a \emph{localizing sequence} if  $\sigma_n$ increases stationarily to $T$ as $n$ goes to $+\infty$ and
$M_{\cdot \wedge \sigma_n}$ is a martingale for any $n\in\mathbb{N}^+$.
For later use,  we specify the following spaces under  $\mathbb{P}$.
\begin{itemize}
\item $\mathcal{S}^\infty$: the set of bounded processes in $\mathcal{S}$;
\item $\mathcal{S}^p (p\geq 1)$: the set of $Y\in \mathcal{S}$ with $Y^*\in \mathbb{L}^p$; 
\item $\mathcal{D}$: the set of $Y\in\mathcal{S}$ such that  $\{Y_\tau |  \tau\in\mathcal{T} \}$ is uniformly integrable;
\item $\mathcal{M}$: the space of $\mathbb{R}^d$-valued $\Prog$-measurable processes $Z$ such that  
$\mathbb{P}$-a.s.
$\int_0^T 
|Z_s|^2 ds <+\infty$;
for any $Z\in \mathcal{M}$, $Z\cdot W$ is a continuous local martingale;
\item $\mathcal{M}^p(p>0)$: the set of $Z\in\mathcal{M}$ with
\[
\norm{Z}_{\mathcal{M}^p}:=\mathbb{E}\Big[\Big(\int_0^T |Z_s|^2 ds \Big)^{\frac{p}{2}}\Big]^{\frac{1}{p}\wedge 1} <+\infty;
\]
in particular, $\mathcal{M}^2$  is a Hilbert space;
\item $\mathcal{C}^p(\mathbb{R})$: the space of $p$ times continuously differentiable functions from $\mathbb{R}$ to $\mathbb{R}$;
\item$\mathcal{W}_{1, loc}^2(\mathbb{R})$: the Sobolev space  of measurable maps $u: \mathbb{R}\rightarrow \mathbb{R}$ such that both $u$ and its generalized derivatives $u^\prime, u^{\prime\prime}$ belong to $\mathbb{L}^1_{loc}(\mathbb{R})$.\\
\end{itemize}
 
 The above spaces 
are  Banach (respectively complete) under suitable norms (respectively metrics); we will not present these facts in more detail since
 they are not involved in our study.  We call $(Y, Z)$ a $\mathbb{L}^p$ \emph{solution} to (\ref{bsde}) if $(Y, Z)\in\mathcal{S}^p\times\mathcal{M}^p$.  This definition simply comes from the fact that its existence is ensured by data in $\mathbb{L}^p$.
 Analogously to most papers on $\mathbb{R}$-valued quadratic BSDEs, 
our existence result essentially relies  on 
 the monotone stability result of quadratic BSDEs; see, e.g.,  Kobylanski \cite{K2000},  Briand and  Hu \cite{BH2008} or Section \ref{section22}, Chapter \ref{sbsde}.
  
\section{Functions of Class  $\mathcal{I}$}
\label{I}
In this section, we introduce the basic ingredients used to treat the quadratic generator in  (\ref{sc}).
We recall that $\mathcal{I}$ is  the set of 
integrable  functions  from $\mathbb{R}$ to $\mathbb{R}$  which are bounded on any compact subset of $\mathbb{R}$.

{\bf $u^f$ Transform.} For any $f\in \mathcal{I}$, define   $u^f: \mathbb{R}\rightarrow \mathbb{R}$ and $M^f$ by
\begin{align*}
u^f (x)&: =  \int_0^x \exp{\Big(2\int_0^y f(u)du\Big)}dy,  \\
M^f&: = \exp\Big( 2\int_{-\infty}^{\infty} |f(u)| du\Big).
\end{align*}
Obviously, $1\leq M^f < + \infty$.
Moreover, the following properties hold by simple computations. 
\begin{enumerate}
\item [\rm{(i)}] $u\in \mathcal{C}^1(\mathbb{R}) \cap \mathcal{W}_{1, loc}^2 (\mathbb{R})$ and 
$
u^{\prime\prime}(x) =2f(x)u^\prime(x)$  a.e.;
if  $f$ is continuous, then $u \in \mathcal{C}^2 (\mathbb{R})$;
\item [\rm{(ii)}] $u$ is strictly increasing and bijective from $\mathbb{R}$ to $\mathbb{R}$;
\item [\rm{(iii)}] $u^{-1} \in  \mathcal{C}^1(\mathbb{R}) \cap \mathcal{W}_{1, loc}^2 (\mathbb{R})$; if  $f$ is continuous, then $u^{-1} \in \mathcal{C}^{2}(\mathbb{R})$;
\item [\rm{(iv)}] 
$
\frac{|x|}{M}  \leq |u(x)| \leq M|x|
$ 
and
$\frac{1}{M}\leq u^\prime(x)\leq M.$
\end{enumerate}

{\bf $v^f$ Transform.} For any $f\in \mathcal{I}$, define $v^f: \mathbb{R}\rightarrow \mathbb{R}^+$ by
\begin{align*}
v^f(x): = \int_0^{|x|} u^{(-f)}(y) 
\exp\Big(2\int_0^y f(u)du     \Big)
  dy.
\end{align*}
Set $v:=v^f$. Simple computations  give
\begin{enumerate}
\item [(i)]$v\in \mathcal{C}^1(\mathbb{R}) \cap \mathcal{W}_{1, loc}^2 (\mathbb{R})$  and
$v^{\prime\prime}(x)-2f(|x|)|v^\prime(x)|=1$ a.e.;
if  $f$ is continuous, then 
$v\in \mathcal{C}^2(\mathbb{R});$
\item [(ii)] $v(x)\geq 0, \sgn(v^\prime(x))=\sgn(x)$ and $v^{\prime\prime}(0)=1;$
 \item [(iii)]
 $
 \frac{x^2}{2M^2}\leq v(x) \leq \frac{M^2x^2}{2}$ and $\frac{|x|}{M^2}\leq|v^\prime(x)|\leq M^2 |x|$.
\end{enumerate}
In the sequel of our study, $u^f$ and $v^f$ exclusively stand for the above transforms associated with $f\in\mathcal{I}$. Hence in situations without ambiguity, we  denote $u^f, v^f, M^f$ by $u, v, M$, respectively.
\section{Krylov Estimate and  the It\^{o}-Krylov Formula}
\label{itokrylov}
The first auxiliary result  is
 the Krylov estimate.
Later, it is used to 
 prove an  It\^{o}'s-type formula for functions in $\mathcal{C}^1(\mathbb{R})\cap \mathcal{W}_{1, loc}^2(\mathbb{R})$. This helps to deal with 
 (possibly discontinuous) quadratic generators.
As the second application, we derive a generalized It\^{o} formula for $y\mapsto |y|^p(p\geq 1)$ which is not smooth enough for $1\leq p <2$. This is a basic tool to study $\mathbb{L}^p(p\geq 1)$ solutions.

 To allow the existence of a local time  in  particular situations, we  study equations of type
\begin{align}
Y_t = \xi + \int_t^T F(s, Y_s, Z_s) ds +\int_t^T dC_s -\int_t^T Z_s dW_s, \label{fc}
\end{align}
where $C$ is a continuous adapted process of finite variation. We denote its total variation process by $V_\cdot(C)$.  Likewise,  sometimes  we denote   (\ref{fc}) by 
$(F, C, \xi)$. The solution to (\ref{fc}) is defined analogously to that to   (\ref{bsde}).

Now we prove the Krylov estimate for (\ref{fc}).
A more complicated version  not needed  for our study can be found in Bahlali et al \cite{BEO2014}. 

\begin{Lemma}[Krylov Estimate] \label{krylovestimate} 
For any measurable function $\psi: \mathbb{R}\rightarrow \mathbb{R}^+$,
\begin{align}
\mathbb{E} \Big[  \int_0^{ \tau_{m}} \psi(Y_s)|Z_s|^2 ds \Big] \leq 6m \norm{\psi}_{\mathbb{L}^1([-m, m])}, \label{psifun}
\end{align}
where  $\tau_m$ is a stopping time defined by 
\begin{align*}
\tau_{m} : = \inf\Big \{t \geq 0: |Y_t| + V_t(C) + \int_0^t |F(s, Y_s, Z_s)| ds \geq m \Big\} \wedge T.
\end{align*}
\end{Lemma}
\begin{Proof}
Without loss of generality we assume
$\norm{\psi}_{\mathbb{L}^1([-m, m])} <+\infty$. 
For   each $n\in\mathbb{N}^+$, set 
\begin{align*}
\tau_{m,n} &: = \tau_{m} \wedge \inf \Big\{t\geq 0: \int_0^t |Z_s|^2 ds \geq n \Big\}.
\end{align*} 
Let $a\in [-m, m]$. By Tanaka's formula, 
\begin{align}
(Y_{t\wedge \tau_{m,n}} - a)^{-}  & =  (Y_0 - a)^- - \int_0^{t\wedge\tau_{m,n}} \mathbb{I}_{\{Y_s < a\} }dY_s + \frac{1}{2}L_{t\wedge \tau_{m,n}}^a(Y)\nonumber\\ 
& =  (Y_0 - a)^- + \int_0^{t\wedge \tau_{m,n}} \mathbb{I}_{\{Y_s < a\}} F(s, Y_s, Z_s)ds +\int_0^{t\wedge \tau_{m, n}} \mathbb{I}_{\{ Y_s < a\}}dC_s \nonumber\\
&-\int_0^{t\wedge\tau_{m,n}}\mathbb{I}_{\{Y_s < a\}} Z_s dW_s 
+\frac{1}{2} L_{{t\wedge \tau_{m,n}}}^a (Y), \label{ito-}
\end{align}
where $L^{a}(Y)$ is the local time of $Y$ at $a$. To estimate the local time, we put it on the left-hand side and the rest terms on the right-hand side. 
Since $x \mapsto (x-a)^-$ is Lipschitz-continuous,  we deduce from the definition of $\tau_{m, n}$ that
\[
(Y_0 - a)^- - (Y_{t\wedge \tau_{m, n}} - a)^- \leq | Y_0 -Y_{t\wedge{\tau_{m, n}}}| \leq 2m.
\]
Meanwhile, the definition of $\tau_m$ also implies that the sum of the $ds$-integral and $dC$-integral is bounded by $m$. Hence, we have
\[
\mathbb{E}\big[L_{t\wedge \tau_{m,n}}^a(Y)\big] \leq 6m. 
\]
By Fatou's lemma applied to the sequence indexed by $n$, 
\[
\sup_{a\in[-m,  m]}\mathbb{E}\big[L_{t\wedge \tau_{m}}^a(Y)\big] \leq 6m.
\]
We then use  
time occupation formula for continuous semimartingales (see Chapter VI., Revuz and  Yor  \cite{RY1999})  and the above inequality to obtain
\begin{align*}
\mathbb{E}\Big[ \int_0^{T\wedge\tau_{m}} \psi(Y_s)|Z_s|^2 ds  \Big] 
 &= 
\mathbb{E}\Big[ \int_{-m}^{m} \psi(x)L_{T\wedge \tau_m}^x(Y)dx  \Big]\\
&= \int_{-m}^m \psi(x) \mathbb{E}\big[L_{T\wedge \tau_{m}}^x (Y)\big] dx\\
&\leq 6m \norm{\psi}_{\mathbb{L}^1([-m, m])}.
\end{align*}
\qed
\end{Proof}

As an immediate consequence of Lemma 
\ref{krylovestimate}, we have $\mathbb{P}$-a.s.
\begin{align}
\int_0^T \mathbb{I}_{\{ Y_s \in A \}}|Z_s|^2 ds = 0, \label{Aneg}
\end{align}
 for any $A\subset\mathbb{R}$ with null Lebesgue measure.  This will be used later several times.

  Given Lemma \ref{krylovestimate}, we turn to the main results of this section.

\begin{thm}[It\^{o}-Krylov Formula]\label{ik} If $(Y, Z)$ is a solution to $(F, C, \xi)$, then
 for any  $u\in\mathcal{C}^1(\mathbb{R})\cap \mathcal{W}_{1, loc}^2(\mathbb{R})$,  we have $\mathbb{P}$-a.s. for all $t\in [0, T]$, 
\begin{align}
u(Y_t) = u(Y_0) + \int_0^t u^\prime (Y_s) dY_s +\frac{1}{2} \int_0^t u^{\prime\prime}(Y_s)|Z_s|^2 ds. \label{ikformula} 
\end{align}

\end{thm}
\begin{Proof}
We use $\tau_{m}$  defined in Lemma \ref{krylovestimate} (Krylov estimate). 
Note that  $\tau_{m}$ increases stationarily to $T$ as $m$ goes to $+\infty$. It is therefore sufficient to prove the equality for $u(Y_{t\wedge\tau_m})$. To this end we use an approximation procedure. We consider $m$ such that  $\mathbb{P}$-a.s. $m\geq |Y_0|$.
Let $u_n$ be a sequence of functions in  $\mathcal{C}^2(\mathbb{R})$ satisfying
\begin{enumerate}
\item [(i)] $u_n$ converges uniformly to $u$ on $[-m, m]$; 
\item [(ii)] $u_n^\prime$ converges uniformly to $u^\prime$ on $[-m, m]$;
\item [(iii)]$u_n^{\prime\prime}$ converges in $\mathbb{L}^1([-m, m])$ to $u^{\prime\prime}$.
\end{enumerate}
By It\^{o}'s formula,
\[
u_n(Y_{t\wedge \tau_{m}}) = u_n(Y_0) + \int_0^{t\wedge \tau_{m}} u_n^\prime (Y_s) dY_s +\frac{1}{2} \int_0^{t\wedge \tau_{m}} u_n^{\prime\prime}(Y_s)|Z_s|^2 ds.
\]
Due to  (i) and $|Y_{t \wedge \tau_m}| \leq m$ , 
$u_n (Y_{\cdot \wedge \tau_{m}})$  converges to $ u(Y_{\cdot \wedge \tau_{m}})$ $\mathbb{P}$-a.s. 
uniformly on $[0, T]$
as $n$ goes to $+\infty$; the second term converges in probability to 
\[
\int_0^{t \wedge \tau_{m}} u^\prime (Y_s) dY_s
\]
by (ii)
 and dominated convergence for stochastic integrals; the last term converges in probability to 
\[
\frac{1}{2}\int_0^{t\wedge \tau_{m}} u^{\prime\prime}(Y_s) |Z_s|^2 ds
\]
due to  (iii) and Lemma \ref{krylovestimate}.
Indeed,  Lemma \ref{krylovestimate} implies
\begin{align*}
\mathbb{E}\Big[\int_0^{\tau_m} |u^{\prime\prime}_n-u^{\prime\prime}|(Y_s)|Z_s|^2 ds  \Big]
\leq 6m \norm{u^{\prime\prime}_n-u^{\prime\prime}}_{{\mathbb{L}^1{([-m, m])}}}.
\end{align*} Hence  collecting these convergence results  gives (\ref{ikformula}). 
By the continuity of both sides of (\ref{ikformula}), the quality also holds $\mathbb{P}$-a.s. for all $t\in[0, T]$.
\qed
\end{Proof}

To study $\mathbb{L}^p(p\geq 1)$ solutions we now prove an It\^o's-type formula for $y \mapsto |y|^p(p\geq 1)$ which
is not smooth enough for $1\leq p<2$.  The proof for 
 multidimensional It\^{o} processes can be found, e.g., in  Briand et al \cite{B2003}. In contrast to their approach, 
we give a novel and simpler proof  for BSDE framework but point out that it can be also extended to  It\^{o} processes. 
\begin{Lemma}\label{lpinequality}
Let $p \geq 1$. If $(Y, Z)$ is a solution to  $(F, C, \xi)$,
then
\begin{align}
|Y_t|^p &+ \frac{p(p-1)}{2} \int_t^T  \mathbb{I}_{\{ Y_s \neq 0\}}|Y_s|^{p-2}  |Z_s|^2ds  \nonumber\\
&   = |\xi|^p - p \int_t^T \sgn(Y_s)|Y_s|^{p-1} dY_s-\mathbb{I}_{\{ p =1 \}}\int_t^T dL_s^0(Y), \label{lemma23}
\end{align}
where $L^0(Y)$ is the local time of 
$Y$ at $0$. 
\end{Lemma}
\begin{Proof}
(i). $p = 1 $. This is immediate from Tanaka's formula. 

(ii). $ p> 2$. $y \mapsto |y|^p \in \mathcal{C}^2(\mathbb{R})$. Hence this is immediate from It\^{o}'s formula.

(iii). $p=2$.  $y \mapsto |y|^p \in \mathcal{C}^2(\mathbb{R})$. 
Due to (\ref{Aneg}), 
$\int_0^\cdot|Y_s|^{p-2}  |Z_s|^2ds$
is indistinguishable from $\int_0^\cdot  \mathbb{I}_{\{ Y_s \neq 0\}}|Y_s|^{p-2}  |Z_s|^2ds$. 
Then the inequality is immediate from It\^o's formula.
    
(iv). $1<p<2$. We use an approximation argument. Define 
\[
u_\epsilon (y): = \big(y^2 + \epsilon^2\big)^{\frac{1}{2}}.
\]
Then for any $\epsilon >0$, we have $u_\epsilon^p\in \mathcal{C}^2(\mathbb{R}).$ By It\^{o}'s formula, 
 \begin{align}
 u^p_\epsilon (Y_t) &=  u^p_\epsilon (\xi) - p \int_t^T Y_s u^{p-2}_\epsilon (Y_s) d Y_s -   \frac{1}{2} \int_t^T \big( p u_\epsilon^{p-2}(Y_s)+p(p-2)|Y_s|^2u_\epsilon^{p-4}(Y_s)\big)
   |Z_s|^2   ds.  \label{uepsilon}
 \end{align} Now we send $\epsilon$ to $0$.
 $u_\epsilon(y) \longrightarrow |y|$ pointwise implies
 $
u_\epsilon(Y_t)^p \longrightarrow
|Y_t|^p$ and  $\ u_\epsilon(\xi)^p
\longrightarrow |\xi|^p$ pointwise on $\Omega$.
  Secondly, 
$yu_\epsilon^{p-2}(y) \longrightarrow \sgn(y)|y|^{p-1}$ pointwise implies 
by dominated convergence for stochastic integrals that
\begin{align*}
&\int_t^T Y_s \sgn(Y_s) u^{p-2}_\epsilon (Y_s) dY_s {\longrightarrow}\int_t^T |Y_s|^{p-1}dY_s \ \text{in probability.}
\end{align*}
To prove that the $ds$-integral in (\ref{uepsilon}) also converges, we split it into two parts and argue their convergence respectively.
 Note that  
 \begin{align}
pu_\epsilon^{p-2}(Y_s)+p(p-2)|Y_s|^2u_\epsilon^{p-4}(Y_s) = p\epsilon^2 u_\epsilon^{p-4}(Y_s) + p(p-1) |Y_s|^2 u_\epsilon^{p-4} (Y_s). \label{2ndorder}
\end{align}
For the second term on the right-hand side of (\ref{2ndorder}), we have
\[
|Y_s|^2 u_\epsilon^{p-4} (Y_s)  =  \mathbb{I}_{\{Y_s \neq 0\}}  |Y_s|^{p-2} \Big|\frac{|Y_s|}{u_\epsilon(Y_s)} \Big|^{4-p}.
\]
Since
$
\frac{|y|}{u_\epsilon(y)} {\nearrow} \mathbb{I}_{\{ y \neq 0\}}
$, 
 monotone convergence gives
\[
\int_t^T |Y_s|^2 u_\epsilon^{p-4}(Y_s) |Z_s|^2ds
\longrightarrow 
\int_t^T \mathbb{I}_{\{Y_s \neq 0 \}}|Y_s|^{p-2}|Z_s|^2  ds \  \text{pointwise in}\ \Omega. 
\]
It thus remains to prove  the $ds$-integral concerning the first term on the right-hand side of (\ref{2ndorder})
  converges to $0$. To this end,  
 we use Lemma \ref{krylovestimate} (Krylov estimate)
 and the same localization procedure.
This gives
\begin{align*}
\mathbb{E} \Big[ \int_0^{\tau_m} \epsilon^2 u_\epsilon^{p-4}(Y_s) |Z_s|^2ds  \Big] &\leq 6m \epsilon^2 \int_{-m}^{m} (x^2 + \epsilon^2)^{\frac{p-4}{2}} dx\\
& \leq  12m\epsilon^2\int_0^m  (x^2 + \epsilon^2)^{\frac{p-4}{2}} dx\\
& \leq  12\cdot 2^{\frac{4-p}{2}}m\epsilon^2 \int_{0}^{m} (x+\epsilon)^{p-4} dx\\
& \leq  12\cdot 2^{\frac{4-p}{2}}m\epsilon^2 \int_{\epsilon}^{m+\epsilon} x^{p-4} dx\\
&= \frac{12\cdot 2^{\frac{4-p}{2}} m}{p-3} \big (\epsilon^2  (m+\epsilon)^{p-3} - \epsilon^{p-1}\big),
\end{align*}
which, due to $1<p$, converges to $0$ as $\epsilon $ goes to $0$. Hence 
$
\int_0^{\cdot}\epsilon^2 u_\epsilon^{p-4}(Y_s)|Z_s|^2 ds 
$
 converges $u.c.p$ to $0$. 
Collecting all convergence results above gives (\ref{lemma23}). By the continuity
 of each term in (\ref{lemma23}), the equality
also holds $\mathbb{P}$-a.s. for all $t\in [0, T]$.
\qed
\end{Proof}
\section{$\mathbb{L}^p(p \geq 1)$ Solutions to Purely Quadratic BSDEs}
\label{degenerate}
Before turning to the main results, we partially extend the existence and uniqueness result for purely quadratic BSDEs studied by  
Bahlali et al \cite{BEO2014}. Later, we present their natural extensions and the motivations of our work. These BSDEs are called purely quadratic, since the generator takes the form 
$F(t, y, z)=f(y)|z|^2$.   
The solvability simply comes from the function $u^f$ defined in Section \ref{I} which transforms  better known BSDEs to  $(f(y)|z|^2, \xi)$ by It\^{o}-Krylov formula. 
\begin{thm}
\label{degeneratethm}
Let $f\in \mathcal{I}$ and $\xi \in \mathbb{L}^p (p \geq 1)$. Then there exists a unique solution to
\begin{align}
Y_t = \xi + \int_t^T f(Y_s)|Z_s|^2ds -\int_t^T Z_s dW_s \label{pqbsde}.
\end{align}
Moreover, if $p>1$, the solution belongs to $\mathcal{S}^p\times\mathcal{M}^p$; if $p=1$, the solution belongs to $\mathcal{D}\times\mathcal{M}^q$ for any $q\in(0, 1)$.
\end{thm}
\begin{Proof}
Let $u: = u^f$.
Then $u, u^{-1}\in
\mathcal{C}^1(\mathbb{R})\cap \mathcal{W}_{1, loc}^2(\mathbb{R})$. The existence and uniqueness result can be seen as a one-on-one correspondence between solutions to  BSDEs.

(i). Existence. $|u(x)| \leq M |x|$ implies  
$u(\xi) \in \mathbb{L}^p$.
By  It\^{o} representation theorem,  
there exists a unique pair $(\widetilde{Y}, \widetilde{Z})$ which solves $(0, u(\xi))$, i.e.,
\begin{align}
d\widetilde{Y}_t = \widetilde{Z}_t dW_t, \ \widetilde{Y}_T = u(\xi). \label{tbsde}
\end{align}
We aim at proving
\begin{align}(Y, Z):=
(u^{-1}(\widetilde{Y}),
\frac{\widetilde{Z}}{u^\prime(u^{-1}(\widetilde{Y}))}) \label{yz}
\end{align}
solves (\ref{pqbsde}).  
 It\^o-Krylov formula applied to $Y_t = u^{-1}(\widetilde{Y}_t)$ yields 
\begin{align}
dY_t =\frac{1}{u^\prime(u^{-1}(\widetilde{Y}_t))} d\widetilde{Y}_t - \frac{1}{2}\Big(\frac{1}{u^\prime(u^{-1}(\widetilde{Y}_t))}\Big)^2 \frac{u^{\prime\prime}(u^{-1}(\widetilde{Y}_t))}{u^\prime(u^{-1}(\widetilde{Y}_t))} |\widetilde{Z}_s|^2 ds. \label{ytt}
\end{align}
To simplify (\ref{ytt}) let us recall that $u^{\prime\prime}(x)= 2f(x)u^{\prime}(x)$ a.e. Hence (\ref{yz}), (\ref{ytt}) and (\ref{Aneg}) give
\[
dY_t = -f(Y_t)|Z_t|^2 dt +  Z_t dW_t, \ Y_T = \xi,
\]
i.e., $(Y, Z)$ solves (\ref{pqbsde}).

(ii). Uniqueness. Suppose $(Y, Z)$ and  $(Y^\prime, Z^\prime)$ are solutions to (\ref{pqbsde}). By It\^o-Krylov formula applied to $u(Y)$ and  $u(Y^\prime)$, we deduce that 
$(u(Y), u^\prime(Y)Z)$ and $(u(Y^\prime), u^\prime(Y^\prime)Z^\prime)$
 solve
$(0, u(\xi))$. But from  (i) it is known that they coincide. Transforming $u({Y})$ and $u(Y^\prime)$ via the bijective function $u^{-1}$ yields the uniqueness result.

(iii). We prove the estimate for the unique solution $(Y, Z)$. For $p>1$, Doob's $\mathbb{L}^p(p>1)$ maximal inequality used to (\ref{tbsde}) implies  $(\widetilde{Y}, \widetilde{Z})\in\mathcal{S}^p\times\mathcal{M}^p$. Hence $({Y}, {Z})\in\mathcal{S}^p\times\mathcal{M}^p$, due to $|u^\prime(x)|\geq \frac{1}{M}$ and $|u^{-1}(x)| \leq M |x|$.
 For $p=1$, $\widetilde{Y}\in \mathcal{D}$ since it is a martingale on $[0, T]$. In view of the above properties of $u$ we have $Y\in\mathcal{D}$. The estimate for $Z$ is immediate from
Lemma 6.1, Briand et al \cite{B2003} which is a version of $\mathbb{L}^p(0<p<1)$ maximal inequality for martingales.
\qed
\end{Proof}
\begin{remm}
If $\xi$ is a general $\mathcal{F}_T$-measurable random variable,
 Dudley representation theorem (see Dudley \cite{d1977}) implies that 
 there  still exists a solution to (\ref{tbsde}) and hence a solution to (\ref{pqbsde}). However, the solution in general is not unique.
 
The proof of Theorem \ref{degeneratethm} indicates that 
$f$  being bounded on compact subsets of $\mathbb{R}$ is not needed 
for he existence and uniqueness result of purely quadratic BSDEs.
\end{remm}
\begin{prop}[Comparison]
\label{dcompare}
Let $f, g\in \mathcal{I}$, $\xi, \xi^\prime \in\mathbb{L}^p(p\geq 1)$ and
 $(Y, Z)$, $(Y^\prime, Z^\prime)$ be the unique solutions to $(f(y)|z|^2, \xi)$, $(g(y)|z|^2, \xi^\prime)$, respectively. 
If $f\leq g$ a.e. and $\mathbb{P}$-a.s. $\xi\leq \xi^\prime$, then  $\mathbb{P}$-a.s. $Y_\cdot \leq Y^\prime_\cdot$. 
\begin{Proof}
Again we transform so as to compare better known BSDEs.  Let us fix $t\in [0, T]$ and set $u:=u^f$. For any  stopping time valued  in $[t, T]$, 
 It\^{o}-Krylov formula yields
\begin{align*}
u (Y_{t}^\prime)  
&=
u(Y_{\tau}^\prime) + \int_{t}^\tau \Big(u^\prime(Y_s^\prime)g(Y_s^\prime) |Z_s^\prime|^2
-\frac{1}{2}u^{\prime\prime}(Y_s^\prime)|Z_s^\prime|^2\Big)
 ds 
-\int_{t}^\tau u^\prime (Y_s)Z_s^\prime dW_s.
\\
&=
u(Y_{\tau}^\prime)+\int_{t}^\tau u^\prime (Y_s^\prime) 
\big(g(Y_s^\prime)-f(Y_s^\prime)\big) |Z_s^\prime|^2 ds
-\int_{t}^\tau u^\prime (Y_s)Z_s^\prime dW_s \\
&\geq
u(Y_{\tau}^\prime) -\int_{t}^\tau u^\prime (Y_s)Z_s^\prime dW_s,
\end{align*}
where the last two lines are due to
$u^{\prime\prime}(x)  = 2 f(x)u^\prime(x)$ a.e., $g\geq f$ a.e. and (\ref{Aneg}). In the next step, we want to eliminate the local martingale part by a localization procedure. Note that 
$\int_t^\cdot u^{\prime}(Y_s)Z_s^\prime dW_s$ is a local martingale on $[t, T]$. Set 
$\{\tau_n\}_{n\in\mathbb{N}^+}$
to be its localizing sequence on $[t, T]$. Replacing $\tau$ by $\tau_n$ in the above inequality thus gives $\mathbb{P}$-a.s.
\[
u(Y_t^\prime) \geq \mathbb{E}\big[u(Y_{t\wedge \tau_n}^\prime)\big|\mathcal{F}_t\big].
\]
This implies that, for any $A\in\mathcal{F}_t$,  we have
\[
\mathbb{E}\big[ u(Y_t^\prime)\mathbb{I}_A\big] \geq 
\mathbb{E}\big[ u(Y_{t\wedge\tau_n}^\prime)\mathbb{I}_A  \big].
\]
Since $u(Y^\prime)\in\mathcal{D}$, we can use Vitali convergence theorem to obtain
\[
\mathbb{E}\big[ u(Y_t^\prime)\mathbb{I}_A \big]
\geq \mathbb{E}\big[ u(\xi^\prime)  \mathbb{I}_A\big] =
\mathbb{E}\big[\mathbb{E}\big[ u(\xi^\prime)\big| \mathcal{F}_t \big]  \mathbb{I}_A\big].
\]
Note that this inequality holds for any $A\in\mathcal{F}_t$. Hence, by choosing $A=\{u(Y_t^\prime) < \mathbb{E}[u(\xi^\prime)|\mathcal{F}_t]\}$, we obtain
$
u(Y_t^\prime) \geq \mathbb{E}\big[u(\xi^\prime)\big|\mathcal{F}_t\big].$  Since $\xi^\prime \geq \xi$ and $u$ is increasing, we further have 
$u(Y_t^\prime) \geq 
\mathbb{E}\big[ u(\xi)\big|\mathcal{F}_t\big].$ Let us recall that, by Theorem \ref{degeneratethm}, 
$(u(Y), u^\prime (Y)Z)$ is the unique solution to $(0, u(\xi))$. Hence, $
u(Y_t^\prime) \geq u(Y_t).$
Transforming both sides via the bijective increasing function $u^{-1}$
yields  $\mathbb{P}$-a.s. $Y_t\leq Y_t^\prime$. By the continuity of $Y$ and $Y^\prime$ we have $\mathbb{P}$-a.s. $Y_\cdot\leq Y^\prime_\cdot$.
\qed
\end{Proof}
\end{prop}
\begin{remm}
In   Proposition \ref{dcompare}, we rely on the fact that 
$\mathbb{P}$-a.s.
\begin{align}\int_0^\cdot 
\Big(\frac{1}{2}u^{\prime\prime}(Y_s^\prime) - f(Y_s^\prime) u^{\prime}(Y_s^\prime)\Big)|Z_s^\prime|^2 ds=0, \label{aac}
\end{align}
 even though   
$u^{\prime\prime}(x)  = 2 f(x)u^\prime(x)$  only holds  almost everywhere on $\mathbb{R}$.
Here we prove it. 
 Let  $A$ be the subset of $\mathbb{R}$ on which $u^{\prime\prime}(x)  = 2 f(x)u^\prime(x)$ fails.  Hence, 
\[\int_0^\cdot \mathbb{I}_{\{ Y^\prime_s\in \mathbb{R}\backslash A\}}
\Big|\frac{1}{2}u^{\prime\prime}(Y_s^\prime) - f(Y_s^\prime) u^{\prime}(Y_s^\prime)\Big||Z_s^\prime|^2 ds=0.\]
Meanwhile, by (\ref{Aneg}), we have 
$\mathbb{P}$-a.s.
\[\int_0^\cdot \mathbb{I}_{\{ Y^\prime_s\in A\}}
\Big|\frac{1}{2}u^{\prime\prime}(Y_s^\prime) - f(Y_s^\prime) u^{\prime}(Y_s^\prime)\Big||Z_s^\prime|^2 ds=0.\]
Hence, (\ref{aac}) holds 
$\mathbb{P}$-a.s.
 This fact also applies to Theorem \ref{degeneratethm} and all results in the sequel of our study.  
\end{remm}
To end our discussions on purely quadratic BSDEs we give some examples.
\begin{ex} Let $\xi\in\mathbb{L}^p(p\geq 1)$. Then   Theorem {\rm\ref{degeneratethm}} holds for $(F, \xi)$, where $F$  verifies any one of the following
\begin{itemize}
\item $F(y, z) = \sin(y) \mathbb{I}_{[-\pi, \frac{\pi}{2}]}(y)|z|^2;$
\item $F(y, z)= \big(\mathbb{I}_{[a, b]} -\mathbb{I}_{[c, d]}\big)(y)|z|^2$ for some $a<b$ and $c<d;$
\item
$F(y, z) =\mathbb{I}_{\{y\neq 0\}}\frac{1}{(1+y^2)\sqrt{|y|}} |z|^2 
+ \mathbb{I}_{\{y=0\}}|z|^2.
$.
\end{itemize}
\end{ex}
Theorem \ref{degeneratethm}   and Proposition \ref{dcompare} are
based on a one-on-one correspondence between solutions (respectively the unique solution) to BSDEs. Hence it is natural to generalize as follows.  Set  $f\in\mathcal{I}, u:=u^f$,  $F(t, y, z):=  G(t, y, z) + f(y)|z|^2$  and
\begin{align}
\widetilde{F} (t, y, z): =
u^\prime (u^{-1}(y)) G(t, u^{-1}(y),\frac{z}{u^\prime (u^{-1}(y))}).\label{adds}
\end{align}
If $G$ ensures the existence of a solution to   $(\widetilde{F}, u(\xi))$, we can transform it via ${u}^{-1}$
to  a solution to $(F, \xi)$. An example is that
$G$ is of continuous linear growth in $(y, z)$ where 
the existence of a maximal (respectively minimal) solution to $(\widetilde{F}, u(\xi))$ can be proved in the spirit of Lepeltier and San Martin \cite{LS1997}.

When the generator is continuous in $(y, z)$, a more general situation is    linear-quadratic growth, i.e., 
\begin{align}
|H(t, y, z)| \leq \alpha_t + \beta|y| + \gamma |z| + f(|y|)|z|^2 : =F(t, y, z). \label{lq}
\end{align}
The existence result  then consists of viewing the maximal (respectively minimal) solution to $(F, \xi^+)$ (respectively $(-F, -\xi^-)$) as a priori bounds for solutions to $(H, \xi)$, 
and using a combination of a localization procedure and  the monotone stability result developed by  Briand and  Hu \cite{BH2006}, \cite{BH2008}.
For details the reader shall refer to Bahlali et al \cite{BEO2014}.

However, either an additive structure in (\ref{adds}) or a linear-quadratic growth (\ref{lq})
 is too restrictive 
   and  uniqueness  is not available in general.  Considering this limitation, we devote Section \ref{lpsolution} to  the solvability under  milder structure conditions.  
\section{$\mathbb{L}^p(p > 1)$ Solutions to Quadratic BSDEs}
\label{lpsolution} 
With the preparatory work in Section \ref{lpp}, \ref{I}, \ref{itokrylov}, \ref{degenerate}, we study  $\mathbb{L}^p (p > 1)$ solutions to quadratic BSDEs under general assumptions. We deal with the quadratic generators in the spirit of   
 Bahlali et al	\cite{BEO2014}, 
derive the estimates in the spirit of Briand et al \cite{B2003} and prove the existence and uniqueness result in the spirit of  Briand et al \cite{BH2006}, \cite{BH2008}, \cite{BLS2007}. This section can also be seen as a generalization of these works.  
The following assumptions on $(F, \xi)$  ensure the  a priori estimates and an existence result. 
\begin{as}\label{lpas1}
Let $p\geq 1$.
There exist  $\beta \in \mathbb{R}, \gamma \geq 0$, an $\mathbb{R}^+$-valued $\Prog$-measurable process $\alpha$, $f(|\cdot|)\in\mathcal{I}$  and a continuous nondecreasing function $\varphi: \mathbb{R}^+ \rightarrow \mathbb{R}^+$ with $\varphi(0)=0$ such that     
$  |\xi|  + |\alpha |_T \in \mathbb{L}^p$ and   $\mathbb{P}$-a.s.
  \begin{enumerate} 
\item [(i)] for any $t\in[0, T]$, $(y, z)\longmapsto F(t, y, z)$ is continuous; 
  \item [(ii)] $F$ is ``monotonic''  at  $y=0$, i.e.,  for any $(t,y, z)\in [0, T]\times\mathbb{R}\times\mathbb{R}^d$,
    \[\sgn (y) F(t, y, z) \leq  \alpha_t + \beta |y| + \gamma |z| + f(|y|) |z|^2; \]
  \item [(iii)]  for any $(t, y, z)\in[0, T]\times\mathbb{R}\times\mathbb{R}^d$,
  \[|F(t, y, z)| \leq \alpha_t + \varphi(|y|) + \gamma |z| + f(|y|)|z|^2.\] 
\end{enumerate}
\end{as} 

It
is worth noticing that 
given \ref{lpas1}(iii) and $f=0$, 
\ref{lpas1}(ii) is a consequence of  $F$ being monotonic at $y=0$. Indeed,  
\[
\sgn (y-0)\big( F(t, y, z) - F(t, 0, z)\big) \leq \beta|y| 
\]
implies
\begin{align*}
\sgn(y) F(t, y, z) &\leq
F(t, 0, z) + \beta|y| \\
&\leq \alpha_t + \beta |y| +\gamma |z|.
\end{align*}
This explains why  we keep  saying that  $F$ is monotonic  at $y=0$, even though $y$ also appears in the quadratic term.
Secondly, our results don't rely on the specific choice of $\varphi$. 
Hence the growth condition in $y$ can be arbitrary as long as \ref{lpas1}(i)(ii) hold.
Assumptions of this type for different settings can also be found in, e.g.,  \cite{P1999}, \cite{BC2000}, \cite{B2003}, \cite{BLS2007}, \cite{BH2008}.
 Finally,  $f$ can be discontinuous; $f(|\cdot|)$ being $\mathbb{R}^+$-valued appears more naturally in the growth condition.
\begin{Lemma}
[A Priori Estimate (i)]
\label{lpestimate1} Let $p \geq 1$ and  {\rm \ref{lpas1}} hold for $(F, \xi)$. If $(Y, Z) \in \mathcal{S}^p \times \mathcal{M} $ is a solution to $(F, \xi)$, then 
\begin{align*}
\mathbb{E}\Big[\Big( \int_0^T |Z_s|^2 ds \Big)^{\frac{p}{2}}\Big] + \mathbb{E}\Big[ \Big(\int_0^T  f(|Y_s|) |Z_s|^2 ds \Big)^{p} \Big] 
\leq c\Big( \mathbb{E}\big[(Y^*)^p + |\alpha|_T^p\big]\Big),
\end{align*}
where  $c$ is a  constant  only depending on
$T, M^{f(|\cdot|)}, \beta, \gamma, p$. 
\end{Lemma}
\begin{Proof} 
Set $v:=v^{f(|\cdot|)}$ and $M:=M^{f(|\cdot|)}$.
For any $\tau\in\mathcal{T}$,  It\^o-Krylov formula  yields
\begin{align}
v(Y_0) &= v(Y_\tau)
+\int_0^\tau  v^\prime(Y_s) F(s, Y_s, Z_s) ds \nonumber
\\ 
&-\frac{1}{2}\int_0^\tau v^{\prime\prime}(Y_s)|Z_s|^2 ds 
- \int_0^\tau v^\prime(Y_s)Z_s dW_s.\label{itov1}
\end{align}
Due to $\sgn(v^\prime(x))=\sgn(x)$ and \ref{lpas1}(ii),  we have
\begin{align}
v^\prime(Y_s)F(s, Y_s, Z_s)
\leq 
|v^\prime(Y_s)|\big(\alpha_t +\beta|Y_s|+\gamma |Z_s|+ f(|Y_s|)|Z_s|^2\big). \label{itov0}
\end{align}
Recall that
$v^{\prime\prime}(x)-2f(|x|)|v^\prime(x)|=1$ a.e.
Hence (\ref{itov1}), (\ref{itov0}) 
and (\ref{Aneg})
give
\begin{align*}
\frac{1}{2}\int_0^\tau |Z_s|^2 ds \leq v(Y_\tau)
+\int_0^\tau |v^\prime(Y_s)|
\big(\alpha_s +\beta|Y_s| +\gamma |Z_s| \big) ds -\int_0^\tau v^\prime(Y_s)Z_s dW_s.
\end{align*}
Moreover, since $v(x)\leq \frac{M^2x^2}{2}$ and  $|v^\prime(x)|\leq M^2|x|$,   this inequality gives
\begin{align}
\int_0^\tau |Z_s|^2 ds \leq c_1 (Y^*)^2
+c_1\int_0^\tau  |Y_s|
\big(\alpha_s +|Y_s| + |Z_s| \big) ds -2\int_0^\tau v^\prime(Y_s)Z_s dW_s, \label{itov2}
\end{align}
where $c_1:= 2{M^2}(1 \vee \beta\vee \gamma)$.
Note that  in (\ref{itov2}),
\begin{align*}
\int_0^{\tau} |Y_s| \alpha_s ds &\leq   \frac{1}{2}(Y^*)^2 + \frac{1}{2}|\alpha|_{T}^2, \\
c_1 \int_0^{\tau} |Y_s| |Z_s| ds &\leq  \frac{1}{2}c^2_1 T\cdot (Y^*)^2 + \frac{1}{2}\int_0^{\tau} |Z_s|^2 ds.
\end{align*}
Hence (\ref{itov2}) yields
\[
\int_0^{\tau}  |Z_s|^2 ds  \leq (3c_1+c^2_1 T) (Y^*)^2 + c_1|\alpha|_T^2- 4\int_0^{ \tau} v^\prime(Y_s)Z_s dW_s.
\]
This estimate implies that for any $p\geq 1$,
\begin{align}
\mathbb{E}\Big[\Big(\int_0^{\tau}  |Z_s|^2 ds \Big)^{\frac{p}{2}}\Big] \leq c_2 \mathbb{E}\Big[(Y^*)^p + |\alpha|_T^p + \Big|\int_0^{ \tau} v^\prime (Y_s) Z_s dW_s\Big|^{\frac{p}{2}}\Big], \label{lpestimate:1}
\end{align}
where $c_2 := 3^{\frac{p}{2}}\big((3c_1+c^2_1 T) \vee 4\big)^{\frac{p}{2}}$.
Define for each $n\in\mathbb{N}^+$, 
$
\tau_n := \inf \big\{  
t\geq 0:\int_0^t |Z_s|^2 ds
\geq n
\big\} \wedge T.
$
 We then replace $\tau$ by $\tau_n$ and use  Davis-Burkholder-Gundy inequality to obtain
\begin{align*}
c_2\mathbb{E}\Big[ \Big(\int_0^{ \tau_n}v^\prime(Y_s)Z_s dW_s \Big)^\frac{p}{2}\Big]& \leq c_2c({p})M^p
\mathbb{E}\Big[ \Big(\int_0^{ \tau_n} |Y_s|^2 |Z_s|^2 ds\Big)^\frac{p}{4}\Big] \\
& \leq \frac{1}{2}c^2_2c(p)^2 M^{2p}\cdot \mathbb{E} \big[ (Y^*)^p\big]  + \frac{1}{2} \mathbb{E}\Big [\Big(\int_0^{\tau_n}  |Z_s|^2 ds\Big)^{\frac{p}{2}}\Big]
\\
&< +\infty.
\end{align*}
We explain that in this inequality, $c(p)$ denotes the constant in Davis-Burkholder-Gundy inequality which only depends on $p$. With this estimate, 
we  come back to (\ref{lpestimate:1}).
Transferring the quadratic term  to the left-hand side  of (\ref{lpestimate:1}) and using Fatou's lemma, we obtain
\[
 \mathbb{E}\Big [\Big(\int_0^{T}  |Z_s|^2 ds\Big)^{\frac{p}{2}} \Big] \leq c \Big( \mathbb{E} \big[ (Y^*)^p+|\alpha|_T^{p}\big] \Big),
\]
where $c : = c_2^2c(p)^2M^{2p} + 2c_2.$

To estimate  $\int_0^T f(|Y_s|)|Z_s|^2 ds$ we use $u: = u^{2f(|\cdot|)}$. This helps to transfer $\int_0^T f(|Y_s|)|Z_s|^2 ds$ to the left-hand side so that standard estimates can be used. 
The proof is omitted since it is not relevant to our study. 
\qed
\end{Proof}

We continue our study 
by sharpening Lemma \ref{lpestimate1} for $p>1$. 
We follow  Proposition 3.2, Briand et al \cite{B2003}  and extend it to quadratic BSDEs. As an important byproduct, we obtain 
the a priori bound for solutions which is crucial to the construction of a solution. 
\begin{Lemma}[A Priori Estimate (ii)]
\label{lpestimate2}
 Let $p > 1$ and   {\rm \ref{lpas1}}  hold for $(F, \xi)$. If $(Y, Z) \in \mathcal{S}^p \times \mathcal{M} $ is a solution to $(F, \xi)$, then	 
\[
\mathbb{E}\big[(Y^*)^p\big] + \mathbb{E}\Big[ \Big( \int_0^T |Z_s|^2 ds   \Big)^{\frac{p}{2}}\Big]
+ \mathbb{E}\Big[ \Big( \int_0^T f(|Y_s|)|Z_s|^2 ds   \Big)^{p}\Big]
  \leq c \Big(\mathbb{E}\big[|\xi|^p + |\alpha|_T^p\big] \Big).
\]
In particular,  
\[
\mathbb{E}\Big[ \sup_{s \in [t, T]}|Y_s|^p \Big|\mathcal{F}_t \Big] \leq c \mathbb{E} \big[|\xi|^p + |\alpha|_{t, T}^p \big| \mathcal{F}_t\big].
\]
In both cases,  $c$ is a constant  only depending on
$T, M^{f(|\cdot|)}, \beta, \gamma, p$.
\end{Lemma} 
\begin{Proof}
Let $u: =u^{f(|\cdot|)}$ and $M : = M^{f(|\cdot|)}$, and  denote 
 $u(|Y_t|), u^\prime(|Y_t|), u^{\prime\prime}(|Y_t|)$ by $u_t, u_t^\prime, u_t^{\prime\prime}$, respectively.
By Tanaka's formula applied  to $|Y_t|$ and  It\^{o}-Krylov formula  applied to $u_t$,
\begin{align*}
 u_t &= u_T + \int_t^{T}\sgn(Y_s) u^\prime_sF(s, Y_s, Z_s)ds -\frac{1}{2}\int_t^{T} \mathbb{I}_{\{Y_s\neq 0\}} u^{\prime\prime}_s|Z_s|^2 ds \\
&- \int_t^{ T} \sgn (Y_s) u^\prime_s Z_s dW_s - \int_t^{ T}u^{\prime}_s dL_s^0(Y),
\end{align*}
where $L^0(Y)$ is the local time of $Y$ at $0$.  Lemma \ref{lpinequality} applied to $u_t$ then gives 
\begin{align*}
|u_t|^p & + \frac{p(p-1)}{2} \int_t^T 
 \mathbb{I}_{\{ u_s \neq 0  \}} 
\mathbb{I}_{\{Y_s \neq 0 \}}|u_s|^{p-2}|u^\prime_s|^2 |Z_s|^2  ds \\
& = |u_T|^p + p \int_t^T \sgn(u_s)|u_s|^{p-1}  \Big (  \sgn(Y_s)u^\prime_sF(s, Y_s, Z_s) -\frac{1}{2} \mathbb{I}_{\{Y_s\neq 0\}}u^{\prime\prime}_s |Z_s|^2\Big)ds\\ &-p\int_t^T\sgn(u_s)|u_s|^{p-1} u^{\prime}_sdL^0_s(Y) 
 - p\int_t^T  \sgn(u_s)\sgn (Y_s)|u_s|^{p-1}u^\prime_s Z_s dW_s.
\end{align*}
To simplify this equality we recall that $\sgn(u_s)= \mathbb{I}_{\{ u_s \neq 0\}} =\mathbb{I}_{\{ Y_s\neq 0\}}$  and 
$u^{\prime\prime}(x) = 2f(x)u^\prime(x)
$ a.e. Hence
\begin{align*}
|u_t|^p  &+ \frac{p(p-1)}{2} \int_t^T \mathbb{I}_{\{Y_s \neq 0 \}}|u_s|^{p-2}|u^\prime_s|^2 |Z_s|^2  ds \\
& \leq |u_T|^p + p \int_t^T \mathbb{I}_{\{Y_s \neq 0 \}}|u_s|^{p-1} u^\prime_s \big ( \alpha_s + \beta |Y_s| + \gamma |Z_s|\big) ds \\
& - p\int_t^T \sgn (Y_s)|u_s|^{p-1} u^\prime_s Z_s dW_s.
\end{align*}
Let $\{c_n \}_{n\in\mathbb{N}^+}$ be constants to be determined. 
 Since $\frac{|x|}{M}\leq u(|x|)\leq M|x|$ and $\frac{1}{M}\leq  u^\prime(|x|)\leq M$, this inequality yields
\begin{align}
|Y_t|^p  &+ c_1 \int_t^T \mathbb{I}_{\{Y_s \neq 0\}}|Y_s|^{p-2} |Z_s|^2 ds  \nonumber \\
&  \leq M^p |\xi|^p + M^p \int_t^T \mathbb{I}_{\{Y_s \neq 0\}}|Y_s|^{p-1} \big ( \alpha_s + |\beta| |Y_s| + \gamma |Z_s|\big)  ds \nonumber \\
&- p \int_t^T \sgn (Y_s)|u_s|^{p-1}  u^\prime_s Z_s dW_s,\label{up1}
\end{align}
where $c_1 := \frac{p(p-1)}{2M^p} >0$. Observe that in (\ref{up1}),
\[
M^p  \gamma \mathbb{I}_{\{Y_s \neq 0\}}|Y_s|^{p-1}|Z_s| \leq  \frac{M^{2p}\gamma^2}{2c_1}|Y_s|^p + \frac{c_1}{2}\mathbb{I}_{\{Y_s\neq 0\}}|Y_s|^{p-2} |Z_s|^2.
\]
We then use this inequality  to (\ref{up1}).  Set 
 $c_2: = M^p \vee \big( M^p \beta + \frac{M^{2p}\gamma^2}{2c_1}\big)$, 
\[
X:  = c_2 \Big( |\xi|^p + \int_0^T |Y_s|^{p-1} \big ( \alpha_s +  |Y_s| \big)  ds\Big),
\]
and $N$ to be the local martingale part of (\ref{up1}). Hence (\ref{up1}) gives
\begin{align}
|Y_t|^p +\frac{c_1}{2} \int_t^T \mathbb{I}_{\{Y_s \neq 0\}}|Y_s|^{p-2} |Z_s|^2 ds \leq X - N_T+ N_t. \label{lp2}
\end{align}
We claim that $N$ is a martingale. Let $c(1)$ be the constant in Davis-Burkholder-Gundy inquality for $p=1.$ We have 
\begin{align*}
\mathbb{E} \big[N^*\big] &\leq c(1) \mathbb{E} \big[\langle N \rangle _T^{\frac{1}{2}}\big] \\
&\leq c(1)M^p \mathbb{E} \Big[ \Big(\int_0^T  |Y_s|^{2p-2} |Z_s|^2 ds \Big)^{\frac{1}{2}}\Big]  \\
& \leq  \frac{c(1)M^p}{p} \Big((p-1)\mathbb{E}\big[(Y^*)^p\big] +  \mathbb{E}\Big[ \Big( \int_0^T |Z_s|^2 ds \Big)^{\frac{p}{2}}\Big]\Big) \\
&< +\infty,
\end{align*}
where the last two lines come from  Young's inequality and Lemma \ref{lpestimate1} (a priori estimate (i)). Hence $N$ is a martingale. 
 Coming back to (\ref{lp2}), we deduce that 
\begin{align}
\mathbb{E}\Big[  \int_0^T \mathbb{I}_{\{Y_s\neq 0\}}|Y_s|^{p-2} |Z_s|^2 ds \Big] \leq \frac{2}{c_1} \mathbb{E} [X]. \label{2x}
\end{align}
Now we estimate $Y$ via $X$. To this end,  taking supremum over $t\in [0, T]$ and using Davis-Burkholder-Gundy inequality to (\ref{lp2})  give
\begin{align}
\mathbb{E}\big[(Y^*)^p\big] \leq \mathbb{E}[ X] + c(1) \mathbb{E} \big[\langle N \rangle_T^\frac{1}{2}\big]. \label{lp3}
\end{align}
Here $c(1)$ denotes the constant in Davis-Burkholder-Gundy inequality for $p=1$.
The second term  in (\ref{lp3}) yields by Cauchy-Schwartz inequality that 
\begin{align*}
 c(1)\mathbb{E} [\langle N \rangle_T^\frac{1}{2}] & \leq  c(1)M^p \mathbb{E}\Big[ (Y^*)^{\frac{p}{2}} \Big( \int_0^T\mathbb{I}_{\{Y_s \neq 0\}} |Y_s|^{p-2} |Z_s|^2  ds\Big)^{\frac{1}{2}}\Big] \\
 &\leq \frac{1}{2}\mathbb{E} \big[(Y^*)^p\big] + \frac{c(1)^2M^{2p}}{2} \mathbb{E}\Big[\int_0^T \mathbb{I}_{\{Y_s \neq 0\}}|Y_s|^{p-2}|Z_s|^2  ds\Big].
\end{align*}
Using   (\ref{2x}) to this inequality gives the estimate of
$\langle N\rangle^{\frac{1}{2}}$  via $Y$ and $X$.   With this estimate we come back to 
 (\ref{lp3}) and obtain
\[
\mathbb{E}[(Y^*)^p] \leq 2\Big(1+\frac{2c({1})^2M^{2p}}{c_1}\Big) \mathbb{E}[X].
\]
Set $c_3: = 2c_2\big( 1+ \frac{c(1)^2M^{2p}}{2}\big)  .$
This inequality  yields
\begin{align}
\mathbb{E}[(Y^*)^p] \leq c_{3} \Big(\mathbb{E}\big[|\xi|^p \big] + \mathbb{E}\Big[ \int_0^T |Y_s|^{p-1}\alpha_s ds\Big]+ \mathbb{E}\Big[ \int_0^T |Y_s|^p ds \Big]\Big).\label{up2}
\end{align} 
 Young's inequality used to the second term on the right-hand side of this inequality gives
\[
c_{3} \int_0^T |Y_s|^{p-1}  \alpha_s ds \leq \frac{1}{2} (Y^*)^p + \frac{c_3}{p}\Big(\frac{2}{c_3 q}\Big)^{\frac{p}{q}}|\alpha|_T^p,
\]
where $q$ is the conjugate index of $p$.
Set $c_4: = 2\Big(c_3 \vee\frac{c_3}{p}\big(\frac{2}{c_3 q}\big)^{\frac{p}{q}}  \Big)$.  (\ref{up2}) and this inequality then yield
\[
\mathbb{E}\big[(Y^*)^p\big] \leq  c_4\Big( \mathbb{E}\big[|\xi|^p + |\alpha|_T^p\big]+  \mathbb{E}\Big[ \int_0^T \sup_{u\in[0, s]} |Y_u|^p ds \Big]\Big),
\]
By Gronwall's lemma, 
\[
\mathbb{E}\big[(Y^*)^p\big] \leq c_4 \exp(c_4 T) \mathbb{E}\big[|\xi|^p + |\alpha|_T^p\big].
\]
Finally by Lemma \ref{lpestimate1} we conclude that there exists a constant $c$ only depending on $T, M, \beta, \gamma, p$ such that 
\[
\mathbb{E}\big[(Y^*)^p\big] + \mathbb{E}\Big[ \Big( \int_0^T |Z_s|^2 ds   \Big)^{\frac{p}{2}}\Big]
+ \mathbb{E}\Big[ \Big( \int_0^T f(|Y_s|)|Z_s|^2 ds   \Big)^{p}\Big]
  \leq c \mathbb{E}\big[|\xi|^p + |\alpha|_T^p\big] .
\]
To prove the remaining statement, we view any fixed $t\in [0, T]$ as  the initial time, reset
\[
X:  = c_2\Big( |\xi|^p +  \int_t^T |Y_s|^{p-1} \big ( \alpha_s +  |Y_s| \big)  ds\Big)
\]
and replace all estimates by conditional estimates. 
\qed
\end{Proof}

An immediate consequence of Lemma \ref{lpestimate2} is that 
\[
|Y_t| \leq \Big(c \mathbb{E} \big[|\xi|^p + |\alpha|_{T}^p \big| \mathcal{F}_t\big] \Big)^{\frac{1}{p}},
\]
i.e., $Y$  has an a priori bound  which is a  continuous supermartingale. 

With this estimate    
we are  ready to construct a $\mathbb{L}^p(p>1)$ solution
via inf-(sup-)convolution
as in  Briand et al   \cite{BH2006}, \cite{BH2008},  \cite{BLS2007}.
A localization procedure where the a priori bound plays a crucial role is used  and the monotone stability  result takes the limit. 
\begin{thm}[Existence]\label{lpexistence} Let $p > 1$ and {\rm\ref{lpas1}} hold for $(F, \xi)$.
Then there exists a solution to $(F, \xi)$ in $\mathcal{S}^p\times \mathcal{M}^p$.
\end{thm}
\begin{Proof}
We introduce the notations used throughout the proof. Define the process
\[
X_t := \Big(c \mathbb{E}\big[|\xi|^p + |\alpha|_T^p \big|\mathcal{F}_t \big]\Big)^{\frac{1}{p}},
\]
where $c$ is the constant  defined in Lemma \ref{lpestimate2}. Obviously
$X$ is continuous by It\^{o}  representation theorem. Moreover, for  each $m, n \in \mathbb{N}^+$, set 
\begin{align*}
\tau_m &: = \inf \big\{
t\geq 0: |\alpha|_t +X_t \geq m      \big\} \wedge T, \\
\sigma_n &:= \inf \big\{ 
t\geq 0: |\alpha|_t \geq n\big\} \wedge T.
\end{align*}
 It then follows from the continuity of $|\alpha|_\cdot$ and $X$ that $\tau_m$ and $\sigma_n$ increase stationarily to $T$ as $m, n$ goes to $+\infty$, respectively. To apply a double approximation procedure we  define
\begin{align*}
F^{n, k}(t, y, z): &=\mathbb{I}_{\{t\leq \sigma_n \}} \inf_{y^\prime, z^\prime}\big \{ F^+(t, y^\prime, z^\prime) + n|y-y^\prime|+n|z-z^\prime| \big\} \\
&-\mathbb{I}_{\{t\leq \sigma_k \}} \inf_{y^\prime, z^\prime} \big\{ F^-(t, y^\prime, z^\prime) + k|y-y^\prime|+k|z-z^\prime| \big\},
\end{align*}
and $
\xi^{n , k}: = \xi^+ \wedge n - \xi^- \wedge k.$ 

Before proceeding to the proof we give some useful facts.  
By  Lepeltier and San Martin \cite{LS1997}, 
$F^{n, k}$ is  Lipschitz-continuous in $(y, z)$; as  $k$ goes to $+\infty$,  $F^{n, k}$ converges decreasingly 
 uniformly on compact sets to a limit denoted by $F^{n,\infty}$; as $n$ goes to $+\infty$, $F^{n, \infty}$ converges increasingly uniformly on compact sets to $F$. Moreover, $\big||F^{n, k}(\cdot, 0, 0) |\big|_T$ and $\xi^{n, k}$ are bounded. 
 
 Hence,  by  Briand et al \cite{B2003},
there exists a unique solution $(Y^{n, k}, Z^{n, k}) \in \mathcal{S}^p\times\mathcal{M}^{p}$ to $(F^{n, k},  \xi^{n, k})$; by comparison theorem, $Y^{n, k}$ is increasing in $n$ and decreasing in $k$. We are about to take the limit by the  monotone stability result.

However,  $\big||F^{n, k}(\cdot, 0, 0)|\big|_T$ and $Y^{n, k}$
are not uniformly bounded in general. 
To overcome this difficulty, we  use Lemma \ref{lpestimate2}
and  work on random time interval
where $Y^{n,k}$ and $\big||F^{n, k}(\cdot, 0, 0)|\big|_\cdot$  are uniformly bounded.
This is the motivation to 
 introduce $X$ and $\tau_m$. To be more precise, the localization procedure is as follows. 

Note that $(F^{n, k}, \xi^{n, k})$ satisfies  \ref{lpas1} associated with $(\alpha, \beta, \gamma, \varphi, f)$. Hence by  Lemma \ref{lpestimate2} (a priori estimate (ii)),
\begin{align}
|Y_t^{n, k}| &\leq \Big(c \mathbb{E}\big[ |\xi^{n, k}|^p + |\mathbb{I}_{[0, \sigma_n \vee \sigma_k]}\alpha|_T^p \big| \mathcal{F}_t\big]\Big)^{\frac{1}{p}}\nonumber\\
&\leq 
X_t. \label{lpx1}
\end{align}
In view of the definition of $\tau_m$,  we deduce that 
\begin{align}
|Y^{n, k}_{t\wedge \tau_m} |\leq X_{t\wedge \tau_m} \leq m.\label{lpbound}
\end{align}
Hence $Y^{n, k}$ is uniformly bounded on $[0, \tau_m]$. Secondly, 
 given $(Y^{n, k}, Z^{n, k})$ which solves $(F^{n, k}, \xi^{n, k})$, it is  immediate  that 
 $(Y^{n, k}_{\cdot\wedge\tau_m}, \mathbb{I}_{[0, \tau_m]}Z^{n, k})$ solves 
$(\mathbb{I}_{[0, \tau_m]}F^{n, k}, Y^{n, k}_{\tau_m})$. 
To make the monotone stability result adaptable, we 
use a truncation procedure. Define 
\[
\rho(y): = -\mathbb{I}_{\{y<-m\}} m + \mathbb{I}_{\{|y| \leq m\}} y +\mathbb{I}_{\{y>m\}}  m.
\]
Hence from (\ref{lpbound})  $(Y^{n, k}_{\cdot\wedge\tau_m}, \mathbb{I}_{[0, \tau_m]}Z^{n, k})$ meanwhile solves $(\mathbb{I}_{[0, \tau_m]}(t)F^{n, k}(t, \rho(y), z), Y^{n, k}_{{\tau_m}})$. Secondly,  we have
\begin{align*}
|\mathbb{I}_{[0, \tau_m]}(t)F^{ n, k} (t, \rho(y), z)|&\leq \mathbb{I}_{\{t\leq \tau_m\}} \big(\alpha_t +  \varphi(|\rho(y)|) + \gamma |z| + f(|\rho(y)|)|z|^2\big)\\
&\leq\mathbb{I}_{\{t\leq \tau_m\}} \Big( \alpha_t+ \varphi(m) + \gamma |z|  + \sup_{|y|\leq m}f(|\rho(y)|) |z|^2\Big)\\
&\leq \mathbb{I}_{\{t\leq \tau_m\}} \Big(\alpha_t+\varphi (m) + \frac{\gamma^2}{4} +  \big(\sup_{|y|\leq m}f(|\rho(y)|) +1\big) |z|^2\Big),
\end{align*}
where $\sup_{|y|\leq m}f(|\rho(y)|)$ is bounded  for each $m$ due to $f(|\cdot|)\in\mathcal{I}$. Moreover, the definition of $\tau_m$ implies $|\alpha|_{\tau_m} \leq m$.
Hence we can use the monotone stability result (Kobylanksi \cite{K2000}, Briand and  Hu \cite{BH2008} or Theorem \ref{monotonestability}) 
to obtain   $(Y^{m, n, \infty}, Z^{ m, n, \infty})\in\mathcal{S}^\infty\times\mathcal{M}^2$ which solves $(\mathbb{I}_{[0, \tau_m]}(t)F^{n, \infty}(t, \rho(y), z), \inf_k Y_{\tau_m}^{n, k})$. Moreover,   $Y^{m, n, \infty}_{\cdot \wedge \tau_m}$ is the $\mathbb{P}$-a.s. uniform limit of $Y^{n, k}_{\cdot \wedge \tau_m}$ as $k$ goes to $+\infty$.
These arguments hold for any $m, n\in\mathbb{N}^+$.  

Due to this convergence result we can pass the comparison property to 
 $Y^{m, n, \infty}$. We use the monotone stability result again to the sequence indexed by $n$ to obtain $(\widetilde{Y}^m, \widetilde{Z}^m) \in \mathcal{S}^\infty\times\mathcal{M}^2$ which solves $(\mathbb{I}_{[0, \tau_m]}(t)F(t, \rho(y), z), \ 
\sup_n\inf_k Y^{n, k}_{\tau_m})$.   Likewise, $\widetilde{Y}^m_{\cdot}$ is the $\mathbb{P}$-a.s. uniform limit of $Y^{m, n, \infty}_{\cdot}$  as $n$ goes to $+\infty$. Hence we  obtain from (\ref{lpbound}) that $|\widetilde{Y}^m_{t}| \leq X_{t\wedge \tau_m} \leq m $. Therefore, 
$(\widetilde{Y}^m, \widetilde{Z}^m)$  solves $(\mathbb{I}_{[0, \tau_m]}F, \sup_n\inf_k Y_{\tau_m}^{n, k})$, i.e., 
\begin{align}
\widetilde{Y}^m_{t\wedge \tau_m}= \sup_n\inf_k {Y}_{\tau_m}^{n, k} +
\int_{t\wedge \tau_m}^{\tau_m} F(s, \widetilde{Y}^m_s, \widetilde{Z}_s^m)ds -
\int_{t\wedge \tau_m}^{\tau_m}  \widetilde{Z}_s^m dW_s. \label{locbsde}
\end{align}

We recall that the monotone stability result also implies that 
 $\widetilde{Z}^m$ is the $\mathcal{M}^2$-limit of $\mathbb{I}_{[0, \tau_m]}Z^{n, k}$ as $k, n$ goes to $+\infty$. This fact and previous convergence results give
 \begin{align}
 \widetilde{Y}^{m+1}_{\cdot\wedge \tau_m}  &=
  \widetilde{Y}^{m}_{\cdot\wedge \tau_m} \ \mathbb{P}\text{-a.s.}, \nonumber\\
 \mathbb{I}_{\{t\leq \tau_m\}}\widetilde{Z}^{m+1}_t &= \mathbb{I}_{\{t\leq \tau_m\}}\widetilde{Z}^m_t \ dt\otimes d\mathbb{P}\text{-a.e.} \label{piece}
 \end{align} 
Define $(Y, Z)$ on $[0, T]$ by
\begin{align*}
Y_t := \mathbb{I}_{\{t\leq \tau_1\}} \widetilde{Y}_t^1  + \sum_{m\geq 2} \mathbb{I}_{]\tau_{m-1}, \tau_m]}\widetilde{Y}^m_t, \\
Z_t := \mathbb{I}_{\{t\leq \tau_1\}} \widetilde{Z}_t^1  + \sum_{m\geq 2} \mathbb{I}_{]\tau_{m-1}, \tau_m]}\widetilde{Z}^m_t.
\end{align*}
By (\ref{piece}),  we have $Y_{\cdot \wedge \tau_m} = \widetilde{Y}^m_{\cdot\wedge \tau_m}$
and
$\mathbb{I}_{\{t\leq \tau_m\}}Z_t = \mathbb{I}_{\{t\leq \tau_m\}} \widetilde{Z}_t^m$. Hence we can rewrite  
 (\ref{locbsde}) as
\[
Y_{t\wedge \tau_m} = \sup_n\inf_k Y^{n, k}_{\tau_m} +\int_{t\wedge \tau_m}^{\tau_m}  F{(s, Y_s, Z_s )}ds -\int_{t\wedge \tau_m}^{\tau_m} Z_s dW_s.
\]
 By sending $m$ to $+\infty$, we deduce  that $(Y, Z)$ solves $(F, \xi)$.  Since $(Y^{n, k}, Z^{n, k})$ verifies Lemma \ref{lpestimate2}, we use Fatou's lemma to prove that $(Y, Z)\in\mathcal{S}^p\times\mathcal{M}^p$.
\qed
\end{Proof}

Theorem \ref{lpexistence} proves the existence of a $\mathbb{L}^p(p>1)$ solution under  \ref{lpas1} which to our knowledge the most general asssumption. 
For example, 
 \ref{lpas1}(ii) allows one to get rid of monotonicity in $y$ which is required by, e.g.,  Pardoux \cite{P1999} and Briand et al \cite{BC2000}, \cite{B2003}, \cite{BLS2007}. Meanwhile, in contrast to these works, 
the generator can also be quadratic by setting $f(|\cdot|)\in\mathcal{I}$. 
Hence Theorem \ref{lpexistence} provides a unified way to construct solutions to both non-quadratic and quadratic BSDEs via the monotone stability result.

 On the other hand, Theorem \ref{lpexistence} is an extension of  Bahlali et al  \cite{BEO2014} which only studies BSDEs with $\mathbb{L}^2$ integrability and linear-quadratic growth. However, in contrast to their work,  \ref{lpas1} is not sufficient in our setting to ensure the existence of a maximal or minimal solution, since the  double approximation procedure  makes the comparison between solutions impossible.
 
However, to prove the existence of a maximal or minimal solution is no way impossible. Since we have $X$ as the a priori bound for solutions, 
we can convert the question of existence into the question of existence  for quadratic BSDEs with double  barriers. This problem has been solved by introducing the notion of generalized BSDEs; see Essaky and Hassani \cite{E2013}.
 
\begin{remm}
One may ask that as we use a localization procedure, whether $f$ being   bounded or  integrable  only on compact subsets of 
$\mathbb{R}$ 
rather than of class $\mathcal{I}$ sufficies to ensure  the existence result. It turns out that in data in  $\mathbb{L}^p$  is not sufficient for such a generalization, and exponential moments integrability is required.  
Hence, our existence result shall be seen as merely complementary to  the quadratic BSDEs studied by
 Briand and  Hu \cite{BH2006}, \cite{BH2008} rather than a complete generalization. 

 Below is an illustrating example  with $f =  1$ which clearly doesn't belong to $\mathcal{I}$. Similar version can be found in Briand et al \cite{BLS2007}.
\end{remm}
\begin{ex}\label{ex1}

 There exists a solution $(Y, Z)$ in $\mathcal{S}^2\times \mathcal{M}^2$ 
to 
\begin{align}
Y_t = \xi + \int_t^T |Z_s|^2 ds  - \int_t^T Z_s dW_s\label{bsde2}
\end{align} 
 if and only if \[\mathbb{E}\big[ \exp(2\xi) \big] < +\infty.\]
\end{ex}
\begin{Proof} (i). $\Longrightarrow$.
Let $(Y, Z) \in \mathcal{S}^2 \times \mathcal{M}^2$ be a solution to (\ref{bsde2}). By It\^{o}'s formula, 
\[
\exp(2 Y_t) = \exp (2 Y_0) + \int_0^t \exp(2Y_s)Z_s dW_s.
\]
Now we define 
$
\tau_n : = \inf \{t\geq 0: Y_t \geq n\}  
$
for each $n\in\mathbb{N}^+$.
$\mathcal{F}_0$ being trivial implies that $Y_0$ is  a constant.
 Hence $\int_0^{\cdot\wedge \tau_n} \exp (2Y_s) Z_s dW_s$ is a bounded martingale, and 
\[
\mathbb{E}\big[\exp(2Y_{T\wedge\tau_n})\big] = \mathbb{E}\big[\exp(2Y_0)\big].
\]
By Fatou's lemma we obtain $
\mathbb{E} [\exp (2 \xi)]  <+\infty. $

(ii). $\Longleftarrow.$ Assume  $\mathbb{E}\big[\exp (2 \xi)\big]  <+\infty$. 
 Thanks to It\^o representation theorem, we can define  $(\widetilde{Y}, \widetilde{Z})\in\mathcal{S}\times\mathcal{M}$ by 
\[
\widetilde{Y}_t : = \mathbb{E}\big[\exp (2\xi)\big| \mathcal{F}_t \big] = \widetilde{Y}_0 + \int_0^t \widetilde{Z}_s dW_s.
\]
Set $
(Y, Z) := ( \frac{1}{2} \ln \widetilde{Y}, \frac{\widetilde{Z}}{2\widetilde{Y}}).
$
By It\^{o}'s formula applied to ${Y}$, we easily deduce that  $(Y, Z)$ solves (\ref{bsde2}). It thus remains to prove $(Y, Z)\in\mathcal{S}^2\times \mathcal{M}^2$. Since $x\mapsto\ln (x)$ is concave and increasing, 	Jensen's inequality yields
\begin{align}
Y_t = \frac{1}{2}\ln \big(\mathbb{E}\big[\exp(2\xi)\big|\mathcal{F}_t\big]\big) \geq \mathbb{E}\big[\xi \big| \mathcal{F}_t\big]\geq 0. \label{positive}
\end{align}
Hence $Y$ is nonnegative. 
For each $n\in\mathbb{N}^+$,  define $\tau_n : = \inf \big\{ t\geq:   \int_0^t |Z_s|^2 ds \geq n    \big\}$. 
$(Y, Z)$ being a solution to
 (\ref{bsde2}) implies that 
\begin{align*}
\int_0^{T\wedge \tau_n} |Z_s|^2 ds &= Y_0 - Y_{T\wedge \tau_n} + \int_0^{T\wedge \tau_n} Z_s dW_s \\
& \leq Y_0  +  \int_0^{T\wedge \tau_n} Z_s dW_s.
\end{align*}
Hence (\ref{positive}) gives
\begin{align}
\mathbb{E}\Big[ \Big(\int_0^{T\wedge \tau_n} |Z_s|^2 ds \Big)^2 \Big] &\leq 2 Y_0^2 + 2\mathbb{E}\Big[  \Big(\int_0^{T\wedge \tau_n} Z_s dW_s \Big)^2    \Big]. \label{z2}
\end{align}
 Moreover, by Jensen's inequality applied to the left-hand side of (\ref{z2}), 
\begin{align*}
\mathbb{E}\Big[  \int_0^{T\wedge \tau_n} |Z_s|^2  ds    \Big]^2 \leq 2Y_0^2 + 2\mathbb{E}\Big[  \int_0^{T\wedge \tau_n} |Z_s|^2 ds\Big], 
\end{align*}
Using $2a \leq \frac{a^2}{2} + 2$ to the last term of this inequality gives
\[
\mathbb{E}\Big[  \int_0^{T\wedge \tau_n} |Z_s|^2 ds\Big] < 4Y_0^2 + 4.
\]
Hence, Fatou's lemma yields $Z\in \mathcal{M}^2$. We then use this result and  Fatou's lemma to (\ref{z2}) to obtain
\[
\mathbb{E}\Big[ \Big(\int_0^{T} |Z_s|^2 ds \Big)^2 \Big] < +\infty.
\]
Finally we deduce from  (\ref{bsde2}) that
\begin{align*}
\mathbb{E}\big[(Y^*)^2\big]&\leq 3\mathbb{E}\big[|\xi|^2\big] + 3\mathbb{E}\Big[\Big(\int_0^T |Z_s|^2 ds \Big)^2\Big] + 3 \mathbb{E}\Big[\Big(\Big|\int_\cdot^T Z_s dW_s\Big|^*\Big)^2\Big] \\
& <+\infty.
\end{align*}
Hence $(Y, Z) \in \mathcal{S}^2\times \mathcal{M}^2$.
\qed
\end{Proof}

Let us turn to the uniqueness result. Motivated by Briand and  Hu \cite{BH2008} or  Da Lio and  Ley \cite{LL2006} from the point of view of PDEs, we impose a convexity condition so as to use $\theta$-techinique which proves to be convenient to treat quadratic generators.  
We start from  comparison theorem and then move to uniqueness and stability result.
To this end, the following assumptions on $(F, \xi)$ are needed. 
\begin{as}\label{lpas3}  Let $p > 1.$
There exist $\beta_1, \beta_2 \in \mathbb{R}$, $\gamma_1, \gamma_2 \geq 0$, an $\mathbb{R}^+$-valued $\Prog$-measurable process $\alpha$,  a continuous nondecreasing function $\varphi: \mathbb{R}^+ \rightarrow \mathbb{R}^+$ with $\varphi(0)=0$, 
$f(|\cdot|)\in\mathcal{I}$  and $F_1, F_2: \Omega\times[0, T]\times \mathbb{R}\times\mathbb{R}^d\rightarrow \mathbb{R}$ 
which are $\Prog\otimes\mathcal{B}(\mathbb{R})\otimes\mathcal{B}(\mathbb{R}^d)$-measurable
such that $F= F_1+ F_2$, 
$|\xi| + |\alpha|_T\in\mathbb{L}^p$
and  $\mathbb{P}$-a.s.
 \begin{enumerate}
 \item [(i)] for any $t\in[0, T]$, 
$(y, z)\longmapsto F(t, y, z)$ is continuous; 
\item[(ii)] $F_1(t, y, z)$ is  monotonic in $y$ and Lipschitz-continuous in $z$,  and $F_2(t, y, z)$ is monotonic at $y=0$ and of linear-quadratic growth in $z$, 
 i.e., for any  $t\in [0, T], y, y^\prime \in \mathbb{R}, z, z^\prime \in \mathbb{R}^d$,
\begin{align*}
\sgn (y -y ^\prime)\big(F_1(t, y, z)- F_1(t, y^\prime, z)\big)&\leq \beta_1 |y- y^\prime|, \\
 \big|F_1(t, y, z)- F_1(t, y, z^\prime)\big| &\leq  \gamma_1 |z-z^\prime|, \\
 \sgn(y)F_2(t, y, z) &\leq \beta_2|y|+ \gamma_2 |z| + f(|y|)|z|^2;
\end{align*}
\item [(iii)] for any $t\in[0, T]$, $(y, z)\longmapsto F_2(t, y, z)$  is convex;
\item [(iv)] for any $(t, y, z)\in [0, T]\times\mathbb{R}\times\mathbb{R}^d$, 
\begin{align*}
&|F(t, y, z)|\leq \alpha_t + \varphi(|y|) + (\gamma_1 +\gamma_2)|z|+ f(|y|)|z|^2.
  \end{align*}
\end{enumerate}
\end{as}

Intuitively, \ref{lpas3} specifies an additive structure consisting of two  classes of BSDEs. 
The cases where $F_2 = 0$  coincide with 
classic existence and uniqueness results (see, e.g.,  Pardoux \cite{PP1990} or Briand et al \cite{BC2000}, \cite{B2003}). 
When $F_1 =0$,  the BSDEs include those studied by
 Bahlali et al \cite{BEK2013}. 
Given  convexity as an additional requirement, we can prove an existence and uniqueness result in the presence of both components. This can be seen as a general version of  the additive structure discussed in Section \ref{degenerate} and 
a complement to  the quadratic BSDEs studied 
by Bahlali et al \cite{BEK2013} and
Briand and  Hu \cite{BH2008}.

We start our proof of comparison theorem by observing that  \ref{lpas3} implies \ref{lpas1}. Hence the  existence of a
$\mathbb{L}^p(p>1)$  solution is ensured. 

\begin{thm}[Comparison]\label{lpcompare} Let $p > 1$,
and  $(Y, Z), (Y^\prime, Z^\prime) \in \mathcal{S}^p\times\mathcal{M}$ be solutions to $(F, \xi)$, $(F^\prime, \xi^\prime)$, respectively. If  $\mathbb{P}$-a.s. for any $(t ,y ,z)\in [0, T]\times\mathbb{R}\times\mathbb{R}^d$, 
$F(t, y, z)\leq F^\prime(t, y, z)$, $\xi \leq \xi^\prime$ and $F$ verifies {\rm\ref{lpas3}},
then $\mathbb{P}$-a.s.
$Y_\cdot\leq Y^\prime_\cdot$.
\end{thm}
\begin{Proof}
We introduce the  notations used throughout the proof. 
For any $\theta \in (0, 1)$,
define
\begin{align*}
\delta F_t &:= F(t, Y_t^\prime, Z_t^\prime) - F^\prime(t, Y_t^\prime, Z_t^\prime), \\
\delta_\theta Y&: = Y- \theta Y^\prime,  \\
\delta Y &: = Y- Y^\prime, 
\end{align*} 
and  $\delta_\theta Z, \delta Z$, etc. analogously.
$\theta$-technique applied to the generators yields
\begin{align}
& F(t, Y_t, Z_t) - \theta F^\prime(t, Y_t^\prime, Z_t^\prime) \nonumber\\
&=\big(F(t, Y_t, Z_t) - \theta F(t, Y^\prime_t, Z^\prime_t)\big)+\theta\big(F(t, Y^\prime_t, Z^\prime_t) -F^\prime(t, Y_t^\prime, Z_t^\prime)\big)  \nonumber\\
&= \theta \delta F_t + \big(F(t, Y_t, Z_t) -\theta F(t,  Y_t^\prime, Z_t^\prime)\big) \nonumber \\
& = \theta \delta F_t + \big(F_1(t, Y_t, Z_t\big) -\theta F_1(t, Y_t^\prime, Z_t^\prime)\big) + \big(F_2(t, Y_t, Z_t) -\theta F_2(t, Y_t^\prime, Z_t^\prime)\big).
 \label{lpconvex1}
\end{align}
By  \ref{lpas3}(iii),
\begin{align*}
F_2(t, Y_t, Z_t) &= F_2(t, \theta Y_t^\prime + (1-\theta)\frac{\delta_\theta Y_t}{1-\theta}, \theta Z_t^\prime + (1-\theta)\frac{\delta_\theta Z_t}{1-\theta})\nonumber\\
& \leq \theta F_2(t, Y_t^\prime, Z_t^\prime) + (1-\theta) F_2(t, \frac{\delta_\theta Y_t}{1-\theta}, \frac{\delta_\theta Z_t}{1-\theta}). 
\end{align*}
Hence we have 
\begin{align}
F_2(t, Y_t, Z_t) -  \theta F_2(t, Y_t^\prime, Z_t^\prime) \leq (1-\theta) F_2(t, \frac{\delta_\theta Y_t}{1-\theta}, \frac{\delta_\theta Z_t}{1-\theta}).
\label{lpconvex2}
\end{align}
Let $u$  be the function defined in Section \ref{I} associated with a function of class $\mathcal{I}$ to be determined later. 
 Denote  $u((\delta_\theta Y_t)^+), u^\prime((\delta_\theta Y_t )^+), u^{\prime\prime}((\delta_\theta Y_t)^+)$ by
 $u_t, u_t^\prime$, $u_t^{\prime\prime}$, respectively. It is then known from Section \ref{I} that  $u_t \geq 0$ and $u^\prime_t >0$. 
 For any $\tau \in \mathcal{T}$, 
 Tanaka's formula applied  to $(\delta_\theta Y)^+$, It\^{o}-Krylov formula applied  to $u((\delta_\theta Y_t)^+)$ and Lemma \ref{lpinequality} give
\begin{align}
|u_{t\wedge \tau}|^p  &+ \frac{p(p-1)}{2} \int_{t\wedge \tau}^\tau \mathbb{I}_{\{ \delta_\theta Y_s>0\}}|u_s|^{p-2}|u_s^\prime|^2  |\delta_\theta Z_s|^2 ds \nonumber\\
& \leq |u_\tau |^p + p \int_{t\wedge \tau}^\tau \mathbb{I}_{\{ \delta_\theta Y_s > 0\}}|u_s|^{p-1} \underbrace{\Big(u^\prime_s\big ( F(s, Y_s, Z_s)- \theta F^\prime (s, Y^\prime_s, Z^\prime_s)\big) -\frac{1}{2}u^{\prime\prime}_s |\delta_\theta Z_s|^2 \Big)}_{: =\Delta_s}ds \nonumber \\
& 
- p\int_{t\wedge \tau}^\tau \mathbb{I}_{\{\delta_\theta Y_s >0\}}|u_s|^{p-1} u^\prime_s  \delta_\theta Z_s  dW_s. \label{lpconvex3}
\end{align}
By (\ref{lpconvex1}), (\ref{lpconvex2}),  \ref{lpas3}(ii) and   $\delta F \leq 0$,  we deduce that, on  $\{ \delta_\theta Y_s >0\}$, 
\begin{align*}
 \Delta_s
&\leq u^\prime_s \Big(  F_1 (s, Y_s, Z_s) - \theta F_1 (s, Y_s^\prime, Z_s^\prime) +\beta_2 (\delta_\theta Y_s)^+ + \gamma_2 |\delta_\theta Z_s| +\frac{f(\frac{|\delta_\theta Y_s|}{1-\theta})}{1-\theta}|\delta_\theta Z_s|^2\Big)
-\frac{1}{2}u_s^{\prime\prime} |\delta_\theta Z_s|^2.
\end{align*}  To eliminate the quadratic term, 
we associate $u$ with 
$\frac{f(\frac{|\cdot|}{1-\theta})}{1-\theta}$, i.e.,  
\begin{align*}
u(x): &= \int_0^x \exp\Big (2 \int_0^y \frac{f(\frac{|u|}{1-\theta})}{1-\theta}du\Big)dy \\
&= \int_0^x \exp\Big (2 \int_0^{\frac{y}{1-\theta}}f(|u|)du\Big)dy. 
\end{align*}
Hence,  on $\{ \delta_\theta Y_s >0\}$, the above inequality  gives
\begin{align}
 \Delta_s
 \leq  u^\prime_s \big(
  F_1 (s, Y_s, Z_s) - \theta F_1 (s, Y_s^\prime, Z_s^\prime) +\beta_2 (\delta_\theta Y_s)^+ + \gamma_2 |\delta_\theta Z_s|\big).
  \label{lpconvex4}
\end{align}
We are about to send $\theta$ to $1$, and to this end we give some auxiliary facts. Reset
$
M: = \exp \big(2\int_0^\infty f(u)du\big).
$
Obviously
 $1\leq M< +\infty$.
By dominated convergence, for $x\geq 0,$ we have
\begin{align}
&\lim_{\theta\rightarrow 1} u(x)= M x, \nonumber\\
&\lim_{\theta\rightarrow 1} u^\prime(x) = M \mathbb{I}_{\{ x>0\}} + \mathbb{I}_{\{ x= 0\}}. \label{lim}
\end{align}
Taking (\ref{lpconvex4}) and (\ref{lim}) into account, we  come back to  (\ref{lpconvex3}) and 
 send $\theta$ to $1$. 
  Fatou's lemma used to the $ds$-integral on the left-hand side of  (\ref{lpconvex3}) and 
dominated convergence used to the rest integrals give
\begin{align}
&((\delta Y_{t\wedge \tau})^+)^p + 
\frac{p(p-1)}{2}  \int_{t\wedge \tau}^\tau \mathbb{I}_{\{\delta Y_s >0\}} ((\delta Y_s)^+)^{p-2}|\delta Z_s|^2ds \nonumber\\
  &\leq   ((\delta Y_\tau)^+)^p + p \int_{t\wedge \tau}^\tau \mathbb{I}_{\{ \delta Y_s >0\}}((\delta Y_s)^+)^{p-1} ( 
F_1(s, Y_s, Z_s) -F_1(s, Y_s^\prime, Z_s^\prime) +   
 \beta_2 (\delta Y_s)^+ + \gamma_2 |\delta Z_s|     )  ds \nonumber \\
&
-p\int_{t\wedge \tau}^\tau \mathbb{I}_{\{\delta Y_s >0\}} ((\delta Y_s)^+)^{p-1}\delta Z_s  dW_s. \label{lpconvex5}
\end{align}
 Moreover,  
\ref{lpas3}(ii) implies 
\[
\mathbb{I}_{\{ \delta Y_s > 0\}}\big(F_1(s, Y_s, Z_s) - F_1(s, Y_s^\prime, Z_s^\prime) \big)\leq 
\mathbb{I}_{\{ \delta Y_s > 0\}}\big(\beta_1 (\delta Y_s)^+ + \gamma_1 |\delta Z_s|\big).
\]
We then use this inequality to (\ref{lpconvex5}). To eliminate the local martingale, 
we replace $\tau$ by a localization sequence $\{\tau_n\}_{n\in\mathbb{N}^+}$
 and use the same estimation as  in  Lemma \ref{lpestimate2} (a priori estimate (ii)). 
\[
 ((\delta Y_{t\wedge \tau_n})^+ )^p 
  \leq c\mathbb{E}  \big[ ((\delta Y_{\tau_n})^+ )^p \big| \mathcal{F}_t\big] ,
\]
where $c$ is a constant  only depending on $T, \beta_1, \beta_2, \gamma_1, \gamma_2, p$. Since $Y, Y^\prime \in \mathcal{S}^p$ and  $\mathbb{P}$-a.s. $\xi \leq \xi^\prime$, dominated convergence  yields  $\mathbb{P}$-a.s. $Y_t\leq Y_t^\prime$.
Finally by the continuity of $Y$ and  $Y^\prime$ we conclude that  $\mathbb{P}$-a.s. $Y_\cdot \leq Y^\prime_\cdot$.
\qed
\end{Proof}

As a byproduct, we obtain the following existence and uniqueness result.
\begin{Corollary}[Uniqueness]\label{lpunique} Let {\rm\ref{lpas3}} hold for $(F, \xi)$. Then there exists a unique solution in $\mathcal{S}^p\times \mathcal{M}^p$.
\end{Corollary}
\begin{Proof}
\ref{lpas3} implies \ref{lpas1}. Hence existence result holds. The uniqueness is immediate from   Theorem \ref{lpcompare} (comparison theorem).
\qed
\end{Proof}

It turns out that a 
 stability result also holds given the convexity condition. We  denote $(F, \xi)$ satisfying \ref{lpas3} by
$(F, F_1, F_2, \xi).$ We set $\mathbb{N}^0:= \mathbb{N}^+ \cup \{ 0 \}.$
\begin{prop}[Stability]
\label{lpstability}Let $p> 1$. Let
$(F^n, F_1^n, F_2^n, \xi^n)_{n\in\mathbb{N}^0}$ satisfy {\rm\ref{lpas3}}  associated with 
$(\alpha^n, \beta_1, \beta_2, \gamma_1, \gamma_2, \varphi, f)$, and $(Y^n, Z^n)$ be their unique solutions in $\mathcal{S}^p\times\mathcal{M}^p$, respectively. 
If 
$
\xi^n - \xi {\longrightarrow} 0$ and $\int_0^T |F^n -F^0 |(s, Y_s^0, Z_s^0) ds {\longrightarrow} 0$
in $\mathbb{L}^p$
as $n$ goes to $+\infty$, 
then $(Y^n, Z^n)$ converges to $(Y, Z)$ in $\mathcal{S}^p\times\mathcal{M}^p$.
\end{prop}
\begin{Proof}
We prove the stability result in the spirit of  Theorem \ref{lpcompare} (comparison theorem).
 For any $\theta \in (0, 1)$,   define 
\begin{align*}
\delta F_t^n &:= F^0(t, Y^0_t, Z^0_t) - F^n(t, Y^0_t, Z^0_t), \\
\delta_\theta Y^n&: = Y^0-  \theta Y^n, \\
\delta Y^n&: = Y^0-   Y^n,
\end{align*} 
and  $\delta_\theta Z^n, \delta Z^n$, etc. analogously.
We  observe the $\theta$-difference of the generators. Likewise, \ref{lpas3}(iii) implies that 
\begin{align*}
&F^0(t, Y^0_t, Z^0_t) - \theta F^n (t, Y^n_t, Z^n_t)\\
&= \delta F_t^n + \big(  F^{n}(t, Y^0_t, Z^0_t) -\theta F^n(t, Y^n_t, Z^n_t)    \big)\\
&\leq 
\delta F_t^n + \big(F_1^n (t, Y_t^0, Z_t^0) - \theta F_1^n(t, Y_t^n, Z_t^n)\big) +
(1-\theta) F^n_2 (t, \frac{\delta_\theta Y_s^n}{1-\theta}, \frac{\delta_\theta Z_s^n}{1-\theta}).
\end{align*}
We first prove convergence of $Y^n$ and later use it to show  that $Z^n$ also converges.

(i).  
By exactly the same arguments as in Theorem 
\ref{lpcompare} but keeping $\delta F_t^n$ along the deductions,
we obtain
\begin{align}
&((\delta Y_t^n)^+)^p + 
\frac{p(p-1)}{2} \int_t^T\mathbb{I}_{\{ \delta Y_s^n >0\}}  ((\delta Y_s^n)^+)^{p-2}|\delta Z_s^n|^2ds \nonumber\\
  &\leq  ((\delta\xi^n)^+)^p+  p \int_t^T\mathbb{I}_{\{\delta Y_s^n >0\}}  ((\delta Y_s^n)^+)^{p-1} \big(  |\delta F_s^n| + (\beta_1 + \beta_2)(\delta Y_s^n)^+ + (\gamma_1+\gamma_2)|\delta Z_s^n|     \big) ds \nonumber\\
&
-p\int_t^T \mathbb{I}_{\{\delta Y_s^n >0\}}((\delta Y_s^n)^+)^{p-1}\delta Z_s^n  dW_s, \label{convex4}
\end{align}
 By the same way of estimation as in  Lemma \ref{lpestimate2} (a priori estimate (ii)),  we obtain
\[
\mathbb{E}\big[\big(((\delta Y^n)^+)^*\big)^p\big] \leq c\Big( \mathbb{E} \big[\big((\delta \xi^n)^+\big)^p\big] +  \mathbb{E}\big[ \big||\delta F^n_\cdot|\big|_T^p \big] \Big),
\] 
where $c$ is a constant only depending on
$T, \beta_1, \beta_2, \gamma_1, \gamma_2, p$.
Interchanging $Y^0$ and $Y^n$ and analogous deductions then yield
\[
\mathbb{E}\big[\big(((-\delta Y^n)^+)^*\big)^p\big]  \leq c\Big( \mathbb{E} \big[\big((-\delta \xi^n)^+\big)^p\big] +  \mathbb{E}\big[ \big||\delta F^n_\cdot|\big|_T^p \big] \Big).
\] 
Hence a combination of the two inequalities implies the convergence of $Y^n.$

(ii).
To prove the convergence of $Z^n$, we combine the arguments in
Lemma \ref{lpestimate1} (a priori estimate (i)) and Theorem \ref{lpcompare}. 
To this end, we introduce the function $v$ defined in Section \ref{I} associated with 
a function of class $\mathcal{I}$ to be determined later. By It\^o-Krylov formula,	
\begin{align}
v(\delta_\theta Y_0^n) &= v(\delta_\theta \xi^n)
+\int_0^T  v^\prime(\delta_\theta Y_s^n) \big(F^0(s, Y_s^0, Z_s^0) -\theta F^n(s, Y_s^n, Z_s^n) \big) ds
\nonumber\\ 
&-\frac{1}{2}\int_0^T v^{\prime\prime}(\delta_\theta Y_s^n)|\delta_\theta Z_s^n|^2 ds 
- \int_0^T v^\prime(\delta_\theta Y_s^n)\delta_\theta Z_s^n dW_s. \label{lpstability1}
\end{align} Note that \ref{lpas3}(ii)(iii) and 
$v^\prime(\delta_\theta Y_s^n) = \sgn(\delta_\theta Y_s^n)|v^\prime(\delta_\theta Y_s^n)|$ give 
\begin{align}
v^\prime(\delta_\theta Y_s^n) \big(F^0(s, Y_s^0, Z_s^0) &-\theta F^n(s, Y_s^n, Z_s^n) \big) \nonumber\\
&\leq
|v^\prime(\delta_\theta Y_s^n)||\delta F_s^n| \nonumber\\
&+ |v^\prime(\delta_\theta Y_s^n)|\sgn(\delta_\theta Y_s^n)\big(
F^n_1(s, Y_s^0, Z_s^0) - \theta F^n_1(s, Y_s^n, Z_s^n)
\big)
\nonumber\\
&+ |v^\prime(\delta_\theta Y_s^n)|
\Big(  \beta_2 |\delta_\theta Y_s^n| 
+\gamma_2 |\delta_\theta Z_s^n|
+ \frac{f(\frac{|\delta_\theta Y_s^n|}{1-\theta})}{1-\theta} |\delta_\theta Z_s^n|^2
\Big). \label{vf}
\end{align}
We associate $v$ with $\frac{f(\frac{|\cdot|}{1-\theta})}{1-\theta}$ so as to eliminate the quadratic term. Note that
\begin{align}
&\lim_{\theta \rightarrow 1} v(x) =\frac{1}{2}|x|^2, 
\nonumber\\
&\lim_{\theta \rightarrow 1} v^\prime(x) =x. \label{limv}
\end{align}
With (\ref{vf}), (\ref{limv}) and \ref{lpas3}(ii), we come back to (\ref{lpstability1})  and send $\theta$ to $1$.   This gives
\begin{align*}
\frac{1}{2}\int_0^T |\delta Z_s^n|^2 ds
&\leq \frac{1}{2}|\delta \xi^n|^2
+\int_0^T |\delta Y_s^n| \big(|\delta F_s^n| + 
(|\beta_1|+|\beta_2|) |\delta Y_s^n| + (\gamma_1 +\gamma_2)|\delta Z_s^n|
 \big) ds\\
 &- \int_0^T  \delta Y_s^n\delta Z_s^n dW_s.
\end{align*}
Now we use the same way of estimation as in Lemma \ref{lpestimate1} to obtain
\begin{align*}
\mathbb{E}\Big[\Big( \int_0^T  |\delta Z_s^n|^2ds\Big)^{\frac{p}{2}}\Big]
\leq 
c \mathbb{E}\big[
((\delta Y^n)^*)^p +
\big||\delta F^n_\cdot|\big|_T^p
\big
],
\end{align*}
where $c$ is a constant only depending on $T, \beta_1, \beta_2, \gamma_1, \gamma_2, p$.
The convergence of $Z^n$ is then immediate from (i).
\qed
\end{Proof}
\begin{remm}
So far we have obtained the  existence and uniqueness of  a  $\mathbb{L}^p(p > 1)$ solution. The solvability for $p =1$ is not included due to the failure of  Lemma \ref{lpestimate2} (a priori estimate (ii)). One may overcome this difficulty by imposing additional structure conditions as in  Briand et al \cite{B2003}, \cite{BH2006}. To save pages the analysis of $\mathbb{L}^1$ solutions is hence omitted. 
\end{remm}

\section{Applications to Quadratic PDEs}
\label{qpde}
In this section, we give an application of our results  to quadratic PDEs. More precisely, we prove the probablistic representation for  the nonlinear  Feymann-Kac formula associated with the BSDEs in our study. Let us consider the following semilinear PDE
\begin{align}
&\partial_t u(t, x) +\mathcal{L}u(t, x) + 
F(t, x, u(t, x), \sigma^\top \nabla_x u(t, x)) =0, \nonumber\\
&u(T, \cdot) = g, \label{pde}
\end{align}
where $\mathcal{L}$ is the infinitesimal generator of the solution  $X^{t_0, x_0}$ to the Markovian SDE
\begin{align}
X_t &= x_0 +\int_{t_0}^t b(s, X_s)ds + \int_{t_0}^t \sigma (s, X_s)dB_s, \label{sde}
\end{align}
for any $(t_0, x_0)\in [0, T]\times \mathbb{R}^n$, 
$t\in [t_0, T]$.  
Denote 
  a solution to the BSDE 
\begin{align}
Y_t = g(X_T^{t_0, x_0})
+\int_t^T F(s, X_s^{t_0, x_0}, Y_s, Z_s)ds -\int_t^T Z_sdW_s, \ t\in[t_0, T],\label{fbsde}
\end{align}
by  $(Y^{t_0, x_0}, Z^{t_0, x_0})$ or $(Y, Z)$ when there is no ambiguity. 
The probablistic representation for nonlinear Feymann-Kac formula consists of proving that,  in Markovian setting, $u(t, x):=Y_t^{t, x}$
is a solution at least in the viscosity sense to  (\ref{pde}) when the source of nonlinearity $F$ is quadratic in $\nabla_x u(t, x)$ and  $g$ is an unbounded function.  To put it more precisely, let us  introduce the FBSDEs.

{\bf The Forward Markovian SDEs.}
Let $b: [0, T]\times\mathbb{R}^n \rightarrow \mathbb{R}^n$, $\sigma: [0, T]\times\mathbb{R}^d\rightarrow \mathbb{R}^{n\times d}$ be continuous functions and assume there exists $\beta \geq 0$ such that $\mathbb{P}$-a.s. for any $t\in[0, T]$, $|b(t, 0)| +|\sigma(t, 0)| \leq \beta$ and $b(t, x), \sigma(t, x)$ are  Lipschitz-continuous in $x$, i.e., $\mathbb{P}$-a.s.  for any $t\in[0, T]$, $x, x^\prime \in\mathbb{R}^n$, 
\[
|b(t, x)- b(t, x^\prime)|
+|\sigma(t, x)- \sigma(t, x^\prime)|\leq \beta |x-x^\prime|.
\]
Then for any $(t_0, x_0)\in [0, T]\times\mathbb{R}^n$,  (\ref{sde}) has a unique solution $X^{t_0, x_0}$ in  $\mathcal{S}^p$ for any $p \geq 1$.  

{\bf The Markovian BSDE.} We continue with the setting of the forward equations above. Set $q\geq 1$.
Let $F_1, F_2: [0, T]\times\mathbb{R}^n\times\mathbb{R}\times\mathbb{R}^d \rightarrow \mathbb{R}$, 
$g: \mathbb{R}^n\rightarrow \mathbb{R}$ be continuous functions, $\varphi: \mathbb{R}^+ \rightarrow \mathbb{R}^+$  a continuous nondecreasing function with $\varphi(0)=0$ and
$f(|\cdot|)\in \mathcal{I}$, and assume moreover  $F = F_1 + F_2$ such  that
\begin{enumerate}
\item[(i)] 
 $F_1(t, x, y, z)$ is  monotonic in $y$ and  Lipschitz-continuous in $z$, and $F_2(t, x, y, z)$ is monotonic at $y=0$ and of linear-quadratic growth in $z$, 
 i.e., for any $(t, x)\in[0, T]\times\mathbb{R}^n$, $y, y^\prime \in \mathbb{R}, z, z^\prime \in \mathbb{R}^d$,
\begin{align*}
\sgn (y -y ^\prime)\big(F_1(t, x, y, z)- F_1(t,  x, y^\prime, z)\big)&\leq \beta |y- y^\prime|, \\
 \big|F_1(t, x, y, z)- F_1(t, x, y, z^\prime)\big| &\leq  \beta |z-z^\prime|, \\
 \sgn(y)F_2(t, x, y, z) &\leq \beta|y| + \beta |z| + f(|y|)|z|^2;
\end{align*}
\item[(ii)] $(y, z) \longmapsto F_2(t, x, y, z)$ is  convex ;
\item[(iii)] for any $(t, x, y, z)\in [0, T]\times\mathbb{R}^n\times\mathbb{R}\times\mathbb{R}^d$, 
\begin{align*}
| F (t, x, y, z) | &\leq  \beta \big(1+ |x|^q + 2|z|\big) + \varphi(|y|)
  + f(|y|)|z|^2, \\
 |g(x)| &\leq \beta  \big(1 + |x|^q\big). 
  \end{align*}
\end{enumerate} 

Since $X^{t_0, x_0}\in \mathcal{S}^p$ for any $p \geq 1$, the  above structure conditions on  $F$ and $g$ allow one to use Corollary \ref{lpunique}
to construct a unique solution $(Y^{t_0, x_0}, Z^{t_0, x_0})$ in $\mathcal{S}^p\times\mathcal{M}^p$ to 
 (\ref{fbsde})
 for any $p > 1$.
Moreover, by standard arguments,
$Y_{t_0}^{t_0, x_0}$
 is deterministic for any $(t_0, x_0) \in [0, T]\times\mathbb{R}^n$. Hence $u(t, x)$ defined as $Y_t^{t, x}$ is a deterministic function. With this fact 
 we now turn to the main result of this section: $u$ is a viscosity solution to  (\ref{pde}). Before our proof let us recall the definition of a viscosity solution. 
 
{\bf Viscosity Solution.} A continuous function $u: [0, T]\times\mathbb{R}^n \rightarrow \mathbb{R}$ is called a viscosity subsolution (respectively supersolution) to (\ref{pde}) if   $u(T, x) \leq g(x)$ (respectively $u(T, x) \geq g(x)$)
and 
for any smooth function  $\phi$ such that  $u-\phi$ reaches the local maximum (respectively local minimum) at $(t_0, x_0)$, we have 
\[
\partial_t \phi(t_0, x_0)
+\mathcal{L}\phi(t_0, x_0)
+ F(t_0, x_0, u(t_0, x_0), \sigma^\top \nabla_x \phi(t_0, x_0))\geq 0\ (\text{respectively} \leq 0).
\]
A function $u$ is called a viscosity solution to (\ref{pde}) if it is both a viscosity subsolution and supersolution. 
\begin{prop}
Given the above assumptions, $u(t, x)$  is continuous with 
\[
|u(t, x)| \leq c\big( 1+ |x|^q\big),
\] where $c$ is a constant. Moreover, $u$ is a viscosity solution to {\rm(\ref{pde})}.
\end{prop}
\begin{Proof}
 Due to the Lipschitz-continuity of $b$ and $\sigma$, $X^{t, x}$ is continuous in $(t, x)$, e.g.,  in mean square sense. The continuity of $u$ is then an immediate consequence of  Theorem \ref{lpstability} (stability). The proof relies on standard arguments and  hence is omitted.
 By  Lemma \ref{lpestimate2} (a priori estimate (ii)), we prove that $u$ satisfies the above polynomial growth. It thus remains to prove that $u$ is a viscosity solution to (\ref{pde}).

Let $\phi$ be a smooth function such that $u-\phi$ reaches local maximum at $(t_0, x_0)$. Without loss of generality we  assume that the local maximum is global and 
$u(t_0, x_0) = \phi(t_0, x_0).$ We aim at proving 
\[
\partial_t \phi(t_0, x_0)
+\mathcal{L}\phi(t_0, x_0)
+ F(t_0, x_0, u(t_0, x_0), \sigma^\top \nabla_x \phi(t_0, x_0))\geq 0.
\]
From (\ref{fbsde}) we obtain
\[
Y_t = Y_{t_0}
-\int_{t_0}^t F(s, X_s^{t_0, x_0}, Y_s, Z_s)ds +\int_{t_0}^t Z_sdW_s. \]
By It\^{o}'s formula, 
\[
\phi(t, X_t^{t_0, x_0}) = 
\phi(t_0, x_0)
+
\int_{t_0}^{t}
\big\{\partial_s\phi+\mathcal{L}\phi \big\} (s, X_s^{t_0, x_0})ds +
\int_{t_0}^{t}
\sigma^{\top} \nabla_x \phi(s, X_s^{t_0, x_0}) dW_s.  
\]
Now we take any $t\in [t_0, T]$.
Note that the existence of a unique solution to (\ref{sde}) and (\ref{fbsde})  implies by Markov property 
that $Y_t =u(t, X_t^{t_0, x_0})$.  Hence,
 $\phi(t, X_t^{t_0, x_0})\geq u(t, X_t^{t_0, x_0}) =    Y_t$.  By touching property, on the set
$
\big\{ \phi(t, X_t^{t_0, x_0}) = Y_t\big\}
$ we have
\begin{align*}
&\partial_t \phi(t, X_t^{t_0, x_0}) +\mathcal{L}\phi(t, X_t^{t_0, x_0}) + F(t, X_t^{t_0, x_0}, Y_t, Z_t) \geq 0\ \ \mathbb{P}\text{-a.s.}, \\
&\sigma^\top \nabla_x\phi(t, X_t^{t_0, x_0}) -Z_t =0\ \  \mathbb{P}\text{-a.s.}
\end{align*}
Now we set $t=t_0$. We have
$\phi(t_0, X_{t_0}^{t_0, x_0}) = \phi(t_0, x_0) = u(t_0, x_0) = Y_{t_0}$. Moreover, 
the above equality implies 
$Z_{t_0} = \sigma^\top \nabla_x \phi(t_0, x_0)$.  Plugging the two equalities into the above inequality gives
\[
\partial_t \phi(t_0, x_0)
+\mathcal{L}\phi(t_0, x_0)
+ F(t_0, x_0, u(t_0, x_0), \sigma^\top \nabla_x \phi(t_0, x_0))\geq 0.
\]
 Hence   $u$ is a viscosity subsolution to (\ref{pde}).  $u$ being a viscosity supersolution and thus a viscosity solution can be proved analogously.
\qed
\end{Proof}

\chapter{Quadratic Semimartingale BSDEs}
\label{qsbsde}
\section{Preliminaries}
\label{section20}
The objectives of our study in this chapter are quadratic BSDEs driven by continuous local martingales.   We fix the time horizon $T>0$, and work on a filtered probability space $(\Omega, \mathcal{F}, (\mathcal{F}_t)_{t\in [0, T]}, \mathbb{P})$ satisfying the usual conditions of right-continuity and $\mathbb{P}$-completeness. 
 $\mathcal{F}_0$ is the $\mathbb{P}$-completion of the trivial $\sigma$-algebra.
Any measurability  will refer to the filtration $(\mathcal{F}_t)_{t\in [0, T]}$.  In particular,  $\Prog$ denotes the progressive $\sigma$-algebra on $\Omega \times [0, T]$.
We assume the filtration  is \emph{continuous}, in the sense that all  local martingales have $\mathbb{P}$-a.s. continuous sample paths.  
 $M = (M^1,..., M^d)^\top$  stands for a  fixed $d$-dimensional   continuous local martingale.  By \emph{continuous semimartingale setting} we mean:  $M$ doesn't have to be a Brownian motion;  the filtration is not necessarily generated by $M$ which is usually seen as the main source of randomness. Hence in various concrete situations there may be a continuous  local martingale $N$ strongly
orthogonal to $M$. As mentioned in the introduction, 
we exclusively study $\mathbb{R}$-valued BSDEs. They can be written as
\[
Y_t = \xi + \int_t^T\big( \mathbf{1}^\top
d\langle M\rangle_s F(s, Y_s, Z_s)
+ g_s d\langle N\rangle_s\big) -\int_t^T \big(Z_s dM_s + dN_s\big),
\]
where $\mathbf{1}:=(1, ..., 1)^\top$, $\xi$ is an $\mathbb{R}$-valued $\mathcal{F}_T$-measurable random variable, $F: \Omega \times [0, T]
\times \mathbb{R} \times \mathbb{R}^d \rightarrow \mathbb{R}^d$ is a $\Prog\otimes \mathcal{B}(\mathbb{R})\otimes
\mathcal{B}(\mathbb{R}^d)$-measurable random function and $g$ is an $\mathbb{R}$-valued $\Prog$-measurable bounded process.  $\int_0^\cdot (Z_sdM_s + dN_s)$, sometimes denoted by $Z\cdot M + N$,
refers to the vector stochastic integral; see Shiryaev and Cherny \cite{SC2002}.
The equations defined in this way encode the matrix-valued process $\langle M\rangle$  which is not amenable to analysis. Therefore we rewrite the BSDEs by factorizing
$\langle M \rangle$. This procedure separates the matrix property from its nature as a measure. It can also be regarded as a reduction of dimensionality.

There are many ways to factorize $\langle M \rangle$; see, e.g., Section III. 4a, Jacod and  Shiryaev \cite{JS1987}. We can and  choose
$A: = \arctan\big( \sum_{i=1}^d {\langle M^i\rangle}\big)$.   By Kunita-Watanabe inequality, we deduce the absolute continuity of  $\langle M^i, M^j \rangle$  
with respect to $A$. Note that such choice makes $A$ continuous, increasing and bounded.
Moreover, by  Radon-Nikod\'{y}m theorem and Cholesky decomposition, there exists a matrix-valued $\Prog$-measurable process $\lambda$  
such that $\langle M\rangle  = (\lambda^\top \lambda) \cdot A.$ 
As will be seen later, our results
don't rely on the specific choice of $A$ but only on its boundedness. In particular, if $M$ is a $d$-dimensional Brownian motion, we may choose $A_t =t$ and $\lambda$ to be the  identity matrix.

 The second advantage of  factorizing $\langle M\rangle$  is that
\[
\mathbf{1}^\top d\langle M\rangle_s F(s, Y_s, Z_s) =
\mathbf{1}^\top \lambda_s^\top \lambda_s F(s, Y_s, Z_s)dA_s,
\]
where $f(t, y, z) := \mathbf{1}^\top \lambda_s^\top \lambda_s F(s, y, z)$  is $\mathbb{R}$-valued. Such reduction of dimensionality makes it easier to formulate the difference of two equations as frequently appears in comparison theorem and uniqueness.
Hence, we may reformulate the BSDEs as follows.

{\bf BSDEs: Definition and Solutions.}
Let $A$ be an $\mathbb{R}$-valued continuous  nondecreasing {bounded} adapted process such that 
$\langle M \rangle = (\lambda^\top \lambda) \cdot A$ for some matrix-valued $\Prog$-measurable process $\lambda$, $f: \Omega \times [0, T]
\times \mathbb{R} \times \mathbb{R}^d \rightarrow \mathbb{R}$
  a $\Prog\otimes \mathcal{B}(\mathbb{R})\otimes
\mathcal{B}(\mathbb{R}^d)$-measurable random function,
$g$  an $\mathbb{R}$-valued $\Prog$-measurable bounded process and $\xi$  an $\mathbb{R}$-valued $\mathcal{F}_T$-measurable random variable.  The semimartingale BSDEs are written as
\begin{align}
Y_t= \xi + &\int_t^T \big(f(s, Y_s, Z_s)dA_s + g_s d\langle N\rangle_s\big) -
\int_t^T \big(Z_s dM_s + dN_s\big).  \label{sbsde}
\end{align}
 We call a process $(Y, Z, N)$   or $(Y, Z\cdot M + N)$ a \emph{solution} to (\ref{sbsde}), if $Y$ is an $\mathbb{R}$-valued continuous adapted process, $Z$ is an $\mathbb{R}^d$-valued $\Prog$-measurable process and $N$  is 
an $\mathbb{R}$-valued 
  continuous local martingale strongly orthogonal to $M$, such that $\mathbb{P}$-a.s. $\int_0^T Z_s^\top d\langle M\rangle_s Z_s <+\infty$ and 
$\int_0^T |f(s, Y_s, Z_s)|dA_s <+\infty$, and  (\ref{sbsde}) holds   
 $\mathbb{P}$-a.s. for all $t\in [0, T]$,

Note that the factorization of $\langle M\rangle$ gives $
 \int_0^\cdot Z_s^\top d\langle M\rangle_s Z_s  = \int_0^\cdot |\lambda_s Z_s|^2 dA_s.
$
Hence we don't  distinguish these two integrals  in all situations. 
$\int_0^T Z_s^\top d\langle M\rangle_s Z_s <+\infty\ \mathbb{P}$-a.s. ensures  that $Z$ is integrable with respect to $M$ in the sense of vector stochastic integration.
As a result, $Z\cdot M$ is a continuous local martingale. 
$M$ and $N$ being continuous and strongly orthogonal implies that
$\langle M^i, N\rangle = 0$ for  $i = 1, ..., d$.
We call $f$ the \emph{generator}, $\xi$ the \emph{terminal value} and   $(\xi, \int_0^T|f(s, 0, 0)|dA_s)$ the \emph{data}. 
In our study, the integrability property of the data determines the estimates for a solution. The conditions imposed on the generator are called the \emph{structure conditions}. For notational convenience, we sometimes write $(f, g, \xi)$ instead of (\ref{sbsde}) to denote the above BSDE.  Finally, 
(\ref{sbsde})
 is called \emph{quadratic} if $f$ has at most quadratic growth in $z$ or $g$ is not indistinguishable from $0$.

Regarding the existence results,  most literature requires $g$ to be a constant;  see, e.g.,   \cite{EK1997}, \cite{M2009}, \cite{MW2012}. The reason is that 
 $g\cdot \langle N\rangle$ can be eliminated via exponential transform only if $g$ is a constant. 
Tevzadze \cite{T2008} allows $g$ to be any bounded process but their results are less general in several aspects.  We also point out that
in mathematical finance, $g$ usually appears as a constant; see, e.g., \cite{MS2005}, \cite{BMS2007}, \cite{HS2011}, \cite{FMW2012}.

We take a further step by studying bounded and unbounded solutions 
to BSDEs associated with 
 any bounded process $g$, and with  monotonicity at $y =0$ and at most quadratic growth in $z$.  The conditions  to our knowledge are the most general compared to existing literature.
We start from bounded solutions to Lipschitz-quadratic BSDEs (see Section \ref{section21}) and then extend the results to general quadratic BSDEs (see Section \ref{section22}, \ref{section23}). 

Let us close this section by introducing all required notations  for this chapter. $\ll$ stands for the strong order of nondecreasing processes, stating that the difference is nondecreasing.
For any random variable or process $Y$, we say $Y$ has some property if this is true except on a $\mathbb{P}$-null subset of $\Omega$. Hence we omit ``$\mathbb{P}$-a.s'' in situations without ambiguity. Define $\sgn(x)= \mathbb{I}_{\{x\neq 0 \}}\frac{x}{|x|}$. For any random variable $X$, define $\norm{X}_\infty$ to be its essential supremum.
For any c\`adl\`ag adapted process $Y$, set $Y_{s, t}: = Y_t -Y_s$ and $Y^* : = \sup_{t\in [0, T]} |Y_t|$.
For any $\Prog$-measurable process $H$, set $|H|_{s,t}:= \int_s^t H_u dA_u$ and $|H|_t : = |H|_{0, t}$. 
 $\mathcal{T}$ stands for the set of all stopping times valued in $[0, T]$ and $\mathcal{S}$ denotes the space of continuous adapted processes. 
For later use we specify the following spaces under $\mathbb{P}$.
\begin{itemize}
\item $\mathcal{S}^\infty$: the space of  bounded processes  $Y\in\mathcal{S}$ with $\norm{Y}:=\norm{Y^*}_\infty$; $\mathcal{S}^\infty$ is a  Banach space;

	\item $\mathcal{M}$: the set of  continuous local martingales starting from $0$; for any $\mathbb{R}^d$-valued $\Prog$-measurable process $Z$ with $\int_0^T Z_s^\top d\langle M\rangle_s Z_s < +\infty$,  $Z\cdot M \in \mathcal{M}$; 
\item $\mathcal{M}^p(p\geq 1)$: the set 
of $\widetilde{M} \in \mathcal{M}$ with 
\[
\norm{\widetilde{M}}_{\mathcal{M}^p}:=\big(\mathbb{E}\big[ \langle \widetilde{M}\rangle_T^{\frac{p}{2}}\big] \big)^{\frac{1}{p}} < +\infty;
\]
 in particular, $\mathcal{M}^2$ is a Hilbert space;
\item $\mathcal{M}^{BMO}$: the 
 set  of BMO martingales $\widetilde{M}\in \mathcal{M}$ \ with 
\[ \norm{\widetilde{M}}_{BMO} :=\sup_{\tau \in \mathcal{T}}
\big\|\mathbb{E}\big[
\langle \widetilde{M}\rangle_{\tau, T}
\big|\mathcal{F}_\tau    \big]^{\frac{1}{2}}\big\|_\infty;
\]
$\mathcal{M}^{BMO}$ is a Banach space.
\end{itemize}

 $\mathcal{M}^2$ being a Hilbert space is crucial to proving convergence of the  martingale parts in the monotone stability result of quadratic BSDEs (see, e.g.,  Kobylanski \cite{K2000},  Briand and  Hu \cite{BH2008},  Morlais \cite{M2009} or Section \ref{section22}). 
 Other spaces  are also Banach  under suitable norms; we will not present these facts in more detail since they are not involved in our study.
  
Finally, for any local martingale $\widetilde{M}$, we call 
$\{\sigma_n \}_{n\in\mathbb{N}^+}\subset \mathcal{T}$ a \emph{localizing sequence} if $\sigma_n$ increases stationarily to $T$ as $n$  goes to $+\infty$ and $\widetilde{M}_{\cdot\wedge \sigma_n}$ is a martingale for any $n\in\mathbb{N}^+$.

\section{Bounded Solutions to Lipschitz-quadratic BSDEs}
\label{section21}
This section takes one step in   solving   quadratic BSDEs and consists  in the study of equations with Lipschitz-continuous generators. In contrast to  El Karoui and  Huang \cite{EK1997}, we allow the presence of  $g\cdot \langle N\rangle$.
We point out that similar results
for linear-quadratic generators have been studied by
Tevzadze \cite{T2008}, but the case of Lipschitz-continuity is not available in that work.   Due to its importance for regularizations of quadratic BSDEs, we study existence and uniqueness   results for equations of this particular type  in the first step.  To this end, we assume
\begin{as}\label{as1} 
There exist $\beta, \gamma \geq 0$
such that 
${\norm{\xi}}_\infty + \big\|{\big| |f(\cdot, 0, 0) | \big|_T}\big\|_\infty <+\infty$ and 
$f$ is Lipschitz-continuous in $(y, z)$, i.e., $\mathbb{P}$-a.s.  for any $t\in[0, T]$, $y, y^\prime \in\mathbb{R}$, $z, z^\prime \in \mathbb{R}^d$,
\[
|f(t, y, z ) - f(t, y^\prime, z^\prime)| \leq \beta |y-y^\prime| + \gamma |\lambda_t (z-z^\prime)|.
\]
\end{as}

Due to the presence of $g\cdot\langle N\rangle$,  we call the BSDE $(f, g, \xi)$ satisfying \ref{as1}   \emph{Lipschitz-quadratic}. Given \ref{as1}, 
we are about to construct a solution in the  space $\mathscr{B}: = \mathcal{S}^\infty\times \mathcal{M}^{BMO}$ equipped with  the norm 
\[
\norm {(Y, Z\cdot M + N)} : = \big(\norm{Y}^2 + \norm{Z\cdot M + N }^2_{BMO}\big)^\frac{1}{2},
\] 
for $(Y, Z\cdot M  + N)\in \mathcal{S}^\infty\times\mathcal{M}^{BMO}$. Clearly $(\mathscr{B}, \norm{\cdot})$ is Banach. 
As a preliminary result,  we claim that   the existence result holds given sufficiently small data.
\begin{thm}[Existence (i)]\label{existence1} If $(f, g, \xi)$ satisfies  {\rm\ref{as1}} with 
\begin{align}
\norm{\xi}^2_{\infty}+8 \big\|{\big||f(\cdot, 0, 0)|\big|_T}\big\|^2_{\infty} \leq \frac{1}{64}\exp \Big(       -\norm{A} 
\big(8\beta^2 \norm{{A}} + 8\gamma^2\big)
    \Big) \label{p1}
\end{align}
 and $\mathbb{P}$-a.s.  $|g_\cdot|\leq \tilde{g}:= \frac{1}{8}$, then
 there exists a solution  in $(\mathscr{B}, \norm{\cdot})$.
\end{thm}	
\begin{Proof}  To overcome the difficulty arising from  the Lipschitz-continuity, we use Banach fixed point theorem under an equivalent norm.
Set $\rho \geq 0$ to be determined later. 
For any  $X \in\mathbb{L}^\infty, Y\in\mathcal{S}^\infty$ and $\widetilde{M} \in\mathcal{M}^{BMO}$, set
$\norm{X}_{ \infty, \rho} : = 
\norm{e^{\frac{\rho}{2} A_T}X }_\infty$,
 $\norm{Y}_\rho : = \norm{e^{\frac{\rho}{2}A} Y } $ and $\norm{\widetilde{M}}_{BMO, \rho} : = \norm{e^{\frac{\rho}{2} A}\cdot \widetilde{M}}_{BMO}$;   for $(Y, Z\cdot M +N)\in \mathscr{B}$, set
\[
\norm{(Y, Z\cdot M +N)}_\rho :  =  \big(\norm{Y}_\rho^2 + \norm{Z\cdot M + N}_{BMO, \rho}^2\big)^{\frac{1}{2}}.
\]
Since $A$ is bounded, $\norm{\cdot}_\rho$  is equivalent to the original norm for each space. Hence
$(\mathscr{B}, \norm{\cdot}_\rho)$ is also a Banach space.
For any $R\geq 0$,  define
\[
\mathbf{B}_R : = \big\{ (Y, Z\cdot M + N) \in \mathscr{B}: \norm{(Y, Z\cdot M +N)}_\rho \leq R   \big\}.
\]
We show by  Banach fixed point theorem that there exists a unique solution in $\mathbf{B}_{R}$ with $R=\frac{1}{2}.$
To this end, we define $\mathbf{F}: (\mathbf{B}_R, \norm{\cdot}_\rho)\rightarrow (\mathscr{B},  \norm{\cdot}_\rho)
$ such that for any $(y, z\cdot M + n) \in\mathbf{B}_R$,  $(Y, Z\cdot M + N) : =\mathbf{F}((y, z\cdot M + n))$ solves
\[
Y_t = \xi + \int_t^T \big( f(s, y_s, z_s) dA_s + g_s d\langle n \rangle_s \big) -\int_t^T \big(Z_s dM_s + dN_s \big).
\]
Indeed, such $(Y, Z, N)$  uniquely exists due to martingale representation theorem. Moreover, by standard estimates, $(Y, Z\cdot M + N) \in (\mathscr{B}, \norm{\cdot}_\rho)$.

(i). We show
$
\mathbf{F}(\mathbf{B}_R) \subset \mathbf{B}_R.
$
For any $\tau \in \mathcal{T}$,  It\^o's formula applied to $e^{\rho A_\cdot}Y^2_\cdot$ yields
\begin{align}
&e^{\rho A_\tau}|Y_\tau|^2 + \rho \mathbb{E}\Big[ \int_\tau^T e^{\rho A_s}Y_s^2 dA_s \Big|\mathcal{F}_\tau \Big]
+ \mathbb{E}\Big[\int_\tau^T e^{\rho A_s} \Big(Z_s^\top d\langle M\rangle_s  Z_s + d\langle N\rangle_s\Big)\Big|\mathcal{F}_\tau \Big] 
\nonumber\\ 
&\leq
\norm{\xi}^2_{\infty, \rho}  + 2\mathbb{E}\Big[\int_\tau^T e^{\rho A_s}|Y_s ||f(s, y_s, z_s)|dA_s \Big| \mathcal{F}_\tau\Big] + 2\mathbb{E}\Big[ \int_\tau^T e^{\rho A_s}|Y_s||g_s| d\langle n\rangle_s\Big| \mathcal{F}_\tau\Big].
\label{lq1}
\end{align}
By  \ref{as1},
\begin{align*}
|Y_s||f(s, y_s, z_s)| \leq |Y_s||f(s, 0, 0)| + \beta |Y_s| |y_s| + \gamma |Y_s| |\lambda_s z_s|.
\end{align*}
We plug this inequality into (\ref{lq1})
and estimate each term on the right-hand side.  Using $2ab \leq \frac{1}{8}a^2 + 8b^2$  gives  
\begin{align*}
2\mathbb{E}\Big[  \int_\tau^T  e^{\rho A_s} |Y_s||f(s, 0, 0)| dA_s     \Big| \mathcal{F}_\tau\Big]
&\leq 
\frac{1}{8}\norm{Y}^2_\rho +  8 \mathbb{E}\Big[ \int_\tau^T e^{\frac{\rho}{2}A_s} |f(s, 0, 0)|dA_s  \Big|\mathcal{F}_\tau\Big]^{2} \\
&\leq \frac{1}{8}\norm{Y}^2_\rho + 8 \big\|{\big||f(\cdot, 0, 0)|\big|_T}\big\|^2_{\infty, \rho},
\end{align*}
\begin{align*}2
\beta \mathbb{E}\Big[ \int_\tau^T e^{\rho A_s}|Y_s||y_s| dA_s \Big|\mathcal{F}_\tau\Big]
&\leq
\frac{1}{8}\norm{y}^2_\rho + 8\beta^2 \mathbb{E}\Big[\int_\tau^T e^{\frac{\rho}{2}A_s} |Y_s| dA_s    \Big| \mathcal{F}_\tau \Big]^{2}\\
& \leq \frac{1}{8}\norm{y}^2_\rho + 8\beta^2\norm{A} \mathbb{E}\Big[\int_\tau^T e^{
\rho A_s} |Y_s|^2 dA_s    \Big| \mathcal{F}_\tau \Big], 
\end{align*}
\begin{align*}
2\gamma \mathbb{E}\Big[ \int_\tau^T e^{\rho A_s}|Y_s||\lambda_s z_s| dA_s \Big|\mathcal{F}_\tau\Big]
& \leq \frac{1}{8}\norm{z\cdot M}^2_{BMO, \rho} + 8\gamma^2 \mathbb{E}\Big[\int_\tau^T e^{\rho A_s} |Y_s|^2 dA_s    \Big| \mathcal{F}_\tau \Big],
\end{align*}
\begin{align*}
2 \mathbb{E}\Big[\int_\tau^T e^{\rho A_s }|Y_s| |g_s|\langle N \rangle _s \Big| \mathcal{F}_\tau \Big]
&\leq \frac{1}{8}\norm{Y}^2_\rho + 8 \tilde{g}^2\mathbb{E}\Big[ \int_\tau^T e^{\frac{\rho}{2}A_s} d\langle N\rangle_s \Big| \mathcal{F}_\tau        \Big]^2 \\
&\leq \frac{1}{8}\norm{Y}^2_\rho + 8 \tilde{g}^2 \norm{n}^4_{BMO, \rho}.
\end{align*}
Set $\rho :=  8\beta^2 \norm{A} + 8\gamma^2$  so as to eliminate 
$\mathbb{E}\big[ \int_\tau^T e^{\rho A_s}Y_s^2 dA_s \big|\mathcal{F}_\tau \big]$ on both sides. Hence  (\ref{lq1}) gives
\begin{align}
e^{\rho A_\tau}|Y_\tau|^2 
&+ \mathbb{E}\Big[\int_\tau^T e^{\rho A_s} \Big(Z_s^\top d\langle M\rangle_s  Z_s + d\langle N\rangle_s\Big)\Big|\mathcal{F}_\tau \Big]  \nonumber \\
& \leq 
\norm{\xi}^2_{\infty, \rho} +8  \big\| \big| |f(\cdot, 0, 0)|\big|_T\big\|^2_{\infty, \rho} + \frac{1}{4} \norm{Y}^2_\rho  \nonumber\\
&+ \frac{1}{8} \big(\norm{y}^2_\rho + \norm{z\cdot M}_{BMO, \rho}^2 \big)
+ 8 \tilde{g}^2 \norm{n}^4_{BMO, \rho}. \label{lq11}
\end{align}
Taking essential supremum and  supremum over all $\tau \in \mathcal{T}$, and using the inequality
\begin{align*}
\frac{1}{2}\norm {(Y, Z\cdot M +N)}_\rho^2 &\leq   \norm{Y}^2_\rho \vee \norm{Z\cdot M + N}^2_{BMO, \rho}\\
 &\leq  \sup_{\tau \in \mathcal{T}}\Big\| e^{\rho A_\tau}|Y_\tau|^2 +  \mathbb{E}\Big[ \int_\tau^T e^{\rho A_s} \Big(Z_s^\top d\langle M\rangle_s Z_s + d\langle N \rangle_s \Big)\Big| \mathcal{F}_\tau\Big]\Big\|_\infty, 
\end{align*}
we deduce by transferring $\frac{1}{4} \norm{Y}^2_\rho$ to the left-hand side of (\ref{lq11}) that 
\begin{align*}
\norm{(Y, Z\cdot M + N)}^2_\rho &\leq  4\norm{\xi}_{\infty, \rho}^2 + 32  \big\|{\big||f(\cdot, 0, 0)|\big|_T}\big\|^2_{\infty, \rho}+ \frac{1}{2} \big(\norm{y}_\rho^2 + \norm{z\cdot M}^2_{BMO, \rho} \big) +
32  \tilde{g}^2 \norm{n}^4_{BMO, \rho} \\
& \leq 4\norm{\xi}^2_{\infty, \rho} + 32  \big\|{\big||f(\cdot, 0, 0)|\big|_T}\big\|^2_{\infty, \rho}  + \frac{1}{2}R^2  + 32 \tilde{g}^2 R^4. 
\end{align*}
Thanks to (\ref{p1}), $\tilde{g}=\frac{1}{8}$ and $R=\frac{1}{2}$,  we verify from the above estimate that 
\[
\norm{(Y, Z, N)}_\rho \leq R.
\]

(ii).
We  prove $\mathbf{F}: (\mathbf{B}_R, \norm{\cdot}_\rho) \rightarrow (\mathbf{B}_R, \norm{\cdot}_\rho)$ is a contraction mapping. 
By (i),  for $i=1, 2$ and  any $(y^i, z^i\cdot M + n^i)\in \mathbf{B}_R$,   we have  $(Y^i, Z^i\cdot M + N^i) : = \mathbf{F}((y^i, z^i\cdot M + n^i)) \in \mathbf{B}_R$. 
For notational convenience we set 
 $\delta y: = y^1 - y^2$ and  $\delta z, \delta n, \delta \langle n\rangle, \delta Y, \delta Z, \delta N, \delta \langle N\rangle$, etc.  analogously.
By  the deductions  in (i) with  minor modifications, we obtain  
\begin{align}
\frac{1}{2}\norm{(\delta Y, \delta Z\cdot M + \delta N)}^2_\rho &\leq   \frac{1}{8} \big(\norm{\delta y}^2_\rho + \norm{\delta z\cdot M}^2_{BMO, \rho} \big)   + \frac{1}{4}\norm{\delta Y}^2_\rho \nonumber\\
&+ 
4\tilde{g}^2  \sup_{\tau \in \mathcal{T}} \Big \| \mathbb{E} \Big[\int_\tau^T e^{\frac{\rho}{2}A_s} d |\delta \langle n \rangle _s|              \Big| \mathcal{F}_\tau \Big]^2\Big\|_\infty.\label{lq111}
\end{align}
 Kunita-Watanabe inequality and Cauchy-Schwartz inequality  used to the last term gives
\begin{align*}
\mathbb{E} \Big[\int_\tau^T e^{\frac{\rho}{2}A_s} d |\delta \langle n \rangle _s|              \Big| \mathcal{F}_\tau \Big]^2 &\leq
\mathbb{E}\Big[  \int_\tau^T  e^{\frac{\rho}{2}A_s}d\langle \delta n \rangle_s  \Big| \mathcal{F}_\tau    \Big]
 \mathbb{E}\Big[  \int_\tau^T  e^{\frac{\rho}{2}A_s}d\langle n^1 + n^2 \rangle_s  \Big| \mathcal{F}_\tau    \Big]
\\
& \leq  \norm{\delta n}^2_{BMO, \rho}  \cdot  2 \big(\norm{n^1}_{BMO, \rho}^2 + \norm{n^2}_{BMO, \rho}^2\big) \\
&\leq \norm{\delta n}^2_{BMO, \rho} \cdot 4R^2,
\end{align*}
where the last inequality is due to $\norm{(y^i, z^i\cdot M + n^i)}_\rho \leq R$, $i = 1, 2.$
Hence (\ref{lq111}) gives
\begin{align*}
\norm{(\delta Y, \delta Z\cdot M + \delta N)}^2_\rho &\leq \frac{1}{2} \big(\norm{\delta y}^2_\rho + \norm{\delta z\cdot M}^2_{BMO, \rho} \big) +  64 \tilde{g}^2 R^2 \norm{\delta n}^2_{BMO, \rho} \\
&\leq \Big(\frac{1}{2} + 64 \tilde{g}^2 R^2\Big) \norm{(\delta y, \delta z\cdot M + \delta n)}^2_\rho \\
&\leq \frac{3}{4}\norm{(\delta y, \delta z \cdot M + \delta n)}^2_\rho,
\end{align*}
i.e.,  $\mathbf{F}: (\mathbf{B}_R, \norm{\cdot}_\rho) \rightarrow (\mathbf{B}_R, \norm{\cdot}_\rho)$ is a contraction mapping. The existence of a solution in $\mathbf{B}_R$ thus follows immediately from Banach fixed point theorem. Finally, since $\norm{\cdot}$ is equivalent to $\norm{\cdot}_\rho$ for $\mathscr{B}$,  the solution also belongs to $(\mathscr{B}, \norm{\cdot})$.
\qed
\end{Proof}

From now on we denote $(\mathscr{B},\norm{\cdot})$ by $\mathscr{B}$ when there is no ambiguity.
 In the spirit of  Tevzadze \cite{T2008},  we  extend this existence result so as to allow any bounded data. 
 To this end, for any $\mathbb{Q}$ equivalent to $\mathbb{P}$ we define 
$\mathcal{S}^\infty (\mathbb{Q})$ analogously to $\mathcal{S}^\infty$ but under $\mathbb{Q}$. This notation also applies to other spaces.
\begin{thm}[Existence (ii)]\label{existence2} If $(f, g, \xi)$ satisfies  {\rm\ref{as1}}, then there exists a solution to $(f, g, \xi)$ in $\mathscr{B}$.
\end{thm}
\begin{Proof}
\par
\noindent
(i).  We first show that it is equivalent to prove the existence result  given $|g_\cdot|\leq \frac{1}{8}$ $\mathbb{P}$-a.s. Suppose that $g$ is bounded by a positive constant $\tilde{g}$, that is, $|g_\cdot|\leq \tilde{g}$ $\mathbb{P}$-a.s.
Observe that,  for any $\theta >0$, $(Y, Z, N)$ is a solution  to $(f, g, \xi)$ if and only if $(\theta Y, \theta Z, \theta N)$  is a solution to $(f^\theta, g/\theta, \theta \xi)$, where
$f^\theta(t, y, z): = {\theta}f(t, \frac{y}{\theta},\frac{z}{\theta}).$
Obviously $f^\theta$ verifies \ref{as1} with the same Lipschitz coefficients as $f$. 
If we set $\theta  := 8  \tilde{g}$, then ${|g_\cdot/\theta|} \leq \frac{1}{8}$ $\mathbb{P}$-a.s. and hence satisfies the parametrization in Theorem \ref{existence1} (existence (i)). Therefore, we can and do assume  $|g_\cdot|\leq \frac{1}{8}$ $\mathbb{P}$-a.s. without loss of generality. 

(ii). Since $\norm{\xi}_\infty + \big\|{\big||f (\cdot, 0, 0)|\big|_T}\big\|_\infty < + \infty$, we can find  $n\in\mathbb{N}^+$ such that
\[
\xi=\sum_{i=1}^n \xi^i, \ 
f(t, 0, 0) = \sum_{i=1}^n f^i (t, 0, 0),
\]
where, for each $i \leq n$, $\xi^i$ is a  $\mathcal{F}_T$-measurable random variable, $f^i: \Omega \times [0, T]
\times \mathbb{R} \times \mathbb{R}^d \rightarrow \mathbb{R}$ is
$\Prog\otimes \mathcal{B}({\mathbb{R}})\otimes \mathcal{B}({\mathbb{R}^d})$-measurable and \[
\norm{\xi^i}^2_{\infty}+8 \big\|{\big||f^i(\cdot, 0, 0)|\big|_T}\big\|^2_{\infty} \leq \frac{1}{64}\exp \Big(       -\norm{A} 
\big(8\beta^2 \norm{{A}} + 8\gamma^2\big)
    \Big).
\]
Set $f^\prime  (t, y, z): = f(t, y, z)- f(t, 0, 0) $ and $(Y^0, Z^0\cdot M + N^0) \in \mathscr{B}$ such that $\norm{(Y^0, Z^0 \cdot M + N^0)} =0$. 
Now we use a recursion argument in the following way for $i=1, ..., n$.

By Theorem \ref{existence1}, there exists a solution $(Y^i, Z^i\cdot M + \widetilde{N}^{i}) \in \mathscr{B}(\mathbb{Q}^i)$ to the BSDE
\begin{align*}
Y_t^i &= \xi^i + \int_t^T \Big(f^i(s, 0, 0)+f^\prime(s, \sum_{j=0}^{i} Y_s^j, \sum_{j=0}^{i} Z_s^j) -  f^\prime(s, \sum_{j=0}^{i-1} Y_s^j, \sum_{j=0}^{i-1} Z_s^j) \Big)dA_s  \\ 
&+ \int_t^T g_s d \langle \widetilde{N}^{i}\rangle_s  -\int_t^T \big(Z_s^i dM_s +d\widetilde{N}_s^{i} \big),
\end{align*}
where 
\[
\frac{d\mathbb{Q}^i}{d\mathbb{P}}: = \mathcal{E}\Big(2g\cdot \sum_{j=0}^{i-1}N^j\Big)_T.
\]
Note that 
 the equivalent change of measure holds due to the fact that  $N^j \in \mathcal{M}^{BMO}$ for  $j \leq i-1$ and Theorem 2.3,  Kazamaki \cite{K1994}.
By Girsanov transformation and Theorem 3.6,  Kazamaki \cite{K1994},   $N^i := \widetilde{N}^{ i } + 2g \cdot \langle \widetilde{N}^{ i}, \sum_{j=0}^{i-1}N^j \rangle $ and $Z^i\cdot M$ belong to $\mathcal{M}^{BMO}$. This  further implies  $\langle N^i \rangle = \langle \widetilde{N}^i\rangle$ and 
$N^i = \widetilde{N}^{ i } + 2g \cdot \langle N^{ i}, \sum_{j=0}^{i-1}N^j \rangle .$ Hence $(Y^{i}, Z^{i}\cdot M + N^{i})\in \mathscr{B} $ solves
\begin{align*}
Y_t^i &= \xi^i + \int_t^T \Big(f^i (s, 0, 0 ) + f^\prime(s, \sum_{j=0}^{i} Y_s^j, \sum_{j=0}^{i} Z_s^j) -  f^\prime(s, \sum_{j=0}^{i-1} Y_s^j, \sum_{j=0}^{i-1} Z_s^j) \Big)dA_s  \\ 
&+ \int_t^T g_s d \Big(\langle N^i\rangle_s + 2\langle N^i, \sum_{j=0}^{i-1} N^j\rangle_s  \Big)  -\int_t^T \big(Z_s^i dM_s + dN_s^i\big ).
\end{align*}
Hence a recursion argument gives $(Y^i, Z^i, N^i)$ for $i =  1, ..., n$.

Define ${Y}: = \sum_{i=0}^n Y^i, Z:= \sum_{i=0}^n Z^i$ and $N:= \sum_{i=0}^n N^i$. 
 Clearly $(Y, Z\cdot M + N) \in \mathscr{B}$.
We show $(Y, Z, N)$ solves $(f, g, \xi)$.
In view of the definition of $f^\prime$, we 
sum up  the above BSDEs  to  obtain
\[
{Y}_t = \xi + \int_t^T \Big(\big(f(s, 0, 0) + f^\prime(s, {Y}_s, Z_s)\big)dA_s + g_s d\langle N \rangle_s\Big) -\int_t^T  \big( \delta Z_s dM_s + d\delta N_s\big).
\]
To conlcude the proof we use $f^\prime  (s, Y_s, Z_s): = f(s, Y_s, Z_s)- f(s, 0, 0) $.
\qed
\end{Proof}

We continue to show that  comparison theorem  and hence uniqueness also hold given Lipschitz-continuity. Similar results in different settings can be found, e.g., in 
\ \cite{MS2005}, \cite{HIM2005},  \cite{M2009}, \cite{T2008}.
\begin{thm}[Comparison]\label{compare1}
Let $(Y, Z\cdot M + N)$, $(Y^\prime, Z^\prime\cdot M+ N^\prime) \in \mathcal{S}^\infty \times \mathcal{M}^{BMO}$ be solutions to $(f, g, \xi)$, $(f^\prime, g^\prime, \xi^\prime)$, respectively. If 
$\mathbb{P}$-a.s. for any $(t, y, z) \in [0, T]\times\mathbb{R}\times\mathbb{R}^d$, 
$f(t, y ,z)\leq f^\prime(t, y, z)$, $g_t\leq g_t^\prime$, $\xi \leq \xi^\prime$ and $(f,g, \xi)$ verifies  {\rm\ref{as1}}, then $\mathbb{P}$-a.s. $Y_\cdot\leq Y^\prime_\cdot$.
\end{thm}
\begin{Proof}
 Set
$\delta Y := Y- Y^\prime$ and $\delta Z, \delta N, \delta \langle N\rangle, \delta \xi$, etc. analogously. For any $\tau\in\mathcal{T}$, $\mathbb{P}$-a.s. $f\leq f^\prime$ and $g_\cdot\leq g^\prime_\cdot$ imply by It\^{o}'s formula that 
\begin{align}
\delta Y_{t\wedge \tau} &= \delta Y_\tau + \int_{t\wedge \tau}^\tau \big(f(s, Y_s, Z_s)- f^\prime (s, Y_s^\prime, Z_s^\prime)\big)dA_s +\int_{t\wedge \tau}^\tau g_s d\langle N\rangle_s 
-\int_{t\wedge \tau}^\tau g^\prime_s d\langle N^\prime\rangle_s 
\nonumber\\
&-\int_{t\wedge \tau}^\tau \big(\delta Z_s dM_s + d\delta N_s\big ) \nonumber \\
& \leq  \delta Y_\tau + \int_{t\wedge \tau}^\tau\big( f(s, Y_s, Z_s) - f(s, Y^\prime_s, Z^\prime_s) \big)dA_s
  + \int_{t\wedge \tau}^\tau g^\prime_s d\delta \langle N\rangle_s  -\int_{t\wedge \tau}^\tau \big(\delta Z_s dM_s + d\delta N_s\big)\nonumber\\
&  =  \delta Y_\tau + \int_{t\wedge \tau}^\tau \big(\beta_s \delta Y_s + (\lambda_s \gamma_s)^\top (\lambda_s\delta Z_s) \big)dA_s 
  + \int_{t\wedge \tau}^\tau g^\prime_s d\delta \langle N\rangle_s  -\int_{t\wedge \tau}^\tau \big(\delta Z_s dM_s + d\delta N_s\big),\label{lqunique}
\end{align}
where  $\beta$ ($\mathbb{R}$-valued) and  $\gamma$ ($\mathbb{R}^d$-valued) are defined by
\begin{align*}
\beta_s &: =\mathbb{I}_{\{\delta Y_s \neq 0 \}} \frac{f(s, Y_s, Z_s) - f(s, Y_s^\prime, Z_s)}{\delta Y_s},\\
\gamma_s &: = \mathbb{I}_{\{\lambda_s\delta Z_s \neq \mathbf{0}\}}\frac{\big(f(s, Y_s^\prime, Z_s)-f(s, Y_s^\prime, Z_s^\prime)\big)\delta Z_s}{|\lambda_s \delta Z_s|^2},
\end{align*}
and  $\mathbf{0}:=(0, ..., 0)^\top.$ Note that $\gamma$ can be seen as defined in terms of discrete gradient. 
By \ref{as1}, $\beta_\cdot$ and $\int_0^\cdot \gamma^\top_s d\langle M \rangle_s \gamma_s$ are bounded processes,  hence $\gamma \cdot M \in \mathcal{M}^{BMO}.$ 
Given these facts we  use a change of measure to attain the comparison result. 
To this end, we define a BMO martingale
\[
\Lambda : = 
\gamma \cdot M  + g^\prime \cdot (N+N^\prime).
\]
In view of  Theorem 2.3 and Theorem 3.6, Kamazaki \cite{K1994},  we  define 
\[
\frac{d\mathbb{Q}}{d\mathbb{P}}: = \mathcal{E} ( \Lambda)_T.
\]
Hence
$ \delta N - g^\prime \cdot \delta \langle N \rangle$ and 
$\delta Z \cdot M- ( \gamma^\top\lambda^\top \lambda \delta Z) \cdot A $  belong to $\mathcal{M}^{BMO}(\mathbb{Q})$.
Therefore, (\ref{lqunique})  and $\mathbb{P}$-a.s. $\delta \xi \leq 0$ give
\begin{align*}
\delta Y_{t} &\leq \mathbb{E}^{\mathbb{Q}} \big[     \delta \xi\big|\mathcal{F}_{t} \big] + \mathbb{E}^{\mathbb{Q}}
\Big[ \int_{t}^{ T} \beta_s \delta Y_s dA_s \Big| \mathcal{F}_{t}\Big]  \\
&\leq     \mathbb{E}^{\mathbb{Q}}
\Big[ \int_{t}^{ T} \beta_s \delta Y_s dA_s \Big| \mathcal{F}_{t}\Big].
\end{align*}
Hence we obtain by Gronwall's lemma that  $\mathbb{P}$-a.s. $\delta Y_t \leq 0$.
Finally by the continuity of $Y$ and $Y^\prime$, we conclude that $\mathbb{P}$-a.s. $Y_\cdot \leq Y_\cdot^\prime.$
\qed
\end{Proof}

As a byproduct, we obtain the following existence and uniqueness result.
\begin{Corollary} [Uniqueness] \label{unique1}If $(f, g, \xi)$  satisfies  {\rm{\ref{as1}}}, then there exists a unique solution
in $\mathscr{B}$.
\end{Corollary}
\begin{Proof}
This is immediate from 
Theorem \ref{existence2} (existence (ii)) and 
 Theorem \ref{compare1} (comparison theorem). 
\qed
\end{Proof}
\section{Bounded Solutions to Quadratic BSDEs}
\label{section22}
In this section, we   prove a general  monotone stability result for quadratic BSDEs. 
Let us recall that Morlais \cite{M2009} uses  a stability-type argument for the existence result
   after performing an exponential transform which eliminates $g\cdot \langle N\rangle$.
But a  direct stability result is not studied. Our work fills this gap. 

Secondly, as a byproduct, we construct a bounded solution via 
  regularization through Lipschitz-quadratic BSDEs studied in Section \ref{section22}.  This  procedure is also called \emph{ Lipschitz-quadratic regularization} in the following context. 
  To this end we give the assumptions for the whole section. 
\begin{as}\label{as2}
There exist $\beta \geq 0$, $\gamma >0$, an $\mathbb{R^+}$-valued $\Prog$-measurable process $\alpha$  and a continuous nondecreasing function $\varphi : \mathbb{R}^+ \rightarrow \mathbb{R}^+$ with $\varphi(0) = 0$ such that
${\norm{\xi}}_\infty  + \norm{|\alpha|_T}_\infty < +\infty$
and $\mathbb{P}$-a.s.
\begin{enumerate}
\item [(i)] for any $t\in[0, T]$, $(y, z)\longmapsto f(t, y, z)$ is continuous;
\item [(ii)]$f$ is monotonic at $y=0$, i.e., for any $(t, y, z)\in[0, T]\times\mathbb{R}\times\mathbb{R}^d$,
\[ \sgn (y) f(t, y, z) \leq \alpha_t + \alpha_t\beta |y|+ \frac{\gamma}{2}|\lambda_t z|^2;\]
\item [(iii)]  for any $(t, y, z)\in[0, T]\times\mathbb{R}\times\mathbb{R}^d$,
 \[ |f(t, y, z)| \leq \alpha_t + \alpha_t\varphi(|y|) + \frac{\gamma}{2}|\lambda_t z|^2. \]
\end{enumerate}
\end{as}

 We continue as before to call  $(\xi,|\alpha|_T)$ the \emph{data}.
 \ref{as2}(ii) 
  allows one to get rid of the linear growth  in $y$ which is required by  Kobylanski \cite{K2000} and Morlais \cite{M2009}. Assumption of this type for quadratic framework is motivated by Briand and Hu \cite{BH2008}.
Secondly, our results don't rely on the specific choice of $\varphi$.  Hence the growth condition in $y$  can be arbitrary as long as \ref{as2}(i)(ii)   hold. 

Given \ref{as2}, we first prove an a priori estimate. In order to treat $\langle Z\cdot M\rangle$ and $g \cdot \langle N \rangle$ more easily, we assume  $\mathbb{P}$-a.s. $|g_\cdot|\leq \frac{\gamma}{2}$ for the rest of this chapter. 
\begin{Lemma}[A Priori Estimate]\label{aprioriestimate1} If $(f, g, \xi)$  satisfies  {\rm\ref{as2}} and $(Y, Z\cdot M + N) \in \mathcal{S}^\infty \times \mathcal{M}$  is a solution to $(f, g, \xi)$, then
\begin{align*}
\norm{Y} \leq \big\| e^{\beta |\alpha|_T}\big (|\xi| +|\alpha|_T\big)\big\|_\infty \end{align*}
and
\begin{align*}
&\norm {Z\cdot M + N }_{BMO} \leq c_{b},
\end{align*}
where $c_b$ is a constant only depending on $\beta, \gamma, \norm{\xi}_\infty, \norm{|\alpha|_T}_\infty$.
\end{Lemma}
\begin{Proof}  Set
 $u(x): = \frac{\exp({\gamma x})-1 - \gamma x}{\gamma^2}.$ The following auxiliary results will be useful:  
 $u(x) \geq 0, u^\prime(x) \geq 0$ and  $u^{\prime\prime}(x) \geq 1$ for $x\geq 0$;  $u(|\cdot|)\in \mathcal{C}^2(\mathbb{R})$ and $u^{\prime\prime}(x) = \gamma u^\prime(x) + 1$. 
For any  $\tau, \sigma \in \mathcal{T}$, It\^{o}'s formula yields
\begin{align*}
u(|Y_{\tau\wedge \sigma}|) =& u(|Y_{\sigma}|) + \int_{\tau \wedge \sigma}^\sigma u^\prime (|Y_s|)\sgn(Y_s)dY_s - \frac{1}{2}\int_{\tau\wedge \sigma}^\sigma u^{\prime\prime} (|Y_s|)\Big(Z_s^\top d\langle M\rangle_s Z_s + d\langle N\rangle_s\Big).
\end{align*}
By  \ref{as2}(ii),\[ 
\sgn(Y_s) f(s, Y_s, Z_s) \leq \alpha_s +\alpha_s \beta |Y_s|+ \frac{\gamma}{2}|\lambda_s Z_s|^2.\]
Note that $
\frac{\gamma}{2}u^\prime (|Y_s|)  -\frac{1}{2} u^{\prime\prime}(|Y_s|) = - \frac{1}{2}$,  $ 
{g_s}u^\prime (|Y_s|)  -\frac{1}{2} u^{\prime\prime}(|Y_s|) \leq - \frac{1}{2}.
$ and $u^\prime(|Y_s|) \leq \frac{e^{\gamma\norm{Y}}}{\gamma}$.
Hence, using these facts to the above equality yields
\begin{align*}
\frac{1}{2}\int_{\tau\wedge \sigma}^\sigma \Big( Z_s^\top d\langle M\rangle_s  Z_s + d\langle N \rangle_s \Big)&\leq 
\frac{e^{\gamma \norm{Y}}}{\gamma^2} +  \int_{\tau\wedge \sigma}^\sigma u^\prime (|Y_s|)\big(\alpha_s + \alpha_s \beta |Y_s|\big)dA_s \nonumber\\
& - \int_{\tau\wedge \sigma}^{ \sigma} u^{\prime}(|Y_s|)\sgn (Y_s) \big(Z_s dM_s + dN_s\big).  
\end{align*} 
To eliminate the local martingale,
we replace $\sigma$ by its localizing sequence and use Fatou's lemma to the left-hand side. Since  $Y^*$ and $|\alpha|_T$ are bounded random variables,  the right-hand side has a uniform constant  upper bound. Hence, we have 
\begin{align}
\frac{1}{2} \mathbb{E}\big[
\langle Z\cdot M  + N \rangle_{\tau, T} 
 \big|\mathcal{F}_\tau \big] \leq 
 \frac{e^{\gamma \norm{Y}}}{\gamma^2} + 
  \frac{e^{\gamma \norm{Y}}}{\gamma}(1+\beta \norm{Y}) \norm{|\alpha|_T}_\infty.
  \label{bmo}
\end{align}
 Now we turn to the estimate for $Y$. We fix $s\in [0, T]$ and for $t\in [s, T]$, set
\[
H_t : = \exp \Big(  \gamma e^{\beta |\alpha|_{s, t}}|Y_t| + \gamma \int_s^t e^{\beta |\alpha|_{s, u}}\alpha_u dA_u   \Big).
\]
We claim that $H$ is a submartingale. 
By Tanaka's formula,
\[
d|Y_t| = \sgn (Y_t) \big(Z_t dM_t + dN_t\big) -\sgn(Y_t)\big(f(t, Y_t, Z_t)dA_t + g_t d\langle N\rangle_t\big) + dL_t^0(Y),
\]
where $L^0(Y)$ is the local time of $Y$ at $0$. Hence, It\^{o}'s formula yields
\begin{align*}
dH_t &= \gamma H_t e^{\beta |\alpha|_{s, t}}\Big[  \sgn (Y_t)\big (Z_t dM_t + dN_t\big) \\
& + \Big(- \sgn(Y_t)f(t, Y_t, Z_t) +\alpha_t + \alpha_t \beta |Y_t| + \frac{\gamma}{2}e^{\beta |\alpha|_{s, t}}|\lambda_t Z_t|^2\Big)dA_t  \\
&+ \Big(- \sgn(Y_t)g_t +\frac{\gamma}{2}e^{\beta |\alpha|_{s, t}}\Big)d\langle N\rangle_t + dL_t^0(Y) 
\Big].
\end{align*}
By  \ref{as2}(ii) and  $ |g_\cdot|\leq \frac{\gamma}{2}$ again,  $(H_t)_{t\in[s, T]}$ is a bounded submartingale.  Hence, 
\begin{align*}
|Y_s| \leq \frac{1}{\gamma}\ln \mathbb{E} \big[H_T \big| \mathcal{F}_s\big].
\end{align*}
Thanks to the boundedness, we have 
\begin{align*}
\norm{Y} \leq \big\| e^{\beta |\alpha|_T} \big(|\xi| +|\alpha|_T\big)\big\|_\infty.
\end{align*}
Finally we  come back to (\ref{bmo}) 
and obtain the estimate for $Z\cdot M + N$.
\qed
\end{Proof}

Given the norm bound in Lemma \ref{aprioriestimate1}, 
we  turn to the main result of this section: monotone stability result.  
 Later, as an immediate application,  we prove an existence result for quadratic BSDEs by Lipschitz-quadratic regularization. To start we recall that $\mathcal{M}^2$ equipped with the norm $\norm{\widetilde{M}}_{\mathcal{M}^2}:=\mathbb{E}\big[\langle\widetilde{M}\rangle_T\big]^\frac{1}{2}$ for $\widetilde{M}\in\mathcal{M}^2$ is a Hilbert space. 
\begin{thm}[Monotone Stability] \label{monotonestability}
Let $(f^n, g^n, \xi^n)_{n\in\mathbb{N}^+}$  satisfy  {\rm\ref{as2}} associated with $(\alpha, \beta, \gamma, \varphi)$, and  $(Y^n, Z^n \cdot M + N^n)$ be their solutions in $\mathscr{B}$, respectively.  Assume
\begin{enumerate}
\item [{\rm(i)}]
 $Y^n$  is monotonic in $n$ and  $
\xi^n - \xi \longrightarrow 0           \ \mathbb{P}$-a.s. with $\sup_n \norm{\xi^n}_\infty< +\infty${\rm;} 
\item [{\rm(ii)}]$\mathbb{P}$-a.s. for any $t\in [0, T]$, $g^n_t - g_t \longrightarrow 0${\rm;}  
\item [{\rm(iii)}] $\mathbb{P}$-a.s. for any $t\in[0, T]$ and $y^n \longrightarrow y, z^n \longrightarrow z$,
$f^n(t, y^n, z^n)\longrightarrow f(t, y, z)$.
\end{enumerate}

Then there exists $(Y, Z\cdot M + N) \in \mathscr{B}$ such that $Y^n$ converges to $Y$ $\mathbb{P}$-a.s. uniformly on $[0, T]$ and $(Z^n \cdot M + N^n)$ converges to $(Z\cdot M + N)$ in  $\mathcal{M}^2$ as $n$ goes to $+\infty$. Moreover, $(Y, Z, N)$ solves $(f, g, \xi)$.

\end{thm}
\begin{Proof}
Without loss of generality we only consider  $Y^n$ to be increasing in $n$.
By Lemma \ref{aprioriestimate1} (a priori estimate),
\begin{align}
\sup_n  \norm{Y^n}+\sup_n \norm{Z^n \cdot M + N^n}_{BMO}  \leq c_b, \label{e31}
\end{align}
where $c_b$ is a constant  only depending on $\beta, \gamma, \sup_n{\norm{\xi^n}}_\infty, \norm{|\alpha|_T}_\infty$. We  rely intensively on the boundedness  result in (\ref{e31})  to derive the limit.

(i). We prove  the  convergence of the solutions. 
 Due to (\ref{e31}),
there exists a bounded monotone limit 
$
Y_t:= \lim_{n} Y_t^n,
$   a subsequence 
 indexed by $\{n_k\}_{k\in\mathbb{N}^+} \subseteq \mathbb{N}^+$ 
 and $Z\cdot  M +N \in \mathcal{M}^2$
 such that 
 $Z^{n_k}\cdot M + N^{n_k}$ converges weakly 
in $\mathcal{M}^2$
 to $Z\cdot M + N$
as $k$ goes to $+\infty.$ 
 The task is to show $Z\cdot M + N$  is the $\mathcal{M}^2$-limit of the whole sequence.  To this end we  define
$
u(x) := \frac{\exp({8\gamma x})-8\gamma x -1}{64\gamma^2}.
$
Recall that $u(x) \geq 0, u^\prime (x)\geq 0$ and $u^{\prime\prime}(x) \geq 0$ for $x\geq 0$; $u \in \mathcal{C}^2(\mathbb{R})$ and $u^{\prime\prime}(x)= 8\gamma u^\prime(x)+1$.  For any $  m\in \{n_k\}_{k\in\mathbb{N}^+}$ ,$n \in \mathbb{N}^+$ with $m \geq n$, define $\delta Y^{m, n}: = Y^m - Y^n, \delta Y^{n}: = Y - Y^n$ and $\delta Z^{m, n}, \delta Z^n, \delta N^{m, n}, \delta N^n$, etc. analogously.
By It\^{o}'s formula, 
\begin{align}
\mathbb{E}\big[u(\delta Y_0^{m, n})\big]- \mathbb{E}\big[u(\delta \xi^{m, n})\big]
&=\mathbb{E}\Big[\int_0^T  u^\prime(\delta Y_s^{m, n})\big(f^m (s, Y^m_s, Z^m_s)- f^n (s, Y^n_s, Z^n_s)\big)dA_s\Big]
\nonumber\\
&+ \mathbb{E}\Big[\int_0^T u^\prime (\delta Y_s^{m, n}) \big(g_s^m d\langle N^m \rangle_s  
-  g_s^{n} d \langle N^n \rangle_s \big)\Big]\nonumber \\
&-\frac{1}{2}\mathbb{E}\Big[  \int_0^T u^{\prime\prime}(\delta Y_s^{m, n} )  \Big( (\delta Z_s^{m, n})^\top d\langle M\rangle_s (\delta Z_s^{m, n})+ d\langle \delta N^{m, n}\rangle_s \Big) \Big]. \label{e32}
\end{align}
Since $f^m$ and $f^{n}$ verify \ref{as2} associated with $(\alpha, \beta, \gamma, \varphi)$, we have
\begin{align*}
|f^m (s, Y^m_s, Z^m_s) &- f^n (s, Y^n_s, Z^n_s)|  \nonumber\\
&\leq \alpha_s^\prime +\frac{\gamma}{2}|\lambda_s Z^m_s|^2 + \frac{\gamma}{2}|\lambda_s Z^n_s|^2 \nonumber\\
& \leq \alpha_s^\prime+ \frac{3\gamma}{2}\big(|\lambda_s \delta Z_s^{m, n}|^2 + |\lambda_s \delta Z^n_s|^2 + |\lambda_s Z_s|^2\big) + \gamma \big(|\lambda_s \delta Z_s^n|^2 + |\lambda_s Z_s|^2\big) \nonumber \\
& \leq \alpha_s^\prime +\frac{3\gamma}{2}|\lambda_s \delta Z^{m, n}_s|^2 +  \frac{5\gamma}{2} \big( |\lambda_s \delta Z^n_s |^2 + |\lambda_s Z_s|^2 \big),
\end{align*}
where
\[
\alpha_s^\prime :  = 2\alpha_s\big(1+\varphi(c_b) \big) \geq   2\alpha_s + \alpha_s\varphi(|Y^n_s|)+ \alpha_s\varphi(|Y^m_s|).
\]
Moreover,
\begin{align*}
g^m d\langle N^m \rangle - g^n d\langle N^n \rangle  &\ll
\frac{\gamma}{2}d\langle N^m \rangle + \frac{\gamma}{2} d\langle N^n \rangle \\
&\ll
\frac{3\gamma}{2}d\langle \delta N^{m , n}\rangle + \frac{5\gamma}{2}\big(d\langle \delta N^n \rangle
+ d\langle N \rangle  \big).
\end{align*}
Plugging the above inequalities into (\ref{e32}), we deduce that 
\begin{align}
&\mathbb{E}\Big[ \int_0^T  \Big(\frac{1}{2}u^{\prime\prime} -\frac{3\gamma}{2}u^\prime\Big) (\delta Y_s^{m, n})  |\lambda_s \delta Z_s^{m, n}|^2 dA_s	\Big]
+ \mathbb{E}\Big[  \int_0^T \Big(\frac{1}{2} u^{\prime\prime}-\frac{3\gamma}{2}u^\prime\Big)(\delta Y_s^{m, n})   d\langle\delta N^{m, n}\rangle_s     \Big]\nonumber \\
&\leq
\mathbb{E}\big[u(\delta \xi^{m, n})\big] + \mathbb{E}\Big[\int_0^T u^\prime (\delta Y^{m, n}_s) \Big( \alpha_s^\prime +\frac{5\gamma}{2}\big(|\lambda_s \delta Z_s^{n}|^2  + |\lambda_s Z_s|^2 \big)  \Big)  dA_s \Big]\nonumber\\
&+  \mathbb{E}\Big[\int_0^T u^\prime (\delta Y_s^{m, n}) \frac{5\gamma}{2} \big(d\langle \delta N^n \rangle_s + d\langle N \rangle_s    \big) \Big]
 \label{e34} 
\end{align}
Due to the weak convergence result and  convexity of $z\longmapsto |z|^2$, $N\longmapsto \langle N\rangle$, we obtain
\begin{align*}
\mathbb{E}\Big[ \int_0^T  \Big(\frac{1}{2}u^{\prime\prime} -\frac{3\gamma}{2}u^\prime\Big) (\delta Y_s^{ n})  |\lambda_t Z_s^{ n}|^2 dA_s	\Big] &\leq \liminf_{m} 
\mathbb{E}\Big[ \int_0^T  \Big(\frac{1}{2}u^{\prime\prime} -\frac{3\gamma}{2}u^\prime\Big) (\delta Y_s^{m, n})  |\lambda_t Z_s^{ m, n}|^2 dA_s	\Big],\\
 \mathbb{E}\Big[  \int_0^T \Big(\frac{1}{2}u^{\prime\prime}-\frac{3\gamma}{2}u^\prime \Big)(\delta Y_s^{ n})   d\langle\delta N^{ n}\rangle_s     \Big]&\leq
 \liminf_m\mathbb{E}\Big[  \int_0^T \Big(\frac{1}{2}u^{\prime\prime}-\frac{3\gamma}{2}u^\prime\Big)(\delta Y_s^{m, n})   d\langle\delta N^{m, n}\rangle_s     \Big].
\end{align*}
We then come back to 
 (\ref{e34}) and send $m$ to  $+\infty$ along $\{n_k \}_{k\in \mathbb{N}^+}$. Taking the  above inequalities into account and using $u^{\prime}(\delta Y^{m, n}_s)\leq u^{\prime}(\delta Y^{n}_s)$
 to the right-hand side, 
 (\ref{e34}) becomes
\begin{align}
&\mathbb{E}\Big[ \int_0^T  \Big(\frac{1}{2}u^{\prime\prime} -\frac{3\gamma}{2}u^\prime\Big) (\delta Y_s^{ n})  |\lambda_s Z_s^{ n}|^2 dA_s	\Big]\nonumber\\
&+ \mathbb{E}\Big[\int_0^T \Big(\frac{1}{2}u^{\prime\prime}-\frac{3\gamma}{2}u^{\prime}\Big)(\delta Y_s^n) d\langle \delta N^n \rangle_s  \Big]
\nonumber\\&\leq \mathbb{E}\big[u(\delta \xi^{ n})\big]  + \mathbb{E}\Big[\int_0^T u^\prime (\delta Y^{ n}_s) \Big( \alpha^\prime_s +\frac{5\gamma}{2}\big(|\lambda_s \delta Z_s^{n}|^2  + |\lambda_s Z_s|^2 \big)  \Big)  dA_s \Big]\nonumber\\
&+ \frac{5\gamma}{2}\mathbb{E}\Big[\int_0^T u^\prime (\delta Y_s^{n})  \big(d\langle \delta N^n \rangle_s + d\langle N \rangle_s    \big) \Big].
\label{e35}
\end{align}
Since $u^{\prime\prime} (x) - 8\gamma u^\prime (x) =1$,   rearranging terms give
\begin{align}
&\frac{1}{2}E\big[ \big(\delta N_T^n\big)^2 \big]  + \frac{1}{2}\mathbb{E}\Big[ \int_0^T |\lambda_s \delta Z_s^n|^2 dA_s\Big]
 \nonumber\\
&\leq \mathbb{E}\big[u(\delta \xi^{ n})\big]+ \mathbb{E}\Big[\int_0^T u^\prime (\delta Y^{ n}_s) \Big( \alpha^\prime_s +\frac{5\gamma}{2} |\lambda_s Z_s|^2   \Big)  dA_s \Big]+ \frac{5\gamma}{2}\mathbb{E}\Big[ \int_0^T u^{\prime}(\delta Y_s^n)d\langle N \rangle_s\Big].
 \label{e36}
\end{align}
Finally by sending $n$ to $+\infty$ and dominated convergence we deduce the convergence.

(ii).
We  prove $(Y, Z\cdot M + N)\in \mathscr{B}$ and solves $(f, g, \xi)$.
Here we rely on the same arguments as in  Kobylanski \cite{K2000} or  Morlais \cite{M2009} and omit the details here.
In addition, we need to prove the $u.c.p$ convergence of  $g^n\cdot \langle N^n \rangle$, which holds if 
\begin{align*}
\lim_{n\rightarrow \infty} \mathbb{E}\Big[ \Big| \int_0^\cdot \big(g^n_s d\langle N^n \rangle_s - g_s d\langle N\rangle_s \big)\Big|^* \Big] =0.
\end{align*}
Indeed, by Kunita-Watanabe inequality and Cauchy-Schwartz inequality, 
\begin{align*}
 \mathbb{E}\Big[ \Big| \int_0^\cdot \big(g^n_s d\langle N^n \rangle_s &- g_s d\langle N \rangle_s \big) \Big|^* \Big]
= \mathbb{E}\Big[ \Big| \int_0^\cdot \Big(g^n_s d\big(\langle N^n \rangle_s - \langle N \rangle_s \big)+ ( g^n_s-g_s )d\langle N \rangle_s \Big) \Big|^* \Big]\\
&\leq \frac{\gamma}{2} \mathbb{E}\big[     \langle N^n - N\rangle_T  \big]^{\frac{1}{2}} \mathbb{E}\big[     \langle N^n + N\rangle_T  \big]^{\frac{1}{2}} + \mathbb{E}\Big[ \Big|  \int_0^\cdot (g^n_s -g_s)d\langle N\rangle_s \Big|^*   \Big]\\
&\leq  \gamma c_b \mathbb{E}\big[     \langle N^n - N\rangle_T  \big]^{\frac{1}{2}} + \mathbb{E}\Big[  \int_0^T |g^n_s -g_s|d\langle N\rangle_s   \Big].
\end{align*}
We then conclude by  $\mathcal{M}^2$-convergence of $N^n$ and dominated convergence used to the second term. Finally 
$Z\cdot M + N \in\mathcal{M}^{BMO}$ by Lemma \ref{aprioriestimate1} (a priori estimate).

For decreasing $Y^n$, we take $m \in \mathbb{N}^+, n \in \{n_k\}_{k\in \mathbb{N}^+}$ with $n \geq m$ and conclude with exactly the same arguments.
\qed
\end{Proof}

There are several major improvements compared to existing monotone stability results. First of all, in contrast to  Kobylanski \cite{K2000} and  Morlais \cite{M2009}, we get rid of linear growth in $y$ by merely assuming \ref{as2}, and allow $g$ to  be any  bounded process.
Secondly, we treat the convergence in a more direct and general way than Morlais \cite{M2009}.

Another advantage concerns the existence result. 
 Thanks to Section \ref{section21} and Theorem \ref{monotonestability},
 we are able to perform a Lipschitz-quadratic regularization where exponential transform to eliminate
   $g\cdot \langle N\rangle$
 is no longer needed; this is in contrast to Morlais \cite{M2009}.
This  also helps to prove the  existence of unbounded solutions with fewer assumptions; see Section \ref{section23}.

\begin{prop}[Existence]\label{existence3}
If $(f, g, \xi)$ satisfy  {\rm\ref{as2}}, then there exists a solution in $\mathscr{B}$.
\end{prop}
\begin{Proof}
We use a double approximation procedure and use Theorem \ref{monotonestability}  (monotone stability) to take the limit. 
 Define
\begin{align*}
f^{n, k}(t, y, z): &= \inf_{y^\prime, z^\prime} \big\{f^+ (t, y^\prime, z^\prime) + n|y-y^\prime| + n|\lambda_t(z-z^\prime)|\big \}\\
&- \inf_{y^\prime, z^\prime}\big\{ f^- (t, y^\prime, z^\prime) + k|y-y^\prime| + k|\lambda_t(z-z^\prime)|\big\}.
\end{align*}
By
Lepeltier and  San Martin \cite{LS1997}, $f^{n, k}$ is Lipschitz-continuous in $(y, z)$; as $k$  goes to $+\infty$, $f^{n, k}$ converges increasingly   uniformly on compact sets to a limit denoted by $f^{n, \infty}$;  as $n$  goes to $+\infty$, $f^{n, \infty}$ converges increasingly uniformly on compact sets to $f$.

By Corollary \ref{unique1}, there exists a unique solution $(Y^{n ,k}, Z^{n, k}\cdot M + N^{n, k})\in\mathscr{B}$ to $(f^{n, k}, g, \xi)$; by  Theorem \ref{compare1} (comparison theorem), $Y^{n, k}$ is increasing in $n$ and decreasing in $k$, and is uniformly bounded due to
Lemma \ref{aprioriestimate1} (a priori estimate).  We then fix $n$ and use  Theorem \ref{monotonestability} to the sequence indexed by $k$ to obtain a solution
 $(Y^{n},Z^{n}\cdot M + N^{n} ) \in$ $\mathscr{B}$ to $(f^{n, \infty}, g, \xi)$. Due to the $\mathbb{P}$-a.s. uniform convergence of $Y^{n, k}$ we  can pass the comparison property to $Y^{n}$. We use Theorem \ref{monotonestability} again to conclude. 
\qed
\end{Proof}
\begin{remm}
In contrast to  Kobylanski \cite{K2000}, the existence of a maximal or minimal solution is not available (yet) given \ref{as1} as the double approximation procedure makes the comparison between solutions impossible.
\end{remm}

There is also a rich literature on the uniqueness of a bounded solution to quadratic BSDEs; see, e.g.,  \cite{K2000}, \cite{MS2005}, \cite{HIM2005},  \cite{M2009}.  Roughly speaking,
they essentially  rely
a type of locally Lipschitz-continuity and use a change of measure  analogously to Section \ref{section21}. 
 The proof in our setting is exactly the same and hence  omitted to save pages. 

To end this section we briefly present various structure conditions used in different situations.

\begin{asp}{as2} 
 \label{as2'} 
 There exist $\beta \geq 0, \gamma >0$, an $\mathbb{R}^+$-valued $\Prog$-measurable process $\alpha$, 
 and a continuous nondecreasing function $\varphi : \mathbb{R}^+ \rightarrow \mathbb{R}^+$ with $\varphi(0) = 0$ such that
$\mathbb{P}$-a.s.  
\begin{enumerate}
\item [(i)] for any $t\in [0, T]$, $(y, z)\longmapsto f(t, y, z)$ is continuous;
\item [(ii)] $f$ is monotonic  at $y=0$, i.e., for any $(t, y, z)\in[0, T]\times\mathbb{R}\times\mathbb{R}^d$,
\begin{align*}\sgn (y) f(t, y, z) \leq \alpha_t + \beta |y|+ \frac{\gamma}{2}|\lambda_t z|^2;\end{align*}
\item [(iii)]  for any $(t, y, z)\in[0, T]\times\mathbb{R}\times\mathbb{R}^d$,
\[
|f(t, y, z)| \leq \alpha_t + \varphi(|y|)+ \frac{\gamma}{2}|\lambda_t z|^2. \] 
\end{enumerate}
\end{asp}

Given bounded data, \ref{as2'} implies \ref{as2}. Indeed, 
\begin{align*}
\sgn(y) f(t, y, z )&\leq \alpha_t \vee 1 +  (\alpha_t \vee 1)\beta |y| + \frac{\gamma}{2}|\lambda_t z|^2, \\
|f(t, y, z)|  &\leq \alpha_t \vee 1 + (\alpha_t\vee 1)\varphi(|y|)+ \frac{\gamma}{2}|\lambda_t z|^2.
\end{align*}
Hence \ref{as2'} verifies  \ref{as2} associated with $(\alpha \vee 1, \beta, \gamma, \varphi)$. 
However, given unbounded data,   \ref{as2'}
appears to be more natural and convenient.
 This will be discussed in detail in Section \ref{section23}.

In  particular situations where the estimate for $\int_0^T |f(s, Y_s, Z_s)|dA_s$ is needed,  e.g., in analysis of measure change  (see  Section \ref{section24})  or the montone stability of quadratic semimartingales
(see Chapter \ref{qs}), there has to be a linear growth in $y$, i.e., 
\begin{aspp}{as2}
\label{as2''}
There exist  $\beta \geq 0$, $\gamma >0$, an $\mathbb{R}^+$-valued $\Prog$-measurable process $\alpha$ such that $\mathbb{P}$-a.s.
\begin{enumerate}
\item [(i)] for any $t\in [0, T]$, $(y, z)\longmapsto f(t, y, z)$ is continuous;
\item [(ii)]
for any $(t, y, z)\in[0, T]\times\mathbb{R}\times\mathbb{R}^d$,
\[
|f(t, y, z)| \leq \alpha_t + \beta|y| + \frac{1}{2}|\lambda_t z|^2.
\]
\end{enumerate}
Indeed, \ref{as2''}
 enables one to obtain the estimate for
$
\int_0^T |f(s, Y_s, Z_s)|dA_s
$
via
\[
\int_0^T |f(s, Y_s, Z_s)|dA_s 
\leq |\alpha|_T + \beta \norm{A}Y^*
+ \frac{\gamma}{2}\langle Z\cdot M \rangle_T.
\]
\end{aspp}

\section{Unbounded Solutions to Quadratic BSDEs}
\label{section23}
This  section extends  Section \ref{section21}, \ref{section22} to unbounded solutions. We prove an existence result and later  show that the uniqueness holds given convexity assumption as an additional requirement.
We point out that similar results have been obtained by Mocha and Westray \cite{MW2012}, but our results rely on much fewer  assumptions and are more natural.
Analogously to section \ref{section22}, we give an a priori estimate in the first step.
We keep in mind that  
$\mathbb{P}$-a.s. $ |g_\cdot| \leq \frac{\gamma}{2}$ throughout our study.
\begin{Lemma}[A priori estimate]\label{aprioriestimate2} If $(f, g, \xi)$ satisfies  {\rm\ref{as2'}} and
$(Y, Z\cdot M + N)\in \mathcal{S}\times\mathcal{M}$ is a solution to $(f, g, \xi)$ such that the process 
\[
\exp \Big( \gamma e^{\beta A_T}{ |Y_\cdot|} + \gamma \int_0^T e^{\beta A_s}\alpha_s dA_s    \Big)
\]
is of class $\mathcal{D}$, then 
\begin{align}
|Y_s| \leq \frac{1}{\gamma} \ln \mathbb{E}\Big[ \exp\Big(\gamma e^{\beta A_{s, T}}|\xi| + 
\gamma \int_s^T e^{\beta A_{s, u}}{\alpha_u dA_u} \Big)\Big| \mathcal{F}_s \Big]. \label{ey}
\end{align}
\end{Lemma}
\begin{Proof}
We  fix $s\in [0, T]$, and for $t\in [s, T]$, set
\begin{align}
H_t : = \exp \Big(  \gamma e^{\beta A_{s, t}}|Y_t| + \gamma \int_s^t e^{\beta A_{s, u}}\alpha_u dA_u   \Big). \label{H}
\end{align}
We claim that $H$ is a local submartingale.  Indeed, 
by Tanaka's formula
\[
d|Y_t| = \sgn (Y_t) \big(Z_t dM_t + dN_t \big) -\sgn(Y_t)\big(f(t, Y_t, Z_t)dA_t + g_t d\langle N\rangle_t\big) + dL_t^0(Y),
\]
where $L^0(Y)$ is the local time of $Y$ at $0$. Hence, It\^{o}'s formula yields
\begin{align*}
dH_t &= \gamma H_t e^{\beta A_{s, t} }\Big[  \sgn (Y_t) \big(Z_t dM_t + dN_t \big)  \\
& + \Big(- \sgn(Y_t)f(t, Y_t, Z_t) +\alpha_t +  \beta |Y_t| + \frac{\gamma}{2}e^{\beta  A_{s, t} }|\lambda_t Z_t|^2\Big)dA_t  \\
&+ \Big(- \sgn(Y_t)g_t +\frac{\gamma}{2}e^{\beta A_{s, t} }\Big)d\langle N\rangle_t + dL_t^0(Y) 
\Big].
\end{align*}
By  \ref{as2'}(ii), $H$ is a local submartingale. To eliminate the local martingale part, we replace $\tau$ by its localizing sequence on $[s, T]$, denoted by
 $\{\tau_n\}_{n\in \mathbb{N}^+}$. Therefore, 
\begin{align*}
|Y_s| &\leq \frac{1}{\gamma}\ln \mathbb{E} \big[H_{T\wedge \tau_n} \big| \mathcal{F}_s\big]\\
&\leq \frac{1}{\gamma} \ln \mathbb{E}\Big[ \exp\Big(\gamma e^{\beta A_{s, T\wedge \tau_n}}|Y_{T\wedge \tau_n}| + 
\gamma \int_s^{T\wedge \tau_n} e^{\beta A_{s, u}}{\alpha_u dA_u} \Big)\Big| \mathcal{F}_s \Big]. 
\end{align*}
Finally by class $\mathcal{D}$ property
we conclude by sending $n$ to $+\infty$.
\qed
\end{Proof}

We then know from Lemma \ref{aprioriestimate2} that exponential moments integrability on $|\xi| + |\alpha|_T$ is a natural requirement for the existence result. 
\begin{remm}
 \ref{as2'} addresses the issue of integrability better than
  \ref{as2}.  
To show this, let us assume \ref{as2}. We then deduce from Lemma \ref{aprioriestimate1}  and corresponding class $\mathcal{D}$ property that 
\begin{align}
|Y_s| \leq \frac{1}{\gamma} \ln \mathbb{E}\Big[  \exp\Big(\gamma e^{\beta|\alpha|_{s, T}}|\xi| + \gamma \int_s^T e^{\beta |\alpha|_{s, u}} \alpha_u dA_u     \Big)         \Big| \mathcal{F}_s     \Big]. \label{bound2}
\end{align}
Obviously, in (\ref{bound2}), even exponential moments integrability  
is not sufficient to ensure the well-posedness of the a priori estimate.
For more dicusssions on the choice of structure conditions, the reader shall refer to  Mocha and Westray \cite{MW2012}.
\end{remm}

Motivated by the above discussions,  we prove  an existence  result given \ref{as2'} and exponential moments integrability. Analogously to Theorem \ref{existence3}, 
we  use a Lipschitz-quadratic regularization 
and take the limit by the monotone stability result in Section \ref{section22}.
The a priori bound for $Y$ obtained in  Lemma \ref{aprioriestimate2}  is also  crucial to the construction of an unbounded solution. 
 \begin{thm}[Existence]
\label{existence4}
 If $(f, g, \xi)$ satisfies {\rm\ref{as2'}} and $e^{\beta A_T} \big(|\xi|+ |\alpha|_T\big)$ has exponential moment of order $\gamma$, i.e.,  
\[
\mathbb{E}\Big[ \exp\Big( \gamma e^{\beta A_T}\big (|\xi|+ |\alpha|_T \big)      \Big)\Big] < + \infty,
\]
then there exists a solution  verifying {\rm(\ref{ey})}.
\end{thm}
\begin{Proof} We introduce the notations used throughout the proof.
Define the process \[
X_t: = \frac{1}{\gamma}\ln \mathbb{E}\Big[ \exp \Big(\gamma e^{\beta A_T}\big( |\xi| + |\alpha|_T\big ) \Big)\Big| \mathcal{F}_t  \Big].
\]
Obviously $X$  is continuous by the  continuity of the filtration. 
For $m, n \in \mathbb{N}^+$, set 
\begin{align*}
\tau_m &:= \inf \big\{ t\geq 0: |\alpha|_t + X_t \geq m \big\} \wedge T,\\
\sigma_n &: = \inf \big\{ t\geq 0: |\alpha|_t \geq n\big\} \wedge T.
\end{align*}
It then follows from the  continuity of $X$ and  $|\alpha|_\cdot$ that $\tau_m$ and $\sigma_n$ increase  stationarily to $T$ as $m, n$ goes to $+\infty$, respectively.
To  apply a double approximation procedure, we  define 
 \begin{align*}
f^{n, k}(t, y, z)&: =\mathbb{I}_{\{t\leq \sigma_n \}} \inf_{y^\prime, z^\prime}\big \{ f^+(t, y^\prime, z^\prime) + n|y-y^\prime|+n|\lambda_t(z-z^\prime)| \big\} \\
&-\mathbb{I}_{\{t\leq \sigma_k \}} \inf_{y^\prime, z^\prime} \big\{ f^-(t, y^\prime, z^\prime) + k|y-y^\prime|+k|\lambda_t(z-z^\prime)| \big\},
\end{align*}
and $
\xi^{n , k}: = \xi^+ \wedge n - \xi^- \wedge k.$   

Before proceeding to the proof we give some useful facts. 
By  Lepeltier and San Martin \cite{LS1997}, 
$f^{n, k}$ is  Lipschitz-continuous in $(y, z)$; as  $k$ goes to $+\infty$,  $f^{n, k}$ converges decreasingly 
 uniformly on compact sets to a limit denoted by $f^{n,\infty}$; as $n$ goes to $+\infty$, $f^{n, \infty}$ converges increasingly uniformly on compact sets to $F$. Moreover, $\big||f^{n, k}(\cdot, 0, 0) |\big|_T$ and $\xi^{n, k}$ are bounded. 
 
Hence, by Corollary \ref{unique1},  there exists a unique solution $(Y^{n,k}, Z^{n, k}\cdot M + N^{n, k})\in\mathscr{B}$ to $(f^{n, k}, g, \xi^{n, k})$; by Theorem \ref{compare1} (comparison theorem), $Y^{n, k}$ is increasing in $n$ and decreasing in $k$.
Analogously to Proposition  \ref{existence3}, we wish to take the limit by  
 Theorem \ref{monotonestability} (monotone stability). 
 
 However, 
 $|f^{n, k}(\cdot, 0, 0)|_T$ and $\xi^{n, k}$ 
are not uniformly bounded in general.
To overcome this difficulty, we use Lemma \ref{aprioriestimate2} (a priori estimate) and work on random interval where	
$Y^{n,k}$ and $|f^{n, k}(\cdot, 0, 0)|_\cdot$ are uniformly bounded. This is the motivation to introduce $X$ and $\tau_m$. To be more precise, the localization procedure is 
as follows.

Note that   $(f^{n, k}, g, \xi^{n, k})$ verifies  \ref{as2'} associated with $(\alpha, \beta, \gamma, \varphi)$.  $Y^{n, k}$ being bounded implies that it is of class $\mathcal{D}$. 
Hence  from Lemma \ref{aprioriestimate2}  we have 
\begin{align}
|Y^{n, k}_t|&  \leq \frac{1}{\gamma}\ln \mathbb{E}\Big[ \exp \Big(\gamma e^{\beta A_{t, T}} |\xi^{n ,k}| 
+\gamma \int_{t}^{T} e^{\beta A_{t, s} } \alpha_{s} \mathbb{I}_{\{ s\leq \sigma_n \wedge \sigma_k\}}dA_s
 \Big)\Big| \mathcal{F}_t  \Big]
\nonumber\\
&
\leq \frac{1}{\gamma} \ln \mathbb{E}\Big[ \exp\Big(\gamma e^{\beta A_{t, T}}|\xi| + 
\gamma \int_t^T e^{\beta A_{t, T}}{\alpha_s dA_s} \Big)\Big| \mathcal{F}_t \Big] \label{x1} \\
&\leq X_t.  \nonumber
\end{align}
In view of the definition of $\tau_m$, we have 
\begin{align}
|Y^{n, k}_{t\wedge \tau_m} |&\leq X_{t\wedge \tau_m} \leq m, \nonumber \\
\big||f^{n, k}(\cdot, 0, 0)|\big|_{\tau_m} &\leq |\mathbb{I}_{[0, \tau_m]}\alpha|_{\tau_m} \leq m.
 \label{bound}
\end{align}
Hence $\big||f^{n, k}(\cdot, 0, 0) |\big|_\cdot$ and $Y^{n, k}$ are uniformly bounded on $[0, \tau_m]$. Secondly, given $(Y^{n, k}, Z^{n, k}\cdot M + N^{n, k})$ which solves $(f^{n, k}, g, \xi^{n, k})$, it is immediate that 
$(Y^{n, k}_{\cdot \wedge \tau_m}, (Z^{n, k}\cdot M + N^{n, k})_{\cdot \wedge \tau_m})$  solves $(\mathbb{I}_{[0, \tau_m ]}(t)f^{n, k}(t, y, z), g, Y_{\tau_m}^{n, k})$. We then use Theorem \ref{monotonestability}  as in Proposition \ref{existence3} to construct a pair $(\widetilde{Y}^{m}, (\widetilde{Z}^m \cdot M + \widetilde{N}^m))$ which solves 
$(f, g, \sup_n\inf_k Y_{\tau_m}^{n, k})$, i.e., 
\begin{align}
\widetilde{Y}^m_t =\sup_n\inf_k Y_{\tau_m}^{n, k} + \int_{t\wedge \tau_m}^{\tau_m} \big(
F(s,\widetilde{Y}^m_s, \widetilde{Z}^m_s ) dA_s + g_s \langle \widetilde{N}^m \rangle_s \big)
 -\int_{t\wedge \tau_m}^{\tau_m}   \big(\widetilde{Z}_s^m dM_s + d\widetilde{N}_s \big). \label{bsdem}
\end{align}
Moreover, 
  $\widetilde{Y}^m$ is the $\mathbb{P}$-a.s. uniform limit of $Y^{n, k}_{\cdot\wedge \tau_m}$ and 
$\widetilde{Z}^m\cdot M + \widetilde{N}^m$ is the $\mathcal{M}^2$-limit of $(Z^{n, k}\cdot M + N^{n, k})_{\cdot \wedge \tau_m}$ as $k, n$ go to $+\infty$. Hence
\begin{align}
\widetilde{Y}^{m+1}_{\cdot\wedge \tau_m}&= \widetilde{Y}^{m}_{\cdot\wedge \tau_m}\  \mathbb{P}\text{-a.s.}, \nonumber\\
  \mathbb{I}_{\{t\leq \tau_m\}}\lambda_t\widetilde{Z}^{m+1}_t &= \lambda_t\widetilde{Z}^m_t\ dA \otimes d\mathbb{P}\text{-a.e},\nonumber
  \\
  \widetilde{N}^{m+1}_{\cdot\wedge \tau_m} &= \widetilde{N}^m_{\cdot\wedge \tau_m}\  \mathbb{P}\text{-a.s.} \label{ae}
\end{align}
Define $(Y, Z, N)$ on $[0, T]$ by
\begin{align*}
Y_t &:= \mathbb{I}_{\{t\leq \tau_m\}} \widetilde{Y}_t^1  + \sum_{m\geq 2} \mathbb{I}_{]\tau_{m-1}, \tau_m]}\widetilde{Y}^m_t,  \nonumber
\\
Z_t &:= \mathbb{I}_{\{t\leq \tau_m\}} \widetilde{Z}_t^1  + \sum_{m\geq 2} \mathbb{I}_{]\tau_{m-1}, \tau_m]}\widetilde{Z}^m_t, \nonumber
\\
N_t &:= \mathbb{I}_{\{t\leq \tau_m\}} \widetilde{N}_t^1  + \sum_{m\geq 2} \mathbb{I}_{]\tau_{m-1}, \tau_m]}\widetilde{N}^m_t.
\end{align*}
By (\ref{ae}), we have
$Y_{\cdot \wedge \tau_m} =\widetilde{Y}^m_{\cdot \wedge \tau_m}$, $\mathbb{I}_{\{t\leq \tau_m\}}Z_t = \mathbb{I}_{\{t\leq \tau_m\}}\widetilde{Z}^m_t$ and 
$N_{\cdot \wedge \tau_m} =\widetilde{N}^m_{\cdot \wedge \tau_m}$. Hence we can rewrite (\ref{bsdem}) as 
\[
Y_{t\wedge \tau_m} = \sup_n \inf_k Y^{n, k}_{\tau_m} +\int_{t\wedge \tau_m}^{\tau_m} \big( f{(s, Y_s, Z_s )}dA_s + g_s d\langle N\rangle_s \big) -\int_{t\wedge \tau_m}^{\tau_m}\big(Z_s dM_s + dN_s\big).
\]
 By sending $m$ to $ +\infty$, we prove that $(Y, Z, N)$ solves $(f, g, \xi)$. By (\ref{x1}),  we have 
\begin{align*}
|Y_t | =  |\sup_n \inf_k Y_t^{n, k} | \leq \frac{1}{\gamma} \ln \mathbb{E}\Big[ \exp\Big(\gamma e^{\beta A_{t, T}}|\xi| + 
\gamma \int_t^T e^{\beta A_{t, s}}{\alpha_s dA_s} \Big)\Big| \mathcal{F}_t \Big].
\end{align*}
\qed
\end{Proof}

Compared to    Mocha  and  Westray \cite{MW2012}, we prove the existence result under rather milder structure conditions. For example, \ref{as2'}(ii) gets rid of linear growth in $y$ and allows $g$ to be   any bounded process, which has been seen repeatedly throughout this chapter. 
Secondly,  in contrast to their work,  $dA_t \ll c_A dt$, where $c_A$ is a positive constant, is not needed. 
Finally, 
they use a regularization procedure through quadratic BSDEs with bounded data. Hence, more demanding structure conditions are imposed to ensure that   the comparison theorem holds. On the contrary, the Lipschitz-quadratic regularization is more direct and essentially merely relies on \ref{as2'}
which is the most general assumption to our knowledge.
 This coincides with Briand and Hu \cite{BH2008} for Brownian framework.

Due to the same reason as in Proposition \ref{existence3}, the existence of a maximal or minimal solution is not available.

\begin{remm} Analogously to Hu and Schweizer \cite{HS2011}, 
one may easily extend the existence  result to infinite-horizon case. In abstract terms,  given exponential moments integrability on $\exp (\beta A_\infty)|\alpha|_\infty$, 
we regularize through Lipschitz-quadratic BSDEs with increasing horizons and null terminal value. Using a localization procedure  and the monotone stability result as in Theorem \ref{existence4}, we obtain a solution  which solves the infinite-horizon BSDE.
\end{remm}

As a result from Lemma \ref{aprioriestimate2}, we derive the estimates for the local martingale part. To save pages we only consider the following extremal case.
\begin{Corollary}[Estimate] \label{ubdestimate} Let   {\rm\ref{as2'}}  hold for $(f, g, \xi)$ and  $e^{\beta A_T} \big(|\xi|+ |\alpha|_T\big)$ has exponential moments of all orders. Then
any solution $(Y, Z, N)$ verifying {\rm(\ref{ey})} satisfies:
 $Y$ has exponential moments of all order and $Z\cdot M + N\in \mathcal{M}^p$ for all $p\geq 1.$ More precisely,  for all $p >1$, 
\begin{align*}
\mathbb{E}\big[e^{p\gamma Y^*} \big] \leq  \Big(\frac{p}{p-1}\Big)^p\mathbb{E}\Big[ \exp\Big( p\gamma e^{\beta A_T}\big(|\xi|+ |\alpha|_T \big)  \Big)\Big],
\end{align*}
and for all $p\geq 1$, 
\[
\mathbb{E}\Big[ \Big(\int_0^T\Big( Z_s^\top d\langle M\rangle_s Z_s + d\langle N\rangle_s \Big)\Big)^{\frac{p}{2}}    \Big]
\leq c \mathbb{E}\Big[ \exp\Big( 4p\gamma e^{\beta A_T}\big(|\xi| + |\alpha|_T\big)   \Big)    \Big], \]
where $c$ is a constant  only depending on $p, \gamma $.
\end{Corollary}
\begin{Proof}
The proof is exactly the same as Corollary 4.2,  Mocha and Westray \cite{MW2012} and hence omitted. \qed
\end{Proof}

Let us turn to the uniqueness result. 
We modify  Mocha and  Westray \cite{MW2012} to allow $g$ to be  any bounded  process rather than merely a constant. A convexity assumption  is imposed so as to use $\theta$-technique which proves to be convenient to treat quadratic terms. 
We start from comparison theorem and then move to uniqueness and stability result.  Similar results can be found in  Briand and  Hu \cite{BH2008} for Brownian setting
or  Da Lio and  Ley \cite{LL2006} from the point of view of PDEs.
To this end, the following structure conditions on $(f, g, \xi)$ are needed. 
\begin{as}\label{as3} There exist $\beta \geq 0, \gamma >0$  and an $\mathbb{R}^+$-valued $\Prog$-measurable process $\alpha$ such that $\mathbb{P}$-a.s.
\begin{enumerate}
\item [(i)] for any $t\in [0, T]$, $(y, z)\longmapsto f(t, y, z)$ is continuous;
\item [(ii)] $f$ is Lipschitz-continuous in $y$, i.e., for any $(t, z)\in [0, T]\times\mathbb{R}^d$, $y, y^\prime \in \mathbb{R}$, 
\[
|f(t, y, z) -f(t, y^\prime, z)| \leq \beta |y-y^\prime|;
\]
\item [(iii)] for any $(t, y)\in[0, T]\times\mathbb{R}$, $z\longmapsto f(t, y, z)$ is convex;
\item [(iv)] for any $(t, y, z)\in [0, T]\times\mathbb{R}\times\mathbb{R}^d$, 
\[|f(t, y, z)| \leq \alpha_t + \beta |y| + \frac{\gamma}{2}|\lambda_tz|^2.\]
\end{enumerate}
\end{as}

We start our proof of comparison theorem by observing that \ref{as3} implies \ref{as2'}. Hence existence  is ensured given suitable integrability.
Likewise, we keep in mind that $\mathbb{P}$-a.s. $ |g_\cdot| \leq \frac{\gamma}{2}$. 

\begin{thm}[Comparison Theorem]  
 \label{compare2}
Let $(Y, Z\cdot M + N)$, $(Y^\prime, Z^\prime\cdot  M + N^\prime)\in \mathcal{S}\times \mathcal{M}$ be solutions to $(f, g, \xi), (f^\prime, g^\prime, \xi^\prime)$, respectively, and $Y^*, (Y^\prime)^*$, $|\alpha|_T$ have exponential moments of all orders. If $\mathbb{P}$-a.s. for any $(t, y, z)\in [0, T]\times\mathbb{R}\times\mathbb{R}^d$, 
  $f(t, y, z)\leq f^\prime(t, y, z)$, $g_t \leq g^\prime_t$, $g_t^\prime \geq 0$, $\xi \leq \xi^\prime$ and $(f, g, \xi)$ verifies {\rm\ref{as3}}, 
then $\mathbb{P}$-a.s. $Y_\cdot\leq Y_\cdot^\prime$.
\end{thm}
\begin{Proof}
We introduce the notations used throughout the proof. 
For any $\theta \in (0, 1)$, define 
\begin{align*}
\delta f_t &:= f(t, Y_t^\prime, Z_t^\prime)-f^\prime(t, Y_t^\prime, Z_t^\prime),\\
 \delta_\theta Y &:= Y- \theta Y^\prime, \\
  \delta Y &:= Y-  Y^\prime,
\end{align*}
and  $\delta_\theta Z$,  $\delta Z$,  $\delta_\theta N, \delta N$, etc. analogously. Moreover, 
define \[
\rho_t : = \mathbb{I}_{\{ \delta_\theta Y_t \neq 0\}} \frac{f(t, Y_t, Z_t)-f(t, \theta Y^\prime_t, Z_t)}{\delta_\theta Y_t} .
\]
By  \ref{as3}(ii), $\rho$ is bounded by $\beta$ for any $\theta \in (0, 1)$. Hence $|\rho |_T \leq \beta \norm{A}$.
By It\^{o}'s formula, 
\begin{align*}
e^{|\rho|_t}\delta_\theta Y_t &= e^{|\rho|_T}\delta_\theta Y_T +\int_t^T e^{|\rho|_s}F^\theta_sdA_s +\int_t^T e^{|\rho|_s} \big(g_s d\langle N \rangle_s -\theta g^\prime_s d\langle N^\prime \rangle_s\big) \nonumber\\
&-\int_t^T e^{|\rho|_s}\big(\delta_\theta Z_s dM_s+ d\delta_\theta N_s\big), 
\end{align*}
where
\begin{align}
F^\theta_s &= f(s, Y_s, Z_s)- \theta f^\prime (s, Y_s^\prime, Z^\prime_s) -\rho_s \delta_\theta Y_s, \nonumber\\
  &=\theta \delta f_s + \big(f(s, Y_s, Z_s) - f(s, Y_s^\prime, Z_s)\big) + \big(f(s, Y_s^\prime, Z_s)- \theta f(s, Y^\prime_s, Z^\prime_s)\big)  -\rho_s \delta_\theta Y_s. \label{ftheta}
\end{align}
We then use \ref{as3}(ii)(iii) to deduce that 
\begin{align*}
f(s, Y_s, Z_s) - f(t, Y_s^\prime, Z_s)  &= f(s, Y_s, Z_s) - f(s, \theta Y_s^\prime, Z_s)  + f(s, \theta Y_s^\prime, Z_s) - f(s, Y_s^\prime, Z_s) \\
&= \rho_s \delta_\theta Y_s + f(t, \theta Y_s^\prime, Z_s) - f(s, Y_s^\prime, Z_s)  \\
&\leq \rho_s \delta_\theta Y_s + (1-\theta)\beta |Y_s^\prime|,\\
f(s, Y_s^\prime, Z_s)- \theta f(s, Y^\prime_s, Z^\prime_s) &= f(s, Y_s^\prime, \theta Z_t^\prime + (1-\theta)\frac{\delta_\theta Z_s}{1-\theta}) -\theta f (t, Y^\prime_s, Z^\prime_s) \\
& \leq (1-\theta) f(s, Y_s^\prime, \frac{\delta_\theta Z_s}{1-\theta})\\
&\leq (1-\theta)\alpha_s + (1-\theta)\beta |Y^\prime_s| + \frac{\gamma}{2(1-\theta)}|\lambda_s \delta_\theta Z_s|^2.
\end{align*}
We also note that $\mathbb{P}$-a.s. $\delta f_s \leq 0$.
Hence plugging these inequalities into (\ref{ftheta}) gives 
\begin{align}
F^\theta_s\leq
(1-\theta)\big( \alpha_s + 2\beta |Y_s^\prime| \big) + \frac{\gamma}{2(1-\theta)}|\lambda_s \delta_\theta Z_s|^2. \label{theta1}
\end{align}
We then perform an exponential transform to eliminate both quadratic terms. Set
\begin{align*}
c&:=  \frac{\gamma e^{\beta \norm{A}}}{1-\theta},\\
 P_t&: = \exp \big(ce^{|\rho|_t} \delta_\theta Y_t\big).
\end{align*}
By It\^{o}'s formula, 
\begin{align*}
P_t =P_T &+ \int_t^T  cP_s e^{|\rho|_s}\Big(
F_s^\theta -\frac{ce^{|\rho|_s}}{2}|\delta_\theta Z_s|^2
\Big)dA_s \\
& + \int_t^T cP_s e^{|\rho|_s}
 \Big(
g_s d\langle N \rangle_s -\theta g^\prime_s d\langle N^\prime \rangle_s - \frac{ce^{|\rho|_s}}{2}d\langle \delta_\theta N \rangle_s \Big)\\
&-\int_t^T cP_s e^{|\rho|_s} \big(\delta_\theta Z_s dM_s + d \delta_\theta N_s\big).
\end{align*}
For notational convenience, we define
\begin{align*}
G_t &:=cP_t e^{|\rho|_t}\Big(
F_t^\theta -\frac{ce^{|\rho|_t}}{2}|Z_t^\theta|^2
\Big),\\
H_t &:= \int_0^t cP_s e^{|\rho|_s}
 \Big(
g_s d\langle N \rangle_s -\theta g^\prime_s d\langle N^\prime \rangle_s - \frac{ce^{|\rho|_s}}{2}d\langle N^\theta \rangle_s \Big).
\end{align*}
By (\ref{theta1}), we have
\[
G_t = cP_te^{|{\rho}|_t} \Big( (1-\theta)\big(\alpha_t + 2\beta |Y^\prime_t| \big) \Big) \leq P_tJ_t, 
\]
where \[
J_t: = \gamma e^{2\beta \norm{A}}\big( \alpha_t + 2\beta |Y_t^\prime|\big).
\]
We claim that   $H$ can also be eliminated. Indeed, 
\begin{align*}
d \langle \delta_\theta N \rangle & = d\langle N \rangle + \theta^2 d\langle N^\prime \rangle -2\theta d \langle N, N^\prime \rangle\\
&\gg d\langle N\rangle + \theta^2 d\langle N^\prime \rangle -\theta  d\langle N\rangle  -\theta d\langle N^\prime\rangle \\
& =(1-\theta)\big(d\langle N \rangle - \theta d\langle N^\prime \rangle\big)\\
& =  (1-\theta)  d \delta_\theta \langle N \rangle.
\end{align*}
We then come back to $H$ and use this inequality to deduce that   
\begin{align*}
g_t d\langle N \rangle_t -\theta g^\prime_t d\langle N^\prime \rangle_t - \frac{ce^{|{\rho}|_t}}{2}d\langle \delta_\theta N \rangle_t &= 
g^+_td\langle N\rangle_t - g^-_t d\langle N\rangle_t -\theta  g_t^\prime  d\langle N^\prime \rangle_t -\frac{ce^{|\rho|_t}}{2}d\langle \delta_\theta N \rangle_t \\
&\ll g^+_t d \delta_\theta \langle N\rangle_t + \theta (g^+_t - g^\prime_t) d\langle N^\prime \rangle_t  -\frac{ce^{|{\rho}|_t}}{2}d\langle \delta_\theta N\rangle_t \\
&\ll g^+_td\delta_\theta\langle N\rangle_t -\frac{\gamma}{2(1-\theta)}d\langle \delta_\theta N \rangle_t \\
&\ll 0,
\end{align*}
 due to $g_\cdot^+ \leq g_\cdot^\prime$ and $g_\cdot\leq \frac{\gamma}{2}$. Hence $dH_t \ll 0.$
To eliminate $G$, we set
$
D_t: = \exp\big( |J|_t \big).
$
By It\^o's formula, 
\begin{align*}
d(D_tP_t) &= D_t \Big( \big(P_tJ_t -G_t \big)dA_t -dH_t + cP_t e^{|{\rho}|_t }\big(\delta_\theta Z_tdM_t +d\delta_\theta N_t \big)  \Big).
\end{align*}
But  previous results show that 
$(P_tJ_t -G_t)dA_t - dH_t \gg 0$. Hence $DP$ is a local submartingale.
 Thanks to the exponential moments integrability on $|\alpha|_T$ and $(Y^\prime)^*$ (and hence $|J|_T$), 
we use a localization procedure and the same arguments in Proposition \ref{dcompare} to deduce that 
\begin{align}
P_t &\leq \mathbb{E}\Big[\exp\Big( \int_t^{T} J_s dA_s \Big)P_{T}\Big| \mathcal{F}_t  \Big]. \label{psub}
\end{align}
We  come back to the definition of $P_T$ and observe that 
\begin{align*}
\delta_\theta \xi &\leq (1-\theta) |\xi|+ \theta\delta \xi\\
 &\leq (1-\theta) |\xi|.
\end{align*}
Hence (\ref{psub}) gives
\begin{align*}
\exp\Big( \frac{\gamma e^{\beta \norm{A} + |{\rho}|_t}}{1-\theta}\delta_\theta Y_t \Big)& \leq 
\mathbb{E}\Big[ \exp\Big( \int_t^T J_s      dA_s      \Big) 
\exp\big( ce^{|{\rho}|_T}\delta_\theta \xi	      \big) \Big| \mathcal{F}_t
\Big]  \\ 
&\leq
\mathbb{E}\Big[ \exp\Big( \int_t^T J_s   dA_s      \Big) 
\exp\big( \gamma e^{2\beta \norm{A}}|\xi|	      \big) \Big| \mathcal{F}_t
\Big].
\end{align*}
Hence
\begin{align*}
\delta_\theta Y_t \leq \frac{1-\theta}{\gamma}\ln \mathbb{E}\Big[\exp \Big(\gamma e^{2\beta \norm{A} }\Big( |\xi| + \int_t^T \big(\alpha_s + 2\beta |Y_s^\prime|\big) dA_s    \Big)            \Big) \Big| \mathcal{F}_t  \Big].
\end{align*}
Therefore we obtain $\mathbb{P}$-a.s.  $Y_t \leq Y_t^\prime$, by sending $\theta$ to $1$. By  the continuity of $Y$ and $Y^\prime$, we also have $\mathbb{P}$-a.s. $Y_\cdot\leq Y_\cdot^\prime$.
\qed
\end{Proof}
As a byproduct, we can prove the existence of a unique solution given \ref{as3}. 
\begin{Corollary}[Uniqueness]
\label{exunique}
If $(f, g, \xi)$ satisfies {\rm\ref{as3}}, $\mathbb{P}$-a.s. $g_\cdot\geq 0$ and $|\xi|$, $|\alpha|_T$ have exponential moments of all orders, then there exists a unique solution $(Y, Z, N)$ to $(f, g, \xi)$ such that $Y^*$ has exponential moments of all order and $(Z\cdot M + N)\in \mathcal{M}^p$ for all $p\geq 1$.
\end{Corollary}
\begin{Proof}
The existence  of a unique solution in the above sense is immediate from  Theorem \ref{existence4} (existence), Theorem \ref{compare2} (comparison theorem) and Corollary \ref{ubdestimate} (estimate).
\qed
\end{Proof}
\begin{remm}
There are spaces  to sharpen the uniqueness.   The convexity in $z$ motivates one to replace \ref{as3}(iv) by 
\[
-\underline{\alpha}_t - \beta |y| - \kappa |\lambda_t z|
\leq f (t, y, z) \leq \overline{\alpha}_t +\beta |y| + \frac{\gamma}{2}|\lambda_t z|^2.
\]
Secondly, in  view of Delbaen et al \cite{DHR2011}, 
  we may prove uniqueness given weaker integrability, by characterizing
the  solution as the value process of a stochastic control problem. 
\end{remm}

It turns out that a stability result also holds given convexity condition. The proof is a modification of Theorem \ref{compare2} (comparison theorem). We set $\mathbb{N}^0:= \mathbb{N}^+ \cup \{ 0 \}.$
\begin{prop}[Stability] 
Let $(f^n, g^n, \xi^n)_{n\in \mathbb{N}^0}$  with $g^n_\cdot \geq 0$ $\mathbb{P}$-a.s. 
satisfy {\rm\ref{as3}} associated with $(\alpha^n, \beta, \gamma, \varphi)$, and $(Y^n, Z^n, N^n)$ be their unique solutions in the sense of Corollary {\rm\ref{exunique}}, respectively.
If $\xi^n - \xi^0 \longrightarrow 0$, 
$\int_0^T |f^n -f^0 | (s, Y_s^0, Z_s^0) dA_s  \longrightarrow 0$ in probability, $\mathbb{P}$-a.s. $g^n_\cdot -g^0_\cdot \longrightarrow 0$  as $n$ goes to $+\infty$ and for each $p > 0$, 
\begin{align}
&\sup_{n\in \mathbb{N}^0} \mathbb{E}\Big[ \exp\Big(   p \big(|\xi^n\big| + |\alpha^n|_T \big)           \Big)\Big] <+\infty,   \label{bdddata}
 \\
& \sup_{n\in\mathbb{N}^0} |g^n_\cdot| \leq  \frac{\gamma}{2}\ \mathbb{P}\text{-a.s.}  \nonumber
\end{align}
Then  for each $p \geq  1$, 
\begin{align*}
&\lim_n \mathbb{E}\Big[ \exp\big( p|Y^n-Y^0|^* \big)\Big]  =1, \\
&\lim_n \mathbb{E}\Big[  \Big(\int_0^T \Big((Z_s^n -Z_s^0)^\top d\langle M\rangle_s (Z_s^n -Z_s^0)  +d\langle N^n- N^0 \rangle_s   \Big) \Big)^{\frac{p}{2}} \Big] = 0.
\end{align*}
\end{prop}
\begin{Proof} By Corollary \ref{ubdestimate} (estimate), for any $p \geq 1$, 
\begin{align}
\sup_{n\in \mathbb{N}^0}\mathbb{E}\Big[ \exp\big(p(Y^n)^*\big) + \Big(\int_0^T \Big( (Z_s^n)^\top d\langle M \rangle_s Z_s^n + d\langle N^n\rangle_s \Big) \Big)^{\frac{p}{2}}        \Big] <+\infty. \label{bddn}
\end{align}
Hence the sequence of random variables
\begin{align*}
\exp \Big(p|Y^n- Y^0|^* \Big)  + \Big(\int_0^T \Big((Z_s^n - Z^0_s)^\top d\langle M \rangle_s (Z_s^n - Z^0_s) + d\langle N^n -N^0 \rangle_s\Big)\Big)^{\frac{p}{2}}
\end{align*}
is uniformly integrable. Due to Vitali  convergence,  it is hence sufficient to prove that 
\begin{align*}
|Y^n - Y|^* + \int_0^T \Big((Z_s^n -Z_s^0)^\top  d\langle M\rangle (Z_s^n -Z_s^0) + d\langle N^n - N\rangle_s\Big) \longrightarrow 0
\end{align*}
in probability as $n$ goes to $+\infty$. 

(i).  We prove $u.c.p$ convergence of $Y^n - Y^0$. To this end we use $\theta$-technique in  the spirit of Theorem \ref{compare2} (comparison theorem). For any $\theta \in (0, 1)$, define 
\begin{align*}
\delta f_t^n &: = f^0(t, Y_t^0, Z_t^0)  - f^n (t, Y_t^0, Z_t^0),\\
\delta g^n &: = g^0- g^n,\\
 \delta_\theta Y^n &: = Y^0- \theta Y^n, 
\end{align*}
 and  $\delta_\theta Z^n, \delta_\theta N^n$, $\delta_\theta \langle N\rangle^n$, etc. analogously.  Further, set
\begin{align*}
\rho_t &: =\mathbb{I}_{\{ Y_t^0 - Y_t^n \neq 0\}} \frac{f^n (t, Y_t^0, Z^n_t) -f^n (t, Y_t^n, Z_t^n)}{Y_t^0 - Y_t^n}, \\
c&: = \frac{\gamma e^{\beta \norm{A}}}{1-\theta}, \\
P_t^n &: = \exp \big(c e^{|\rho|_t} \delta_\theta Y_t^n \big), \\
 J_t^n&: = \gamma e^{2  \beta \norm{A}}\big(\alpha_t^n + 2\beta |Y^0_t|\big), \\
 D_t^n &: = \exp \Big(\int_0^t J_s^n dA_s \Big).
\end{align*} 
Obviously $\rho$ is bounded by $\beta $ due to \ref{as3}(i).
 The $\theta$-difference implies that 
\begin{align}
&f^0(t, Y^0_t, Z^0_t) - \theta f^n (t, Y_t^n, Z_t^n) \nonumber \\
& = \delta f_t^n  + \big( \theta f^n(t, Y_t^0, Z_t^n)  - \theta f^n(t, Y^n_t, Z^n_t) \big) +  \big(f^n (t, Y^0_t, Z^0_t) - \theta f^n(t, Y_t^0, Z_t^n)\big). \label{0n}
\end{align}
By \ref{as3}(i)(ii), 
\begin{align*}
\theta f^n(t, Y_t^0, Z_t^n) - \theta f^n(t, Y^n_t, Z^n_t)
& =\theta \rho_t (Y_t^0 - Y_t^n)\\
& = \rho_t \big(\theta Y_t^0  - Y_t^0 + Y_t^0	 -\theta Y_t^n\big)\\
&\leq 
 (1-\theta)\beta |Y_t^0| + \rho_t \delta_\theta Y^n_t,\\
f^n (t, Y^0_t, Z^0_t) - \theta f^n(t, Y_t^0, Z_t^n) &\leq  (1-\theta)\alpha_t^n + (1-\theta) \beta |Y_t^0| + \frac{\gamma}{2(1-\theta)} |\delta_\theta Z^n_t|^2.
\end{align*}
Hence  (\ref{0n}) gives
\begin{align}
f^0(t, Y^0_t, Z^0_t) - \theta f^n (t, Y_t^n, Z_t^n)  -\rho_t \delta_\theta Y^n_t \leq \delta f_t^n + (1-\theta) \big(\alpha_t^n + 2 \beta |Y_t^0| \big) +\frac{\gamma}{2(1-\theta)}|\delta_\theta Z_t^n|^2. \label{thetaf}
\end{align}
To analyze the quadratic term concerning $N^0$ and $N^n$, we deduce by the same arguments as in Theorem \ref{compare2}  that 
\begin{align}
g^0_t d\langle N^0 \rangle_t -\theta g^n_t d\langle N^n \rangle_t - \frac{ce^{|\rho|_t}}{2}d\langle 
\delta_\theta N \rangle_t
&= \delta g_t^n d\langle N^0\rangle_t + g^n_t d\delta_\theta\langle N \rangle_t^n  -\frac{ce^{|\rho|_t}}{2}d\langle \delta_\theta N^n \rangle_t \nonumber\\
&\ll g^n_t \Big(d\delta_\theta\langle N\rangle^n_t -\frac{1}{1-\theta}d\langle \delta_\theta N^n\rangle_t \Big) +  \delta g^n_t d\langle N^0 \rangle_t \nonumber\\
&\ll \delta g^n_t d\langle N^0 \rangle_t. \label{thetag}
\end{align}
Given (\ref{thetaf}) and (\ref{thetag}), we use an  exponential transform which is analogous to that in Theorem \ref{compare2}. This  gives
\begin{align*}
P_t^n \leq D_t^n P_t^n \leq \mathbb{E}\Big[ D_T^n P_T^n + \frac{\gamma e^{2\beta \norm{A}}}{1-\theta}\int_t^T D_s^n P_s^n\big(| \delta f_s^n| dA_s       
+ | \delta g_s^n |d\langle N^0\rangle_s  \big)
\Big| \mathcal{F}_t\Big].
\end{align*}
Using $\log x \leq x$ and $Y^0 - Y^n \leq (1-\theta)|Y^n| + \delta_\theta Y^n$, we deduce that 
\[
Y_t^0 - Y_t^n \leq (1-\theta) |Y_t^n| + \frac{1-\theta}{\gamma} \mathbb{E}\Big[ D_T^n P_T^n + \frac{\gamma e^{2\beta \norm{A}}}{1-\theta}\int_t^T D_s^n P_s^n \big(|\delta f_s^n| dA_s       
+  |\delta g_s^n| d\langle N^0\rangle_s\big)  
\Big| \mathcal{F}_t\Big].
\]
Set
\begin{align*}
\Lambda^n(\theta) &: = \exp\Big( \frac{\gamma e^{2\beta \norm{A}}}{1-\theta} \big((Y^0)^* + (Y^n)^*\big)     \Big) \geq P_t^n, \\
\Xi^n (\theta) &:=  \exp\Big( \frac{\gamma e^{2\beta \norm{A}}}{1-\theta} \big(|\xi^0 -\theta \xi^n| \vee |\xi^n-\theta \xi^0|\big)     \Big) \geq P_T^n.
\end{align*}
We then have
\[
Y_t^0 - Y_t^n \leq (1-\theta) |Y_t^n| + \frac{1-\theta}{\gamma} \mathbb{E}\Big[ D_T^n\Xi^n(\theta)  + \frac{\gamma e^{2\beta \norm{A}}}{1-\theta}D_T^n\Lambda^n(\theta)\int_t^T \big(| \delta f_s^n| dA_s       
+   |\delta g_s^n| d\langle N^0\rangle_s  \big)
\Big| \mathcal{F}_t\Big].
\]
Now we use \ref{as3}(ii)(iii) to $f^n$ and proceed analogously to Theorem \ref{compare2}. This gives
\[
Y_t^n -Y_t^0 \leq (1-\theta)|Y_t^0| +   \frac{1-\theta}{\gamma} \mathbb{E}\Big[ D_T^n \Xi^n(\theta) + \frac{\gamma e^{2\beta  \norm{A}}}{1-\theta}D_T^n\Lambda^n(\theta)\int_t^T \big(| \delta f_s^n| dA_s       
+    |\delta g_s^n| d\langle N^0\rangle_s \big)  
\Big| \mathcal{F}_t\Big].
\]
Though looking symmetric, the two inequalities come from slightly different treatments for the $\theta$-difference. The two estimates give 
\begin{align*}
|Y_t^n -Y_t^0| \leq \underbrace{(1-\theta)\big(|Y_t^0| +|Y_t^n|\big)}_\text{$X^1_t$}  &+ \underbrace{\frac{1-\theta}{\gamma}\mathbb{E}\Big[ D_T^n \Xi^n(\theta)\Big|\mathcal{F}_t\Big]}_{\text{$X^2_t$}} \\
&+ 
\underbrace{e^{2\beta \norm{A}}\mathbb{E}\Big[D_T^n\Lambda^n(\theta)\int_0^T \big(| \delta f_s^n| dA_s     +  |\delta g_s^n| d\langle N^0\rangle_s \big)  
\Big| \mathcal{F}_t\Big]}_{\text{$X^3_t$}}.
\end{align*}
We then prove $u.c.p$ convergence of $Y^n-Y^0.$ For any $\epsilon >0$, 
\begin{align}
\mathbb{P}\Big( |Y^n -Y^0|^* \geq  \epsilon \Big)\leq \mathbb{P}\Big( (X^1)^* \geq \frac{\epsilon}{3}    \Big) 
+\mathbb{P}\Big( (X^2)^* \geq \frac{\epsilon}{3}    \Big) 
+\mathbb{P}\Big( (X^3)^* \geq \frac{\epsilon}{3}    \Big). \label{converge}
\end{align}
We aim at showing that each term on the right-hand side of (\ref{converge}) converges to $0$ if we send $n$ to $+\infty$ first and then $\theta$ to $1$.
To this end, we give some useful estimates.
By Chebyshev's inequality,  
\begin{align*}
\mathbb{P}\Big( (X^1)^* \geq \frac{\epsilon}{3}    \Big) \leq \frac{3(1-\theta)}{\epsilon}\mathbb{E}\big[ (Y^0)^* + (Y^n)^* \big],
\end{align*}
where  $\mathbb{E}[ (Y^0)^* + (Y^n)^* ]$ is uniformly bounded.
Secondly, 
Doob's inequality yields
\begin{align}
\mathbb{P}\Big( (X^2)^* \geq \frac{\epsilon}{3}    \Big)  \leq  \frac{3(1-\theta)\gamma}{\epsilon}\mathbb{E}\big[D_T^n    \Xi_T^n    \big]. \label{x3}
\end{align}
Moreover, by Vitali convergence,  the right-hand side of (\ref{x3}) satisfies
\begin{align*}
\limsup_n \mathbb{E}\big[D_T^n \Xi_T^n \big] &\leq \sup_n\mathbb{E}\big[(D^n)^2\big]^{\frac{1}{2}}\cdot 
\limsup_n\mathbb{E}\big[ (\Xi^n)^2\big]^{\frac{1}{2}}
\\
&\leq \sup_n\mathbb{E}\big[(D^n)^2\big]^{\frac{1}{2}}\cdot 
\mathbb{E}\Big[ \exp\Big(2 \gamma e^{2\beta \norm{A}}|\xi^0|\Big)\Big]^{\frac{1}{2}}\\
 & < +\infty.
\end{align*}
Hence, the first term and the second term on the right-hand side of (\ref{converge}) converge to $0$ as $n$ goes to $+\infty$ and $\theta$  goes to $1$.
Finally, we claim that the third term on the right-hand side of (\ref{converge}) also converges. Indeed,  Doob's inequality and H\"{o}lder's inequality give
\begin{align}
\mathbb{P}\Big( (X^3)^* \geq \frac{\epsilon}{3}    \Big) &\leq \frac{3e^{2\beta\norm{A}}}{\epsilon}\mathbb{E}\Big[D_T^n\Lambda^n(\theta)\int_t^T \big(| \delta f_s^n| dA_s     +  |\delta g_s^n| d\langle N^0\rangle_s \big)  
\Big] \nonumber \\
&\leq \frac{3e^{2\beta\norm{A}}}{\epsilon}\mathbb{E}\Big[\big(D_T^n\Lambda^n(\theta) \big)^2 \Big]^{\frac{1}{2}}
\mathbb{E}\Big[\Big(\int_0^T \big(   | \delta f_s^n| dA_s     +  |\delta g_s^n| d\langle N^0\rangle_s \big) \Big)^{2}    \Big]^{\frac{1}{2}}. \label{x3vitali}
\end{align}
Note that 
\[
\int_0^T \big(   | \delta f_s^n| dA_s     +  |\delta g_s^n| d\langle N^0\rangle_s \big)
\leq |\alpha|_T + |\alpha^n|_T 
+ 2\norm{A} (Y^0)^*  + \gamma \langle Z^0\cdot M + N^0\rangle_T.
\]
Hence the left-hand side of this inequality has finite moments of all orders by Corollary \ref{ubdestimate}. 
Therefore, the left-hand side of (\ref{x3vitali}) converges to $0$ as $n$ goes to $+\infty$ due to Vitali convergence. 

Finally, collecting these convergence results for each term in (\ref{converge})  gives the convergence of $Y^n-Y^0$.

(ii).  It remains to prove convergence of the martingale parts. 
By It\^{o}'s formula, 
\begin{align*}
&\mathbb{E}\Big[ \int_0^T \Big( (Z_s^n - Z_s^0)^\top d\langle M \rangle_s (Z_s^n - Z_s^0) + d\langle N^n -N^0\rangle_s \Big)  \Big]\\
&\leq
  \mathbb{E}\big[ \big|\xi^n - \xi^0 \big|^2\big]  + 2\mathbb{E}\Big[ |Y^n-Y^0|^* \int_0^T \big|F^n(s, Y_s^n, Z_s^n) - F^0(s, Y^0_s, Z^0_s)\big| dA_s \Big]
\\
&+ 2\mathbb{E}\Big[ |Y^n-Y^0|^* \Big| \int_0^T \big(g_s^nd\langle N^n \rangle_s - g^0_s d\langle N^0 \rangle_s \big)         \Big|    \Big],
\end{align*}
As before, we conclude by Vitali convergence. 
\qed
\end{Proof}

\section{Change of Measure}
\label{section24}
We show that given exponential moments integrability, the martingale part $Z\cdot M + N$, though not BMO,   defines an equivalent change of measure, i.e., its stochastic exponential is a  strictly positive martingale. 
We don't require convexity which ensures uniqueness. But to derive the estimate for
$\int_0^T f(s, Y_s, Z_s)dA_s$, we use \ref{as2''} where
$f$ is of linear growth in $y$. 
We keep assuming that $\mathbb{P}$-a.s. $|g_\cdot|\leq \frac{\gamma}{2}$. The following result comes from Mocha and Westray \cite{MW2012}. 
\begin{thm}[Change of Measure]
If $(f, g, \xi)$ satisfies {\rm \ref{as2''}} and $\xi$, $|\alpha|_T$ have exponential moments of all orders, then for any solution $(Y, Z, N)$ such that $Y$ has exponential moments of all orders and any
 $|q| > \frac{\gamma}{2}$, $\mathcal{E}\big(q\big(Z\cdot M + N\big)\big)$ is a continuous martingale.
\end{thm}
\begin{Proof}
We start by recalling  Lemma 1.6. and Lemma 1.7., Kazamaki \cite{K1994}: if $\widetilde{M}$ is a martingale such that 
\begin{align}
\sup_{\tau \in \mathcal{T}}\mathbb{E}\Big[
\exp\Big(\eta \widetilde{M}_\tau      + \Big(\frac{1}{2}- \eta\Big)\langle \widetilde{M} \rangle_\tau \Big)
\Big]  < +\infty, \label{criteria}
\end{align} 
for  $\eta \neq 1$, then $\mathcal{E}\big(\eta \widetilde{M}\big)$ is a martingale. Moreover, if (\ref{criteria}) holds for some $\eta^* > 1$ then it holds for any $\eta \in (1, \eta^*)$.

By  Lemma \ref{ubdestimate} (estimate),  $Z\cdot M + N$ is a continuous martingale. 
First of all, we apply the above criterion to $\widetilde{M}: = \tilde{q}(Z\cdot M + N)$
 for some fixed $|\tilde{q}| > \frac{\gamma}{2}$. Define $ \Lambda_t (\eta)$ such that
 \[
 \ln \Lambda_t (\eta) : = 
 \tilde{q}\eta \big( 
(Z\cdot M)_t +N_t\big)  + \tilde{q}^2 \Big(\frac{1}{2} - \eta\Big)\langle Z\cdot M + N \rangle_t. 
 \]
From  the BSDE (\ref{sbsde})
 and \ref{as2''}, we obtain, for any $\tau \in \mathcal{T}$,
\begin{align}
\ln \Lambda_\tau (\eta)
&=
\tilde{q}\eta \Big(
Y_t - Y_0 +\int_0^t \big(f(s, Y_s, Z_s)dA_s + g_s d\langle N\rangle_s
\big) \Big)
+\tilde{q}^2 \Big( \frac{1}{2} -\eta\Big)\langle Z\cdot M  + N\rangle_t \nonumber\\
&\leq
(2 + \beta\norm{A})|\tilde{q}|\eta  Y^* + |\tilde q|\eta |\alpha|_T  + 
|\tilde{q}|\eta \Big( \frac{\gamma}{2} +\frac{|\tilde{q}|}{\eta}\Big( \frac{1}{2}-\eta\Big)   \Big) \langle Z\cdot M + N \rangle_T. \label{lnlambda}
\end{align}
Note that 
\[
 \frac{\gamma}{2} +\frac{|\tilde{q}|}{\eta}\Big( \frac{1}{2}-\eta\Big)  \leq 0  \Longleftrightarrow \eta \geq \frac{|\tilde{q}|}{2|\tilde{q}|-\gamma} =: q_0 \Big(>\frac{1}{2}\Big).
\]
Hence for any $\eta \geq q_0$, (\ref{lnlambda}) gives
\[
\Lambda_\tau (\eta) \leq 
\exp \big( |\tilde{q}|\eta (2 + \beta) Y_* + |\tilde q|\eta |\alpha|_T    \big).
\]
By exponential moments integrability, we have
\[
\sup_{\tau \in \mathcal{T}}
\mathbb{E}\big[ \Lambda_\tau (\eta)  \big] < + \infty.
\]
It then follows from the first statement of the criterion that $\mathcal{E}\big(\tilde{q}\eta (Z\cdot M + N)\big)$ is a martingale for all $\eta \in [q_0, \infty)\backslash\{1\}$. The second statement ensures that it is a martingale for any $\eta >1$. For any $|q| > \frac{\gamma}{2}$, we set 
 $|\tilde{q}|\in (\frac{\gamma}{2},|q|)$, $\eta := \frac{q}{\tilde{q}} > 1$, and apply the result above to conclude that  $\mathcal{E}\big(q \big(Z\cdot M+ N \big)\big)$ is a martingale.
\qed
\end{Proof}

\chapter{Quadratic Semimartingales with Applications to Quadratic BSDEs}
\label{qs}
\section{Preliminaries}
This chapter is a survey of   the monotone stability  result for quadratic BSDEs studied by Barrieu and El Karoui \cite{BEK2013} . Roughly speaking, it  comes from the stability of quadratic semimartingales which are  processes characterizing the solutions to BSDEs. 
 We fix the time horizon $T>0$, and work on a filtered probability space $(\Omega, \mathcal{F}, (\mathcal{F}_t)_{t\in [0, T]}, \mathbb{P})$ satisfying the usual conditions of right-continuity and $\mathbb{P}$-completeness. We also assume that $\mathcal{F}_0$ is the $\mathbb{P}$-completion of the trivial $\sigma$-algebra. Any measurability  will refer to the filtration $(\mathcal{F}_t)_{t\in [0, T]}$. In particular,  $\Prog$ denotes the progressive $\sigma$-algebra on $\Omega \times [0, T]$.
We assume  the filtration  is \emph{continuous}, in the sense that all  local martingales  have $\mathbb{P}$-a.s. continuous sample paths.  As mentioned in the introduction, 
we exclusively study quadratic semimartingales and BSDEs which are  $\mathbb{R}$-valued. 

Now we introduce all required notations for this chapter. 
 $\ll$ stands for the strong order of nondecreasing processes, stating that the difference is nondecreasing.
For any random variable or process $Y$, we say $Y$ has some property if this is true except on a $\mathbb{P}$-null subset of $\Omega$. Hence we omit ``$\mathbb{P}$-a.s.'' in situations without ambiguity. 
For any random variable $X$, define $\norm{X}_\infty$ to be its essential supremum.
For any c\`adl\`ag process $Y$, set $Y_{s, t}: = Y_t -Y_s$ and $Y^* : = \sup_{t\in [0, T]} |Y_t|$
; 
we denote its total variation process by
$|Y|_\cdot$. 
$\mathcal{T}$ stands for the set of all stopping times valued in $[0, T]$ and 
 $\mathcal{S}$ denotes the space of continuous adapted processes. For later use we specify the following spaces under $\mathbb{P}$.
\begin{itemize}
\item $\mathcal{S}^\infty$: the space of bounded process  $Y\in\mathcal{S}$ with $\norm{Y}:=\norm{Y^*}_\infty < +\infty;$
\item $\mathcal{S}^p (p\geq 1)$: the set of $Y\in\mathcal{S}$ with  $Y^*\in \mathbb{L}^p$;
\item $\mathcal{M}$: the space of $\mathbb{R}$-valued  continuous local martingale;
\item $\mathcal{M}^p(p\geq 1)$: the set 
of $M\in\mathcal{M}$ with  
\[
\norm{M}_{\mathcal{M}^p}: =
\big(\mathbb{E}\big[ (M^*)^p \big] \big)^{\frac{1}{p}}< +\infty;
\]
$\mathcal{M}^p$ is  a Banach space.
\end{itemize}
Finally, for any local martingale $M$, we call $\{\sigma_n\}_{{n\in\mathbb{N}^+ }}\subset \mathcal{T}$ a \emph{localizing sequence} if $\sigma_n$ increases stationarily to $T$ as $n$ goes to $+\infty$ and
$M_{\cdot \wedge \sigma_n}$ is a martingale for any $n\in\mathbb{N}^+$.
\section{Quadratic Semimartingales}
\label{qs2}
In this section, we give the notion and characterizations of quadratic semimartingales.

{\bf $\mathcal{Q}(\Lambda, C)$ Semimartingale.}
 Let $Y$ be a {continuous} semimartingale with canonical decomposition $Y_\cdot = Y_0 -  A_\cdot + M_\cdot$, where $A$ is a continuous adapted  process of finite variation and $M$ is a continuous local martingale with quadratic variation $\langle M \rangle$. Moreover, let $\Lambda$ and $C$ 
 be fixed continuous adapted processes of finite variation. 
  We call $Y$ a $\mathcal{Q}(\Lambda, C)$ semimartingale if  structure condition $\mathcal{Q}(\Lambda, C)$ holds
\begin{align*}
d|A|\ll d\Lambda + |Y|dC + \frac{1}{2}d\langle M\rangle. 
\end{align*}

When there is no ambiguity, 
 $Y$ is also called a  $\mathcal{Q}$ semimartingale or quadratic semimartingale.
Obviously, $Y$ is a $\mathcal{Q} (\Lambda, C)$ semimartingale if and only if $-Y$ is a $\mathcal{Q} (\Lambda, C)$  semimartingale.  Throughout our study, $\Lambda$ and $C$ exclusively denote continuous nondecreasing adapted processes in the above form. 
For any optional process ${Y}$, we define \[
D_\cdot^{\Lambda, C} ({Y}): = \Lambda_\cdot + \int_0^\cdot |{Y}_s|dC_s.\] 

We are  about to introduce the optional strong submartingales  and their decomposition which is crucial to  the characterizations of quadratic semimartingales.
We present a general definition which doesn't require the filtration to be continuous.

{\bf Optional Strong Submartingale.}
 An optional process $Y$ is a  strong submartingale if 
 for any $\tau, \sigma \in \mathcal{T}$ with $\tau \leq \sigma$,
$
\mathbb{E}[Y_\sigma |\mathcal{F}_\tau] \geq Y_\tau 
$ 
and $Y_\sigma$ is integrable. 

By Theorem 4, Appendix I, Dellacherie and Meyer \cite{DM2011},
every optional strong submartingale is indistinguishable from  a l\`adl\`ag process. Hence we  assume without loss of generality that all optional strong submartingales are l\`adl\`ag.
We can also define a (local) optional strong submartingale (respectively supermartingale, martingale) in an obvious way. Though not c\`{a}dl\`{a}g in general, a local optional strong submartingale also has a decomposition of Doob-Meyer's type  called \emph{Mertens decomposition}; see Appendix I,  Dellacherie and  Meyer \cite{DM2011}. More precisely, if $Y$ is a local optional strong submartingale, then it admits a unique decomposition $Y_\cdot = Y_0 +A_\cdot + M_\cdot $, where 
$A$ is a nondecreasing predictable process (which is in general only l\`{a}dl\`{a}g)
and
$M$ is a  local martingale.   If $Y$ is c\`{a}dl\`{a}g, Mertens decomposition coincides with Doob-Meyer decomposition.

We then introduce the following process by using the continuity of the filtration.

{\bf $\mathcal{Q}$ Submartingale.}
We call a  semimartingale $Y$ a $\mathcal{Q}$ submartingale if $ Y_\cdot =  Y_0 - A_\cdot + M_\cdot $, where $M$ is 
a continuous local martingale such that $-A + \frac{1}{2}\langle M \rangle$ is a nondecreasing predictable process.

By It\^{o}'s formula,  $Y$ is a $\mathcal{Q}$ submartingale if and only if $e^{Y_\cdot} = e^{Y_0 -A_\cdot +\frac{1}{2}\langle M \rangle_\cdot}  \mathcal{E}(M)_\cdot$ is a local optional submartingale.

For any optional process $Y$, we  define the following optional processes
\begin{align*}
X_\cdot^{\Lambda, C}(Y)  &:=  Y_\cdot + \Lambda_\cdot + \int_0^\cdot |Y_s|dC_s =  Y_\cdot + D_\cdot^{\Lambda, C}(Y) , \\
U_\cdot^{\Lambda, C}(e^Y) &:= e^{Y_\cdot} + \int_0^\cdot e^{Y_s}d\Lambda_s + \int_0^\cdot e^{Y_s}|Y_s|dC_s = e^{Y_\cdot} 
+\int_0^\cdot e^{Y_s}dD_s^{\Lambda, C}(Y).
\end{align*}
If $\Lambda$ and $C$ are fixed, we set $X(Y): = X^{\Lambda, C}(Y)$, $U(e^Y): = U^{\Lambda, C}(e^Y)$ and 
$D := D^{\Lambda, C}(Y) = D^{\Lambda, C}(-Y)$ when there is no ambiguity. This notation also applies to other processes. 

With the above notions and properties, we prove equivalent characterizations of quadratic semimartingales. 
\begin{thm}[Equivalent Characterizations]
 \label{eqv}
$Y$ is  a $\mathcal{Q}(\Lambda, C)$
 semimartingale if and only if both  $X(Y)$ $(respectively \  U(Y)$$)$ and $X (-Y)$ $(respectively \ U(-Y)$$)$ are $\mathcal{Q}$ submartingales $($respectively local optional strong submartingales$)$. 
\end{thm}
\begin{Proof}
(i). $\Longrightarrow$. Suppose $Y$  has canonical
decomposition $Y_\cdot = Y_0 - A_\cdot +M_\cdot$.  Hence, 
\[
X_\cdot(Y) = Y_0-A_\cdot + D_\cdot  + M_\cdot.
\]
By the structure condition $\mathcal{Q}(\Lambda, C)$, 
\[
 -dA \gg -d|A| \gg - dD -\frac{1}{2} d\langle M \rangle.
\]
This implies that 
$-dA+dD + \frac{1}{2}\langle M\rangle \gg 0$. Hence by definition  $X(Y)$ is a $\mathcal{Q}$ submartingale.  $X(-Y)$ being also a $\mathcal{Q}$ submartingale is immediate since $-Y$ is a $\mathcal{Q}({\Lambda, C})$ semimartingale.

(ii). $\Longleftarrow$.
Suppose 
 $X(Y)$ and $X(-Y)$   admit the following  decomposition
\begin{align*}
X_\cdot(Y) &= Y_0 - \overline{A}_\cdot + \overline{M}_\cdot, \\
 X_\cdot(-Y) &= -Y_0 - \underline{A}_\cdot + \underline{M}_\cdot,
\end{align*}
where $-\overline{A} + \frac{1}{2}\langle \overline{M}\rangle$ and $-\underline{A} + \frac{1}{2}\langle \underline{M}\rangle$
are nondecreasing and predictable. Hence $ -\Delta\overline{A}$ and $ -\Delta\underline{A}$ are nonnegative.
 Moreover, the process
\begin{align}
 2D_\cdot = X_\cdot (Y) + X_\cdot (-Y)  =  -\overline{A}_\cdot -\underline{A}_\cdot + \overline{M}_\cdot+  \underline{M}_\cdot \label{xy}
\end{align}
is of finite variation. Therefore  
$\overline{M}  + \underline{M} = 0 $ and  $\langle \overline{M} \rangle = \langle \underline{M}\rangle$. On the other hand,  the continuity of $D$ implies that $\Delta (-\overline{A}-\underline{A}) =0$. 
A combination of this  fact and 
$ -\Delta\overline{A}$, $ -\Delta\underline{A}\geq 0$
 thus shows that $\Delta \overline{A} = 0$ and  $\Delta \underline{A} = 0$. Hence $Y$ is a continuous semimartingale with canonical decomposition 
\[
Y_\cdot = \frac{X_\cdot(Y) - X_\cdot(-Y)}{2} = Y_0 -A_\cdot + \overline{M}_\cdot, 
\]
where
 $-A := \frac{-\overline{A} + \underline{A}}{2}$. It thus remains to show that $A$ satisfies  the structure condition $\mathcal{Q}(\Lambda, C)$.
 From (\ref{xy}) we obtain
\[
2dD +d\langle \overline{M}\rangle =  \Big( \underbrace{- d\overline{A} + \frac{1}{2} d\langle \overline{M}\rangle}_{\gg 0} \Big)  + \Big(\underbrace{- d\underline{A}  +\frac{1}{2}\langle \overline{M}\rangle}_{\gg 0} \Big).
\]
By  Radon-Nikod\'{y}m theorem  
there exists a predictable process  $\alpha$  valued in $[0, 1]$ such that
\begin{align*}
&d \Big(-\overline{A} + \frac{1}{2}\langle \overline{M}\rangle\Big) = \alpha d \big(2D + \langle \overline{M}\rangle\big), \\
&d \Big(-\underline{A} + \frac{1}{2}\langle \overline{M}\rangle\Big) = (1- \alpha )d \big(2D + \langle \overline{M}\rangle\big).
\end{align*} This gives
\[
-dA = (2\alpha -1)d
\big(D + \frac{1}{2}\langle \overline{M} \rangle\big).\]
Hence
\begin{align*}
d|A| \ll dD + \frac{1}{2}d\langle \overline{M}\rangle.
\end{align*}

(iii). It remains to prove the rest statement. Suppose $Y$ is $\mathcal{Q}(\Lambda, C)$ semimartingale, then $U (e^Y)$ is a continuous semimartingale. 
  It\^{o}'s formula applied to $U(e^Y)$ 
  and $X_\cdot(Y) = Y_\cdot + D_\cdot$ imply 
\begin{align*}
dU_\cdot(e^Y) &= 
d e^{Y_\cdot} + e^{Y_\cdot} dD_\cdot \\ 
&=
de^{X_\cdot(Y)-D_\cdot} + e^{Y_\cdot}dD_\cdot  \\
&= 
e^{-D_\cdot}d e^{X_\cdot(Y)} +e^{X_\cdot (Y)}d e^{-D_\cdot} + e^{Y_\cdot}dD_\cdot\\
&= e^{-D}de^{X_\cdot(Y)}.
\end{align*}
Hence, $U(e^Y)$ is a continuous local submartingale by (i). The same arguments also apply to $U(e^{-Y})$. For the converse direction, we show analogously to (ii) that $Y$ is a  continuous semimartingale by Mertens decomposition of $U(e^{Y})$ and $U(e^{-Y})$.  Therefore, $X(Y)$ and $X(-Y)$ are both continuous semimartingales.
Again It\^{o}'s formula used to $e^{X(Y)}$
gives
\[
de^{X_\cdot(Y)} =e^{D_\cdot} dU_\cdot (e^Y).
\] Hence 
 $X(Y)$ is $\mathcal{Q}$ submartingale. The same arguments also apply to $X(-Y)$.
Finally by (ii) we  conclude that 
$Y$ is $\mathcal{Q}(\Lambda, C)$ semimartingale.
\qed
\end{Proof}

For later use, for any  optional process ${Y}$, we define 
\[
\overline{X}_\cdot^{\Lambda, C} ({Y}) :=
e^{C_\cdot} |{Y}_\cdot| + \int_0^\cdot e^{C_s}d\Lambda_s.
\]
Sometimes only the terminal value of this process matters. Hence we use the same notation $X_T^{\Lambda, C}$ to define, 
for any $\mathcal{F}_T$-measurable random variable $\Xi$, 
\[\overline{X}^{\Lambda, C}_T(|\Xi|)
:=
e^{C_T} |\Xi| + \int_0^T e^{C_s}d\Lambda_s.
\]

\begin{prop}
\label{qsub}
If $Y$ is a $\mathcal{Q}(\Lambda, C)$ semimartingale, then $\overline{X}(Y)$ is a continuous $\mathcal{Q}$ submartingale. 
\end{prop}
\begin{Proof}
By
It\^{o}'s formula,
\begin{align*}
d\overline{X}(Y) &= e^C \Big(|Y|dC + d\Lambda  - \sgn(Y)dA  +\sgn(Y) dM + dL      (Y) \Big)\\
&=e^{C}\Big(dD -\sgn(Y)dA + dL(Y)  \Big)
 + e^C \sgn(Y)dM,
\end{align*}
where $L(Y)$ is the local time of $Y$ at $0$. By the structure condition $\mathcal{Q}(\Lambda, C)$, $\overline{X}(Y)$ is a continuous $\mathcal{Q}$ submartingale.
\qed
\end{Proof}

Analogously we deduce that 
\begin{align} e^{C_{u, \cdot}} |Y_\cdot| +\int_u^\cdot e^{C_{u, s}}d\Lambda_s = e^{-C_u} \cdot (\overline{X}_\cdot(Y)
 - \int_0^u e^{C_s}d\Lambda_s)
  \label{usub}
\end{align} is a $\mathcal{Q}$ submartingale on $[u, T]$ starting from $|Y_u|$ if we view  $u\in [0, T]$ as the intial time.

\section{Stability of Quadratic Semimartingales}
\label{qs3}
Let us turn to the main study of this chapter: stability of quadratic semimartingales. To this end,  we give some estimates used later to prove the convergence of quadratic semimartingales, their martingale parts and finite variation parts in suitable spaces. 

Observe that 
 the nonadapted continuous process $\phi^{\Lambda, C} (|Y_T|)$ defined by 
\[
\phi_\cdot: = e^{C_{\cdot, T}}|Y_T| + \int_\cdot^T e^{C_{\cdot, s}}d\Lambda_s 
\]
is a positive decreasing solution to the ODE
\[
d\phi_t = - (d\Lambda_t +|\phi_t|dC_t), \ \phi_0 = \overline{X}_T(|Y_T|), \ 
\phi_T = |Y_T|.
\]
By differentiating $e^{\phi_\cdot}$ we obtain
\begin{align}
e^{\phi_0} = e^{\phi_\cdot} + \underbrace{\int_0^\cdot e^{\phi_s}d\Lambda_s + \int_0^\cdot e^{\phi_s}|\phi_s|dC_s}_{:=A^\phi_\cdot}. \label{ode}
\end{align} 

Let us  make the following   standing assumption for the estimations. 

{\bf Assumption.}  
For a $\mathcal{Q}(\Lambda, C)$ semimartingale $Y$, set 
$\exp \big(\overline{X}^{\Lambda, C}_T(|Y_T|)\big)\in \mathbb{L}^1$ 
 and define 
$\Phi_\cdot^{\Lambda, C}(|Y_T|) : = \mathbb{E} \big[\exp \big(\phi_\cdot^{\Lambda, C}(|Y_T|)\big)\big|\mathcal{F}_\cdot\big]$.

\begin{thm}\label{classd}
Set $(\phi, \Phi) := (\phi^{\Lambda, C}(|Y_T|), \Phi^{\Lambda, C}(|Y_T|))$. 
\begin{enumerate}
\item[{\rm(i)}]
 $\Phi$ is a   continuous positive supermartingale of class $\mathcal{D}$ with canonical decomposition $\Phi_\cdot = - A^\Phi_\cdot + M^\Phi_\cdot$, where $M^{\Phi}$ is a continuous martingale and
$A^\Phi_\cdot = \int_0^\cdot \Phi_s d\Lambda_s +\int_0^\cdot 
\mathbb{E} \big[ \exp({\phi_{s}})|\phi_{s}|\big|\mathcal{F}_s\big]dC_s.$
\item [{\rm(ii)}]
$U_\cdot(\Phi) = \Phi_\cdot + \int_0^\cdot \Phi_s d\Lambda_s + \int_0^\cdot \Phi_s \ln (\Phi_s)dC_s $ is a  continuous  positive supermartingale of class $\mathcal{D}$ with canonical decomposition $U_\cdot (\Phi) = - A_{\cdot}^U + M_\cdot^U$, where $M^U = M^\Phi$ and $A^U_\cdot = \int_0^\cdot \big(\mathbb{E}\big[\exp({\phi_s})|\phi_s|\big|\mathcal{F}_s\big] - \Phi_s \ln (\Phi_s) \big)dC_s.$
\item [{\rm(iii)}] If in addition $ \exp(|Y_\cdot|) \leq \Phi_\cdot$, then the processes $U(e^Y)$ and $U(e^{-Y})$ are continuous submartingales of class $\mathcal{D}$ dominated by $U(\Phi)$.
\end{enumerate}
\end{thm}
\begin{Proof} (i).
For any $\tau, \sigma \in\mathcal{T}$, $\tau \leq \sigma$, $\phi_\cdot$ being decreasing yields
\begin{align*}
\mathbb{E}\big[\exp(\phi_0)\big|\mathcal{F}_\tau\big]   \geq 
\mathbb{E} \big[\exp(\phi_\tau)\big|\mathcal{F}_\tau\big] 
=\Phi_\tau \geq
\mathbb{E}\big[\mathbb{E} \big[\exp (\phi_\sigma)\big|\mathcal{F}_\sigma\big] \big|\mathcal{F}_\tau\big]
= \mathbb{E}\big[\Phi_\sigma \big|\mathcal{F}_\tau\big].
\end{align*}
Hence $\Phi$ is a supermartingale of class $\mathcal{D}$ which is also the optional projection of $\exp({\phi_{\cdot}})$.  Moreover, since $\Lambda$ and $C$ are continuous,  nondecreasing and adapted, the dual predictable projection of $A^\phi$ in (\ref{ode}) is $A^{\Phi}$ with $\mathbb{E}\big[A^{\phi}_{t, T} - A^{\Phi}_{t, T}\big|\mathcal{F}_t\big]=0$.  Hence  $\widetilde{M}: = \mathbb{E}\big[A_T^{\Phi}- A_T^\phi \big|\mathcal{F}_t\big] =  \mathbb{E}\big[A_t^{\phi} - A_t^{\Phi}\big|\mathcal{F}_t\big]$ is a martingale. 
Hence (\ref{ode}) gives
\[
\mathbb{E}\big[\exp({\phi_0})\big|\mathcal{F}_t]=
\Phi_t + \mathbb{E}\big [  A_t^{\phi} \big|\mathcal{F}_t\big]=\Phi_t  + \widetilde{M}_t + A_t^{\Phi}.
\]
 Then (i) is immediate by setting 
 $M^\Phi_\cdot : = \mathbb{E}\big[\exp({\phi_0})\big|\mathcal{F}_\cdot\big] - \widetilde{M}_\cdot.$ 
 The continuity simply comes from 
 the continuity of $A^{\Phi}$ and $M^{\Phi}$.
 
(ii). By (i), $U(\Phi)$ is a  continuous positive semimartingale  with canonical decomposition $U_\cdot(\Phi)=   - A^U_\cdot + M^U_\cdot$, where $M^U = M^\Phi$ and $A^U_\cdot = \int_0^\cdot \big(\mathbb{E} \big[e^{\phi_s}|\phi_s| \big|\mathcal{F}_s \big]-\Phi_s \ln (\Phi_s)\big)dC_s.$
 By   Jensen's inequality, $A^U_\cdot$ is nondecreasing, hence
$U_\cdot(\Phi)$ is a supermartingale.
 The class $\mathcal{D}$ property comes from the fact that $U(\Phi)$ is dominated by $M^\Phi$.

(iii). This directly comes from $e^{|Y_\cdot|} \leq \Phi_\cdot$, (ii)  and  characterizations of $\mathcal{Q}(\Lambda, C)$ semimartingales (Theorem \ref{eqv}).
\qed
\end{Proof}
\begin{remm}
 A sufficient and necessary condition to verify  $e^{|Y_\cdot|} \leq \Phi_\cdot$ in Theorem \ref{classd}(iii) is that $\exp(\overline{X}(Y))$ is of class $\mathcal{D}$. Indeed, if $\exp (\overline{X}(Y))$ is of class $\mathcal{D}$,
$\exp (|Y_\cdot|) \leq \Phi_\cdot$
is immediate from the $\mathcal{Q}$ submartingale property in (\ref{usub}). For the converse direction we assume $\exp(|Y_\cdot|) \leq \Phi_\cdot = \mathbb{E}\big[\exp (\phi_\cdot)\big|\mathcal{F}_\cdot \big].
$
Taking power $e^{C_{\cdot}}$ on both sides, using  Jensen's inequality to the right-hand side and finally 
multiplying both sides by $\exp\big(\int_0^\cdot e^{C_s}d\Lambda_s\big)$ yields
$\exp (\overline{X}_\cdot(Y)) \leq 
\mathbb{E}\big[ \exp (\phi_0) \big| \mathcal{F}_\cdot\big].$
Hence $\exp (\overline{X}(Y))$ is of class $\mathcal{D}$.

In applications $e^{|Y_\cdot|}\leq \Phi_\cdot$ is natural and often satisfied. For example, in BSDE framework, $e^{|Y_\cdot|} \leq \Phi_\cdot$ can be seen as an estimate for the solution process $Y$.

If  $\Xi$ is a  $\mathcal{F}_T$-measurable random variable such that $\exp\big( \overline{X}_T^{\Lambda, C} (|\Xi|)\big) \in \mathbb{L}^1$,  
then Theorem \ref{classd} still holds. Hence it is a common property of $\mathcal{Q}(\Lambda, C)$ semimartingales whose terminal values are bounded by $|\Xi|.$ This fact will be used to prove the stability result.
\end{remm}

Given stronger integrability condition on
$\overline{X}_T^{\Lambda, C}(|Y_T|)$
we can prove a maximal inequality for ${\mathcal{Q}}(\Lambda, C)$ semimartingales. The proof essentially relies on Proposition \ref{qsub} which states that $\overline{X}^{\Lambda, C}(Y)$ is a $\mathcal{Q}$ submartingale dominating $Y$.  To this end we
define $\psi_p = x^p$ for $p\neq 1$ and 
$\psi_1 =x\ln x - x +1$ for $x\in \mathbb{R}^+$.    
\begin{Lemma}[Maximal Inequality] \label{es0}
Let $p\geq 1.$
If $Y$ is a $\mathcal{Q}(\Lambda, C)$ semimartingale such that 
$\psi_p(\overline{X}_T(|Y_T|))\in\mathbb{L}^1$, then 
\begin{enumerate}
\item [{\rm(i)}]
$
\mathbb{E}\big[\exp(pY^*)\big]^{\frac{1}{p}}
$ is dominated by some increasing function of 
\[\mathbb{E}\Big[\psi_p\Big(\exp\big(\overline{X}_T(|Y_T|)\big)   \Big) \Big].\]
\item [{\rm(ii)}] For any $0<q<1$, $
\mathbb{E}\big[\exp(qY^*)\big]
$
is dominated by some increasing function of \[\psi_q \Big(\mathbb{E}\Big[ 
\exp \big(\overline{X}_T(|Y_T|)\big)
\Big]\Big).
\]
\end{enumerate}
\end{Lemma}
\begin{Proof}
The proof is based on various maximal inequalities and omitted here since it is not relevant to our study of the  stability result.  For details the reader shall refer to 
 Proposition 3.4 and Proposition 3.5,  Barrieu and  El Karoui \cite{BEK2013}.
\qed
\end{Proof}

Given the above estimates  we are ready to introduce a stable family of quadratic semimartingales. The stability result consists in proving convergence of quadratic semimartingales, their finite variation parts and martingale parts.

{\bf $ \mathcal{S}_{\mathcal{Q}}(|\Xi|, \Lambda, C)$ Class.}  Let $\Xi$ be a $\mathcal{F}_T$-measurable random variable
with 
$\exp \big(\overline{X}_T(|\Xi|)\big)\in\mathbb{L}^1$.  Define  $\mathcal{S}_{\mathcal{Q}}(|\Xi|, \Lambda, C)$ to be the set of $\mathcal{Q}(\Lambda, C)$ semimartingales $Y$ with $|Y_T| \leq |\Xi|$ such that  
 $\exp (|Y_\cdot|) \leq \Phi_\cdot (|Y_T|)$.
 By the remark after Theorem \ref{classd}, this inequality is equivalent to $\exp \big(\overline{X}_\cdot(Y)\big)$ being of class $\mathcal{D}$.
Define $\mathbf{P} : = \big\{p\in \mathbb{R}^+: \mathbb{E}\big[\exp \big(p \overline{X}_T(|\Xi|)\big)\big]  < + \infty\big\}$ 
and $p^*: = \sup \mathbf{P}$.
   It is obvious that $1\in\mathbf{P}$ and $p^* \geq 1$. 
   
To prepare for the stability result, we  
give some estimates for the finite variation parts and the martingale parts
in the next two lemmas.
\begin{Lemma}[Estimate] \label{es1}
Let $Y \in \mathcal{S}_{\mathcal{Q}}(|\Xi|, \Lambda, C)$ with canonical decomposition $Y_\cdot = Y_0 -A_\cdot + M_\cdot$.  Set $(\overline{X}_T, \Phi) :=(\overline{X}_T^{\Lambda, C}(|\Xi|),  \Phi^{\Lambda, C}(|\Xi|))$.
\begin{enumerate}
\item [{\rm(i)}] For any $\tau \in \mathcal{T}$, 
\[
\frac{1}{2}\mathbb{E}\big[\langle M\rangle_{\tau, T}\big|\mathcal{F}_\tau\big]
\leq \Phi_\tau \mathbb{I}_{\{\tau < T\}}\leq 
\mathbb{E}\big[\exp (\overline{X}_T) \mathbb{I}_{\{\tau < T\}}\big|\mathcal{F}_\tau\big]. 
\]
In particular, \[
\mathbb{E}\big[\langle M \rangle_T\big] \leq 2\mathbb{E}\big[\exp (\overline{X}_T)\big].
\]
\item[{\rm(ii)}]
If $\Phi_\cdot$ is  bounded, then $M$ is a BMO martingale.
\item [{\rm(iii)}]
For any $p \in \mathbf{P}\cap [1, +\infty)$,  $M\in \mathcal{M}^{2p}$ with
\[
\mathbb{E} \big[\langle M\rangle_T^{p}\big] \leq (2p)^p \mathbb{E}\big[\exp(p\overline{X}_T) \big].
\]
\end{enumerate}
\end{Lemma}
\begin{Proof} (i).  
Define $u(x): = e^{x} - x - 1.$ Hence, $u(x) \geq 0, u^\prime(x) \geq 0$ and $u^{\prime\prime}(x)\geq 1$ for $x \geq 0$; $u\in\mathcal{C} ^2(\mathbb{R})$ and  $u^{\prime\prime} - u^\prime =1$. For any $\tau, \sigma \in \mathcal{T}$, 
 It\^{o}'s formula and the structure condition $\mathcal{Q}(\Lambda, C)$ yield
\begin{align*}
\frac{1}{2}\langle M \rangle_{\tau \wedge \sigma, \sigma} 
\leq 
u(|Y_\sigma|) - u(|Y_{\tau \wedge \sigma}|) + \int_{\tau\wedge \sigma}^\sigma u^\prime(|Y_s|)dD_s - \int_{\tau\wedge \sigma}^\sigma u^\prime (|Y_s|)dM_s.
\end{align*}  
 $\exp \big(\overline{X}^{\Lambda, C}_\cdot(Y)\big)$ being of  class $\mathcal{D}$ implies that  $\exp (|Y_\cdot|)$ is of class $\mathcal{D}$. To eliminate local martingale we replace $\sigma$ by its localizing sequence. 
By Fatou's lemma and class $\mathcal{D}$ property of $\exp (|Y_\cdot|)$,
\begin{align*}
\frac{1}{2}\mathbb{E}\big[\langle M \rangle_{\tau, T}\big|\mathcal{F}_\tau\big]
\leq 
\mathbb{E}\Big[u(|Y_T|) - u(|Y_\tau|) + \int_\tau^T u^{\prime}(|Y_s|)dD_s \Big|\mathcal{F}_\tau\Big].
\end{align*}
By $u^\prime(|Y_s|)\leq \exp({|Y_s|}) \leq \Phi_s^{\Lambda, C}(|Y_T|)\leq \Phi_s$,
\[
\int_\tau^T u^\prime (|Y_s|)dD_s \leq
\int_\tau^T \Phi_s d\Lambda_s + \Phi_s \ln(\Phi_s) dC_s.
\]
Since $U_\cdot(\Phi) = \Phi_\cdot + \int_0^\cdot \Phi_s d\Lambda_s + \int_0^\cdot \Phi_s \ln (\Phi_s)dC_s $  is a supermartingale by Theorem \ref{classd}(ii),
\begin{align}
\frac{1}{2}\mathbb{E} \big[\langle M\rangle_{\tau, T}\big|\mathcal{F}_\tau\big]
&\leq 
\mathbb{E}\Big[ u(|Y_T|) -u(|Y_\tau|)+
\int_\tau^T \Phi_s d\Lambda_s + \Phi_s\ln (\Phi_s)dC_s \Big|\mathcal{F}_\tau\Big] \nonumber\\
&\leq \mathbb{E}\big[u(|Y_T|) -\Phi_T- u(|Y_\tau|) +\Phi_\tau \big|\mathcal{F}_\tau \big] \nonumber\\
&\leq  \Phi_\tau \mathbb{I}_{\{\tau <T \}}\label{psitau}\\
&\leq  \mathbb{E} \big[\exp (\overline{X}_T)\mathbb{I}_{\{\tau <T\}}\big|\mathcal{F}_\tau\big],  \nonumber
\end{align}
where the third inequality is due to $u(|Y_T|) \leq \Phi_T$ and $u(|Y_\tau|) \geq 0$. 

(ii). This is immediate from (i). 

(iii). This is immediate from  Garsia-Neveu lemma (see Chapter VI,  Dellacherie and Meyer \cite{DM2011}) applied to (i).
\qed
\end{Proof}
\begin{Lemma}[Estimate] \label{es2}
Let $Y\in \mathcal{S}_{\mathcal{Q}}(|\Xi|, \Lambda, C)$ with canonical decomposition $Y_\cdot = Y_0 - A_\cdot +M_\cdot$.
Set $(\overline{X}_T, \Phi) :=(\overline{X}_T^{\Lambda, C}(|\Xi|),  \Phi^{\Lambda, C}(|\Xi|))$.
\begin{enumerate} 
\item [{\rm(i)}]
  For any $\tau \in \mathcal{T}$, 
\[
\mathbb{E} \big[|A|_{\tau, T}\big|\mathcal{F}_\tau\big]
\leq 2 \mathbb{E}\big[\exp (\overline{X}_T)\mathbb{I}_{\{\tau <T\}}\big|\mathcal{F}_\tau\big].
\]
In particular, 
\[
\mathbb{E}\big[|A|_T\big] \leq 2\mathbb{E}\big[\exp (\overline{X}_T)\big].
\]
\item [{\rm(ii)}]
If $\Phi$ is bounded  by $c_\Phi$, then for any $\tau \in \mathcal{T}$,
\[
\mathbb{E}\big[|A|_{\tau, T} \big|\mathcal{F}_
\tau\big] \leq 2c_\Phi.
\]
\item [{\rm(iii)}]
For  any $p\in\mathbf{P}\cap [1, +\infty)$, 
the total variation of $A$ satisfies 
\[
\mathbb{E} \big[|A|_T^p\big] \leq (2p)^p \mathbb{E}\big[\exp (p \overline{X}_T)\big].
\]

\end{enumerate}
\end{Lemma}
\begin{Proof}(i).
By  the structure condition $\mathcal{Q}(\Lambda, C)$, 
$ \exp(|Y_\cdot|)\leq \Phi_\cdot$, supermartingale property of $U(\Phi)$
and (\ref{psitau}), we have
\begin{align*}
\mathbb{E}\big[|A|_{\tau, T}\big|\mathcal{F}_\tau\big]
&\leq 
\mathbb{E} \Big[\Lambda_{\tau, T} + \int_\tau ^T |Y_s|dC_s \Big|\mathcal{F}_\tau \Big]
+  \frac{1}{2}\mathbb{E}\big[\langle M\rangle_{\tau, T}\big|\mathcal{F}_\tau \big]\\
&\leq 
\mathbb{E} \Big[\int_\tau^T e^{|Y_s|} \big(d\Lambda_s + |Y_s|dC_s\big)\Big|\mathcal{F}_\tau\Big] + \Phi_\tau \mathbb{I}_{\{\tau <T \}}\\
&\leq 
\mathbb{E}\Big[\int_\tau^T \Phi_s \big(d\Lambda_s + \ln (\Phi_s)dC_s\big)\Big|\mathcal{F}_\tau\Big] +\Phi_\tau \mathbb{I}_{\{\tau <T \}} \\
&\leq 
\mathbb{E}\big[\Phi_\tau -\Phi_T|\mathcal{F}_\tau\big] + \Phi_\tau \mathbb{I}_{\{\tau <T \}}\\
&\leq 
2\mathbb{E}\big[\exp (\overline{X}_T)\mathbb{I}_{\{\tau <T\}}\big|\mathcal{F}_\tau\big].
\end{align*}

(ii). This is 
immediate from (i).

(iii). This is immediate from  Garsia-Neveu lemma.
\qed
\end{Proof}

With the estimates for the finite variation parts and the martingale parts, we are ready to prove the stability result. We start by showing that  $\mathcal{S}_{\mathcal{Q}}(|\Xi|, \Lambda, C)$ is stable by 
$\mathbb{P}$-a.s convergence.

{\bf Stability of $\mathcal{S}_{\mathcal{Q}}(|\Xi|, \Lambda, C)$.}  Let $\{Y^n\}_{n\in\mathbb{N}^+}\subset\mathcal{S}(|\Xi|, \Lambda, C)$ and assume $\mathbb{P}$-a.s. $Y^n_\cdot$  converges    to $Y_\cdot$  on $[0, T]$ as $n$ goes to $+\infty$. By  Theorem \ref{classd}(iii), the continuous   submartingales $U(Y^n)$ and $U(-Y^n)$ are dominated by the positive supermartingale $U (\Phi(|\Xi|))$.
Hence, by dominated convergence, we can pass the submartingale property 
  to $U(Y)$ and $U(-Y)$. Clearly $U(Y)$ and $U(-Y)$ are optional since they are limit of continuous submartingales.
  Then  by 
  characterizations of $\mathcal{Q}(\Lambda, C)$  semimartingales (Theorem \ref{eqv}),  $Y$ is a $\mathcal{Q}(\Lambda, C)$ semimartingale.  Moreover, taking the limit also yields $|Y_T| \leq |\Xi|$ and 
  $\exp(|Y_\cdot|)\leq \Phi_\cdot (|Y_T|)$.
 Hence $Y\in\mathcal{S}_{\mathcal{Q}}(|\Xi|, \Lambda, C).$
 In addition, if  the convergence is monotone,   Dini's theorem implies that the convergence is $\mathbb{P}$-a.s. uniform on $[0, T]$. 

Given the  estimates for  $Y^n, A^n$ and $M^n$, the following theorem states that $A^n$ and $M^n$ also converge in suitable spaces. 
\begin{thm} [Stability] \label{stability1}
Let $\{Y^n\}_{n\in \mathbb{N}^+}\subset\mathcal{S}_{\mathcal{Q}}(|\Xi|, \Lambda, C)$ with canonical decomposition  $Y^n_\cdot = Y_0 -A^n_\cdot + M^n_\cdot$. If ${Y^n}$  is Cauchy for $\mathbb{P}$-a.s. uniform convergence on $[0, T]$, then  the limit process $Y$ belongs to $\mathcal{S}_{\mathcal{Q}}(|\Xi|, \Lambda, C)$. Denote its canonical decomposition by $Y_\cdot = Y_0 - A_\cdot + M_\cdot$.
\begin{enumerate}
\item [{\rm (i)}] For any $1 \leq p<2$, $M^n$ converges to $M$ in $\mathcal{M}^p$.
\item [{\rm (ii)}] If $p^* > 1$, then for any $1\leq p < p^*$,
$M^n$ converges to $M$ in $\mathcal{M}^{2p}$.
\item [{\rm (iii)}] In both cases, $A^n$ converges  at least in $\mathcal{S}^1$ to $A$.
\end{enumerate}
\end{thm}
\begin{Proof}
(i). Set $1\leq p < 2$. Define $\delta Y^{m, n}: = Y^m- Y^n$ and 
$\delta A^{m, n}, \delta M^{m, n}$, etc. analogously. 
 For  each $k \in \mathbb{N}^+$, set
$
\tau_k : = \inf \big\{t\geq 0 : \langle \delta M^{m, n}\rangle_t \geq k  \big \} \wedge T.
$ 
By It\^{o}'s formula, 
\[
|\delta Y_0^{m, n}|^2 + \langle \delta M^{m, n}\rangle_{\tau_k} \leq 
|\delta Y_{\tau_k}^{m, n}|^2 + 2\int_0^{\tau_k} |\delta Y_s^{m, n}|d|\delta A^{m, n}|_s +2 \Big|\int_0^{\tau_k} \delta Y_s^{m, n}d\delta M_s^{m, n}\Big|.
\]
By Davis-Burkholder-Gundy inequality and Cauchy-Schwartz inequality,
\begin{align*}
\mathbb{E}\big[\langle \delta M^{m, n}\rangle_{\tau_k}^{\frac{p}{2}}\big]
&\leq 
\mathbb{E}\big[((\delta Y^{m, n})^*)^p\big] +2^{\frac{p}{2}}\mathbb{E}\big[((\delta Y^{m, n})^*)^{\frac{p}{2-p}}\big]^{\frac{2-p}{2}}\mathbb{E}\big[|\delta A^{m, n}|_T\big]^{\frac{p}{2}}\\
&+ \Big(c(p) \mathbb{E}\big[|\delta ((\delta Y^{m, n})^*)^p\big]^{\frac{1}{2}}\Big)\mathbb{E}\big[\langle \delta M^{m, n}\rangle_{\tau_k}^{\frac{p}{2}}\big]^{\frac{1}{2}}\\
&<+\infty.
\end{align*}
Here $c(p)$ denotes the constant from Davis-Burkholder-Gundy inequality which only depends on $p$.
Using $ab \leq \frac{a^2 + b^2}{2}$, we obtain by transferring $\mathbb{E}\big[\langle \delta M^{m, n}\rangle_{\tau_k}^{\frac{p}{2}}\big]$ to the left-hand side and Fatou's lemma that  
\[
\mathbb{E} \big[\langle \delta M^{m, n} \rangle_{T}^{\frac{p}{2}}\big] \leq
(2+ c(p)^2)\mathbb{E}\big[((\delta Y^{m, n})^*)^p\big] + 2^{\frac{p+2}{2}} \mathbb{E}\big[((\delta Y^{m, n})^*)^{\frac{p}{2-p}}\big]^{\frac{2-p}{2}} \mathbb{E}\big[|\delta A^{m, n}|_T\big]^{\frac{p}{2}}.
\]  
By Lemma \ref{es2}(i), $\mathbb{E}[|\delta A^{m, n}|_T]$ is uniformly bounded. Moreover, 
Lemma \ref{es0}(ii) implies by
de la Vall\'{e}e-Poussin criterion that
 for any $r>0$, $((\delta Y^{m, n})^*)^r$ is uniformly integrable.
Hence Vitali convergence implies that $M^n$ is Cauchy in $\mathcal{M}^{p}$. The $\mathcal{M}^p$-limit of $M^n$  coincides with  $M$ since the canonical decomposition of $Y$ is unique.

(ii).
For any any $\tau, \sigma \in \mathcal{T}$,  It\^{o}'s formula yields 
\begin{align*}
|\delta Y_{\tau \wedge \sigma}^{m, n}|^2 + \langle \delta M ^{m, n}  \rangle_{\tau \wedge \sigma, \sigma}   
&= |\delta Y_\sigma^{m, n}|^2 +
2\int_{\tau\wedge \sigma}^\sigma \delta Y_s^{m, n}d \big(\delta A_s^{m, n} - \delta M_s^{m, n}\big) \\
&\leq ((\delta Y^{m, n})^*)^2 + 2\int_0^T |\delta Y_s^{m, n}|d
|\delta A^{m, n}|_s
-2\int_{\tau\wedge \sigma}^\sigma \delta Y_s^{m, n}d\delta M^{m, n}_s.
\end{align*}
To eliminate the local martingale we replace $\sigma$ by its localizing sequence. Then Fatou's lemma yields
\begin{align*}
\mathbb{E}\big[\langle \delta M ^{m, n}  \rangle_{\tau, T}   
 \big| \mathcal{F}_\tau\big]  \leq 
 \mathbb{E}\Big[ \Big(((\delta Y^{m, n})^*)^2
  + 2(\delta Y^{m, n})^*  |\delta A^{m, n}|_T
\Big)\mathbb{I}_{\{\tau < T\}}\Big|\mathcal{F}_\tau \Big].
\end{align*}
For any $p$ such that $ 1 \leq p  < p^*$, we can find $\epsilon >0$ such that 
$1 \leq p < p+\epsilon < p^*$. 
By Garsia-Neveu lemma and H\"{o}lder inequality, the above estimate gives 
\begin{align*}
\mathbb{E}\big[\langle \delta M^{m, n}\rangle_T^{p}\big]
&\leq p^p\mathbb{E}\Big[\Big(((\delta Y^{m, n})^*)^2 +        (\delta Y^{m, n})^*  |\delta A^{m, n}|_T \Big)^{p}\Big] \\
&\leq 
p^p2^{p-1} \Big( \mathbb{E}\big[((\delta Y^{m, n})^*)^{2p}\big]
+ \mathbb{E}\big[((\delta Y^{m, n})^*)^p   |\delta A^{m, n}|_T^p\big]
\Big)\\
&\leq
p^p2^{p-1}\Big(\mathbb{E}\big[((\delta Y^{m, n})^*)^{2p}\big]
+ \mathbb{E} \big[((\delta Y^{m, n})^*)^{\frac{p(p+\epsilon)}{\epsilon}}\big]^{\frac{\epsilon}{p+\epsilon}}
\mathbb{E} \big[|\delta A^{m, n}|_T^{p+\epsilon}\big]^{\frac{p}{p+\epsilon}}
     \Big).
\end{align*}
Since $p+\epsilon < p^*$, 
$\mathbb{E}\big[|\delta A^{m, n}|_T^{p+\epsilon}\big]$ is uniformly bounded due to Lemma \ref{es2}(iii).  Hence Vitali convergence gives the result.

(iii). This is immediate from  (i)(ii).
\qed
\end{Proof}
\begin{remm}
The case where $p=1$ in Theorem \ref{stability1}(i) is also a consequence of  Barlow and  Protter  \cite{BP1990} which proves convergence of  the martingale parts in $\mathcal{M}^1$ for semimartingales.
\end{remm}
\section{Applications to Quadratic BSDEs}
\label{qs4}
Based on the stability result of quadratic semimartingales obtaind in Section \ref{qs3}, we study the corresponding  monotone stability  result for quadratic BSDEs.  Here we continue with the continuous semimartingale setting in Chapter \ref{qsbsde}.

Recall that  the BSDE $(f, g, \xi)$ is written as follows
\begin{align}
Y_\cdot &= Y_0 - \underbrace{\int_0^\cdot {\big((f(s, Y_s, Z_s) dA_s + g_s d\langle 
N \rangle_s}\big)}_{\text{$:= \widetilde{A}_\cdot$}}
 + \underbrace{\int_0^\cdot \big(Z_s dM_s + dN_s\big)}_{\text{$: = \widetilde{M}_\cdot$}}, \ Y_T  = \xi,  \label{bsde44}
\end{align}
where $ Y_\cdot = Y_0 - \widetilde{A}_\cdot + \widetilde{M}_\cdot$   is the canonical decomposition. Without loss of generality we assume $\mathbb{P}$-a.s. $|g_\cdot|\leq \frac{1}{2}$. Let $\alpha$ be an $\mathbb{R}$-valued $\Prog$-measurable process and $\beta \geq 0$.
  If the structure condition 
 \begin{align}
 |f(t, y, z)|& \leq \alpha_t + \beta |y|+ \frac{1}{2}|\lambda_t z|^2, 
 \label{qssc}
 \end{align}
 holds, 
then 
\begin{align*}
|f(t, Y_t, Z_t)| dA_t &\ll 
\Big(\alpha_t + \beta |Y_t|  + \frac{1}{2}|\lambda_t Z_t|^2 \Big)dA_t\\
& = \Lambda_t + |Y_t|dC_t + \frac{1}{2}d\langle Z\cdot M \rangle_t, 
\end{align*}
where
 $\Lambda := \alpha \cdot A $, $C: = \beta A$. Hence
\[
d|\widetilde{A}| \ll  \Lambda + |Y|dC + \frac{1}{2}d\langle \widetilde{M} \rangle.
\]
Thus if $(Y, Z, N)$ is a solution to (\ref{bsde44}) which satisfies (\ref{qssc}), then $Y$ is a $\mathcal{Q}(\Lambda, C)$ semimartingale. This motivates us to convert
  the machinery of quadratic semimartingales into a monotone stability  result for quadratic BSDEs.

\begin{prop}[Monotone Stability]
\label{qsstability}
Let $(Y^n, Z^n, N^n)$ be solutions to $(f^n, g^n, \xi^n)$ for each $n\in\mathbb{N}^+$, respectively, and 
$Y^n$ be a monotone sequence in $\mathcal{S}_{\mathcal{Q}}(|\Xi|, \Lambda, C)$ which converges $\mathbb{P}$-a.s. to $Y$.  Denote their canonical decomposition  by 
$Y^n_\cdot = Y_0^n - \widetilde{A}^n_\cdot + \widetilde{M}^n_\cdot$, where $\widetilde{M}^n = Z^n\cdot  M + N^n.$
\begin{enumerate}
\item [{\rm(i)}] Then $Y\in\mathcal{S}_\mathcal{Q}(|\Xi|, \Lambda, C)$ 
and the convergence is $\mathbb{P}$-a.s. uniform on $[0, T]$. Denote its canonical decomposition by $Y_\cdot: = Y_0 - \widetilde{A}_\cdot + \widetilde{M}_\cdot$
Then  
$\widetilde{M}^n$ converges in $\mathcal{M}^p$ to $\widetilde{M}$ for any $1\leq p<2$. Moreover, $\widetilde{M}$ admits a decomposition  $\widetilde{M} = Z\cdot M +N.$
\item [{\rm(ii)}]
If $(f^n, g^n, \xi^n)$ satisfies {\rm(\ref{qssc})}
and  $\mathbb{P}$-a.s. for any $t\in [0, T]$, $y^n \longrightarrow y$,  $z^n \longrightarrow z$
$f^n(t, y_n, z_n) \longrightarrow f(t, y, z)$ and  
$g^n_t \longrightarrow g_t$, then $(Y, Z, N)$ solves $(f, g, \lim_n \xi^n)$.
\end{enumerate}
\end{prop}
\begin{Proof}
(i). This is immediate from the stability  result of $\mathcal{S}_{\mathcal{Q}}(\Lambda, C, |\Xi|)$
and 
Theorem \ref{stability1}(i).

(ii). Given the convergence of $Y^n$ and $\widetilde{M}^n$,  it remains 
 to prove 
$\widetilde{A}^n$ converges $u.c.p$ to $\widetilde{A}$, which consists of proving 
\begin{align}
&\int_0^\cdot g^n_s d\langle N^n\rangle_s \longrightarrow \int_0^\cdot g_s d\langle N\rangle_s, \ 
\int_0^\cdot f^n(s, Y_s^n, Z_s^n)dA_s 
\longrightarrow 
\int_0^\cdot f(s, Y_s, Z_s)dA_s  \label{ucp}
\end{align}
$u.c.p$. as $n$ goes to $+\infty$.
To prove the first convergence result, Kunita-Watanabe inequality and Cauchy-Schwartz inequality yield
\begin{align*}
&\mathbb{E}\Big[\Big(\Big| \int_0^\cdot \big(g_s^n d\langle N^n\rangle_s - g_s d\langle N\rangle_s \big)  \Big|^*\Big)^{\frac{1}{2}} \Big]
\\ 
& \leq \mathbb{E}\Big[\Big| \int_0^T |g_s^n| d|\langle N^n\rangle -\langle N\rangle|_s  \Big|^{\frac{1}{2}} \Big]
+ 
\mathbb{E}\Big[\Big( \int_0^T 
|g_s^n - g_s|d\langle N\rangle_s
 \Big)^{\frac{1}{2}}  \Big]\\
   &\leq 
  \frac{1}{2} 
   \mathbb{E}\big[\langle N^n - N \rangle_T^{\frac{1}{2}}\big]^{\frac{1}{2}} \mathbb{E}\big[\langle N^n + N \rangle_T^{\frac{1}{2}}\big]^{\frac{1}{2}}
   + 
   \mathbb{E}\Big[\Big( \int_0^T 
|g_s^n - g_s|d\langle N\rangle_s
   \Big)^{\frac{1}{2}}\Big].
\end{align*}
By Lemma \ref{es1}(i),  $\mathbb{E}\big[\langle N^n + N \rangle_T^{\frac{1}{2}}\big]^{\frac{1}{2}}$ is uniformly bounded. By Theorem \ref{stability1}(i), $N^n$ converges to $N$ in $\mathcal{M}^1$. For the second term, we use dominated convergence. 

To prove the second convergence result in (\ref{ucp}), we use a localization procedure.  Note that
$\Psi_\cdot : = \mathbb{E}\big[\exp \big(\phi_0(\Xi)\big)\big|\mathcal{F}_\cdot \big]$ is a continuous martingale due to the continuity of the filtration. For each $n\in\mathbb{N}^+$, define
$
\sigma_k: = \inf\big\{ t\geq 0: \Psi_t \geq k\big\}.
$
The definition of $\mathcal{S}_{\mathcal{Q}}(|\Xi|, \Lambda, C)$ then implies that 
$\exp (Y^n_{\cdot \wedge \sigma_k}) \leq \Phi_{\cdot\wedge\sigma_k} \leq \Psi_{\cdot \wedge\sigma_k} \leq k$. Secondly, by Lemma \ref{es1}(i), $\widetilde{M}^n$ and hence $Z^n\cdot M$ are uniformly bounded in $\mathcal{M}^2$. Moreover, since $\widetilde{M}^n$ converges to $\widetilde{M}$ in $\mathcal{M}^1$ by Theorem \ref{stability1}(i), we can assume $Z^n$ converges to $Z$ $d\langle M\rangle \otimes d\mathbb{P}$-a.e. and
\[\mathbb{E}\Big[
\int_0^T \sup_n |Z^n_s|^2d\langle M\rangle_s \Big]<+\infty,
\]
by substracting a subsequence; see Lemma 2.5,  Kobylanski \cite{K2000}. By dominated convergence, 
$
\int_0^{\sigma_k} |f^n(s, Y^n_s, Z_s^n)-f(s, Y_s, Z_s)|dA_s
$
converges to $0$ as $n$ goes to $+\infty.$  Hence the second convergence result in (\ref{ucp}) is immediate.
\qed
\end{Proof}

 The  stability result in this section gives a forward point of view to answer the question of convergence.
In contrast to Kobylanski \cite{K2000}, it allows unboundedness  and proves the  stability of $\mathcal{S}_{\mathcal{Q}}(\Lambda, C, |\Xi|)$
 which is later used to show
$\mathcal{M}^p (p\geq 1)$ convergence 
of the martingale parts. 
Nevertheless the structure condition $\mathcal{Q}(\Lambda, C)$ requires a  linear growth in $y$ which
is crucial to the estimate for the finite variation parts; see Lemma \ref{es2}. Hence, given general growth conditions (see, e.g.,  Briand and Hu \cite{BH2008} or Section 
\ref{section22}, \ref{section23})  where the estimate for $A$ is not available, it is more difficult to derive the stability result with the help of quadratic semimartingales. 

To end this section we give an existence example where boundedness  as required by classic existence results is no longer needed. 

{\bf Existence: an Example.} 
 Let the BSDE $(f, 0, \xi)$ satisfy  (\ref{qssc}) and $\exp\big(\overline{X}^{\Lambda, C}_T(|\xi|)\big)\in \mathbb{L}^1$.
We show that there exists a solution to $(f, 0, \xi)$ by  Proposition \ref{qsstability} (monotone stability). For each $n, k\in\mathbb{N}^+$, define
\begin{align*}
f^{n, k}(t, y, z): &= \inf_{y^\prime, z^\prime} \big\{f^+ (t, y^\prime, z^\prime) + n|y-y^\prime| + n|\lambda_t(z-z^\prime)|\big \}\\
&- \inf_{y^\prime, z^\prime}\big\{ (f^- (t, y^\prime, z^\prime) + k|y-y^\prime| + k|\lambda_t(z-z^\prime)|\big\}.
\end{align*}
By Lepeltier and  San Martin \cite{LS1997}, $f^{n, k}$  satisfies  (\ref{qssc})
and 
is Lipschitz-continuous in $(y, z)$. 
Moreover, $\exp\big(\overline{X}^{\Lambda, C}_T(|\xi|)\big)\in \mathbb{L}^1$ implies $\xi, |\alpha|_T \in \mathbb{L}^2$. 
Hence, by  El Karoui and  Huang \cite{EK1997}, there exists a unique solution $(Y^{n, k}, Z^{n, k}, N^{n, k})$ to $(f^{n, k}, 0, \xi)$. To prove $Y^{n, k} \in \mathcal{S}_{\mathcal{Q}}(\Lambda, C, |\xi|)$, 
it remains to show $\exp(|Y^{n, k}_\cdot|)\leq \Phi_\cdot(|\xi^{n, k}|)$.  First of all we assume $\overline{X}_T^{\Lambda, C}(|\xi|)$ is bounded. Then $Y^{n, k}$ is bounded and this inequality holds due to the class $\mathcal{D}$ property of $\overline{X}^{\Lambda, C}(Y^{n, k})$.  Note that the above inequality is stable when  taking the limit in $\overline{X}_T^{\Lambda, C}(|\xi|)$, hence the inequality also holds for  $Y^{n, k}$ with 
$\exp\big(\overline{X}^{\Lambda, C}_T(|\xi|)\big)\in \mathbb{L}^1$.

Given $(Y^{n, k})_{n, k\in\mathbb{N}^+}\subset \mathcal{S}_{\mathcal{Q}}(\Lambda, C, |\xi|)$, we are ready to 
construct a solution by a double approximation procedure. 
By comparison theorem,  $Y^{n, k}$ is decreasing in $k$ and increasing in $n$. Now we fix $n$. $\exp(|Y_\cdot^{n, k}|) \leq \Phi_\cdot(|\xi^{n, k}|)\leq \Phi_\cdot(|\xi|)$ implies that  the limit of $Y^{n, k}$
as $k$ goes to $+\infty$ exists.  We then use
 Proposition \ref{qsstability}  to deduce the existence of a solution 
to  $(f^{n, \infty}, 0, \xi)$. 
  We denote it by  $(Y^{n, \infty}, Z^{n, \infty}, N^{n, \infty})$. 
Thanks to the convergence of $Y^{n, k}$ we can pass the comparison property to $Y^{n, \infty}$. By exactly the same arguments as above, we   construct
a solution to $(f, 0, \xi)$ which is the limit of $(Y^{n, \infty}, Z^{n, \infty}, N^{n, \infty})$ as $n$ goes to $+\infty$
in the sense of Proposition \ref{qsstability}(i).

\bibliography{myrefs}

\begin{thebibliography}{10}

\bibitem{BEO2014}
K~Bahlali, M~Eddahbi, and Y~Ouknine.
\newblock Quadratic bsdes with l2--terminal data existence results, krylov's
  estimate and it\^{o}--krylov's formula.
\newblock {\em arXiv preprint arXiv:1402.6596}, 2014.

\bibitem{BP1990}
Martin~T Barlow and Philip Protter.
\newblock {\em On convergence of semimartingales}.
\newblock Springer, 1990.

\bibitem{B2008}
Pauline Barrieu, Nicolas Cazanave, and Nicole El~Karoui.
\newblock Closedness results for bmo semi-martingales and application to
  quadratic bsdes.
\newblock {\em Comptes Rendus Mathematique}, 346(15):881--886, 2008.

\bibitem{BEK2013}
Pauline Barrieu, Nicole El~Karoui, et~al.
\newblock Monotone stability of quadratic semimartingales with applications to
  unbounded general quadratic bsdes.
\newblock {\em The Annals of Probability}, 41(3B):1831--1863, 2013.

\bibitem{BMS2007}
Giuliana Bordigoni, Anis Matoussi, and Martin Schweizer.
\newblock A stochastic control approach to a robust utility maximization
  problem.
\newblock In {\em Stochastic analysis and applications}, pages 125--151.
  Springer, 2007.

\bibitem{B2003}
Ph~Briand, Bernard Delyon, Ying Hu, Etienne Pardoux, and Lucretiu Stoica.
\newblock $\mathbb{L}^p$ solutions of backward stochastic differential
  equations.
\newblock {\em Stochastic Processes and their Applications}, 108(1):109--129,
  2003.

\bibitem{BC2000}
Philippe Briand and Ren{\'e} Carmona.
\newblock Bsdes with polynomial growth generators.
\newblock {\em International Journal of Stochastic Analysis}, 13(3):207--238,
  2000.

\bibitem{BH2006}
Philippe Briand and Ying Hu.
\newblock Bsde with quadratic growth and unbounded terminal value.
\newblock {\em Probability Theory and Related Fields}, 136(4):604--618, 2006.

\bibitem{BH2008}
Philippe Briand and Ying Hu.
\newblock Quadratic bsdes with convex generators and unbounded terminal
  conditions.
\newblock {\em Probability Theory and Related Fields}, 141(3-4):543--567, 2008.

\bibitem{BLS2007}
Philippe Briand, Jean-Pier Lepeltier, Jaime San~Martin, et~al.
\newblock One-dimensional backward stochastic differential equations whose
  coefficient is monotonic in y and non-lipschitz in z.
\newblock {\em Bernoulli}, 13(1):80--91, 2007.

\bibitem{LL2006}
Francesca Da~Lio and Olivier Ley.
\newblock Uniqueness results for second-order bellman--isaacs equations under
  quadratic growth assumptions and applications.
\newblock {\em SIAM journal on control and optimization}, 45(1):74--106, 2006.

\bibitem{DHR2011}
Freddy Delbaen, Ying Hu, Adrien Richou, et~al.
\newblock On the uniqueness of solutions to quadratic bsdes with convex
  generators and unbounded terminal conditions.
\newblock {\em Ann. Inst. Henri Poincar{\'e} Probab. Stat}, 47(2):559--574,
  2011.

\bibitem{DM2011}
Claude Dellacherie and P-A Meyer.
\newblock {\em Probabilities and Potential, C: Potential Theory for Discrete
  and Continuous Semigroups}.
\newblock Elsevier, 2011.

\bibitem{d1977}
Richard~M Dudley et~al.
\newblock Wiener functionals as it{\^o} integrals.
\newblock {\em The Annals of Probability}, 5(1):140--141, 1977.

\bibitem{EK1997}
N~El~Karoui and SJ~Huang.
\newblock A general result of existence and uniqueness of backward stochastic
  differential equations.
\newblock {\em Pitman Research Notes in Mathematics Series}, pages 27--38,
  1997.

\bibitem{E2013}
EH~Essaky and M~Hassani.
\newblock Generalized bsde with 2-reflecting barriers and stochastic quadratic
  growth.
\newblock {\em Journal of Differential Equations}, 254(3):1500--1528, 2013.

\bibitem{FMW2012}
Christoph Frei, Markus Mocha, and Nicholas Westray.
\newblock Bsdes in utility maximization with bmo market price of risk.
\newblock {\em Stochastic Processes and their Applications}, 122(6):2486--2519,
  2012.

\bibitem{HIM2005}
Ying Hu, Peter Imkeller, Matthias M{\"u}ller, et~al.
\newblock Utility maximization in incomplete markets.
\newblock {\em The Annals of Applied Probability}, 15(3):1691--1712, 2005.

\bibitem{HS2011}
Ying Hu and Martin Schweizer.
\newblock Some new bsde results for an infinite-horizon stochastic control
  problem.
\newblock In {\em Advanced mathematical methods for finance}, pages 367--395.
  Springer, 2011.

\bibitem{JS1987}
Jean Jacod and Albert~N Shiryaev.
\newblock {\em Limit theorems for stochastic processes}, volume 1943877.
\newblock Springer Berlin, 1987.

\bibitem{K1994}
Norihiko Kazamaki.
\newblock {\em Continuous exponential martingales and BMO}.
\newblock Springer, 1994.

\bibitem{K2000}
Magdalena Kobylanski.
\newblock Backward stochastic differential equations and partial differential
  equations with quadratic growth.
\newblock {\em Annals of Probability}, pages 558--602, 2000.

\bibitem{LS1997}
Jean-Pierre Lepeltier and Jaime San~Martin.
\newblock Backward stochastic differential equations with continuous
  coefficient.
\newblock {\em Statistics \& Probability Letters}, 32(4):425--430, 1997.

\bibitem{MS2005}
Michael Mania, Martin Schweizer, et~al.
\newblock Dynamic exponential utility indifference valuation.
\newblock {\em The Annals of Applied Probability}, 15(3):2113--2143, 2005.

\bibitem{MW2012}
Markus Mocha and Nicholas Westray.
\newblock Quadratic semimartingale bsdes under an exponential moments
  condition.
\newblock In {\em S{\'e}minaire de Probabilit{\'e}s XLIV}, pages 105--139.
  Springer, 2012.

\bibitem{M2009}
Marie-Am{\'e}lie Morlais.
\newblock Quadratic bsdes driven by a continuous martingale and applications to
  the utility maximization problem.
\newblock {\em Finance and Stochastics}, 13(1):121--150, 2009.

\bibitem{P1999}
{\'E}tienne Pardoux.
\newblock Bsdes, weak convergence and homogenization of semilinear pdes.
\newblock In {\em Nonlinear analysis, differential equations and control},
  pages 503--549. Springer, 1999.

\bibitem{PP1990}
{\'E}tienne Pardoux and Shige Peng.
\newblock Adapted solution of a backward stochastic differential equation.
\newblock {\em Systems \& Control Letters}, 14(1):55--61, 1990.

\bibitem{RY1999}
Daniel Revuz and Marc Yor.
\newblock {\em Continuous martingales and Brownian motion}, volume 293.
\newblock Springer Science \& Business Media, 1999.

\bibitem{SC2002}
Albert~Nikolaevich Shiryaev and Aleksander~Semenovich Cherny.
\newblock Vector stochastic integrals and the fundamental theorems of asset
  pricing.
\newblock {\em Proceedings of the Steklov Institute of
  Mathematics-Interperiodica Translation}, 237:6--49, 2002.

\bibitem{T2008}
Revaz Tevzadze.
\newblock Solvability of backward stochastic differential equations with
  quadratic growth.
\newblock {\em Stochastic processes and their Applications}, 118(3):503--515,
  2008.

\end{thebibliography}
\bibliographystyle{plain}
\end{document}